\documentstyle[titlepage,leqno,12pt]{article}
\newcommand{\thenewsection}{\bf{\thesection}\rule{1em}{0em}}
\newenvironment{newsection}{\begin{flushleft}\rule{0mm}{3ex}\thenewsection
\bf}{\rule{0mm}{3ex}\noindent\end{flushleft}}

\newcommand{\qed}{\hspace*{\fill}\rule{2mm}{2mm}\vspace{5mm}}

\newcommand{\ad}{\mbox{${\rm ad}\,$}}

\renewcommand{\deg}{\mbox{\rm deg}}

\newcommand{\sgn}{\mbox{\rm sgn}}
\newcommand{\wcw}{\wedge\cdots\wedge}

\newcommand{\norm}[1]{\mbox{$\left\|{#1}\right\|\,$}}

\newcommand{\ra}{\rightarrow}
\newcommand{\lrar}{\longrightarrow}
\newcommand{\llar}{\longleftarrow}

\newcommand{\Ar}[1]{\stackrel{#1}{\relbar\joinrel\relbar\joinrel
    \relbar\joinrel\relbar\joinrel\rightarrow}}

\newcommand{\Arr}[1]{\stackrel{#1}{\relbar\joinrel\relbar\joinrel
    \relbar\joinrel\relbar\joinrel\relbar\joinrel\rightarrow}}

\newcommand{\arr}[1]{\stackrel{#1}{\longrightarrow}}
\newcommand{\larr}[1]{\stackrel{#1}{\longleftarrow}}

\newcommand{\ual}[1]{{\phantom{}\Big\uparrow
\makebox[0pt][r]{$\scriptstyle #1\;$}}}

\newcommand{\uar}[1]{{\phantom{}\Big\uparrow
\makebox[0pt][r]{$\scriptstyle #1\;$}}}

\newcommand{\dal}[1]{{\phantom{}\Big\downarrow
\makebox[0pt][r]{$\scriptstyle #1\;$}}}

\newcommand{\dar}[1]{{\phantom{}\Big\downarrow
\makebox[0pt][l]{$\scriptstyle #1\;$}}}

\newcommand{\Buar}[1]{{\phantom{}\Bigg\uparrow
\makebox[0pt][l]{$\;\scriptstyle #1$}}}

\newcommand{\Bdal}[1]{{\phantom{}\Bigg\downarrow
\makebox[0pt][r]{$\scriptstyle #1\;\;\;$}}}

\newcommand{\Bdar}[1]{{\phantom{}\Bigg\downarrow
\makebox[0pt][l]{$\;\scriptstyle #1$}}}

\newcommand{\Lla}[1]{\stackrel{#1}{\leftarrow\joinrel\relbar\joinrel
    \relbar\joinrel\relbar\joinrel\relbar}}

\newcommand{\str}[3]{\mbox{${#1}_{#2},\dots,{#1}_{#3}$}}
\newcommand{\astr}[2]{\mbox{$a_{#1},\dots ,a_{#2}$}}
\newcommand{\vstr}[2]{\mbox{$v_{#1},\dots ,v_{#2}$}}

\newcommand{\Tor}[4]{\mbox{${\rm Tor}^{#1}_{#2}({#3}, {#4})$}} %% Tor
\newcommand{\ind}{\mbox{\rm ind}}
\newcommand{\larray}{\left(\begin{array}{cc}}
\newcommand{\rarray}{\end{array}\right)}

\newcommand{\ig}{{\rm Im}\,}
\newcommand{\coker}{{\rm coker}\,}

\newcommand{\tr}{\mbox{${\rm tr}$}}
\newcommand{\pr}{\partial}

\newcommand{\om}{\omega}
\newcommand{\Om}{\Omega}

\newcommand{\nt}{\natural}

\newcommand{\End}{{\rm End}}
\newcommand{\mod}{{\rm mod}}

\newcommand{\comp}{\mbox{${\scriptscriptstyle \circ}$}}

\renewcommand{\deg}{{\rm deg}}

\newcommand{\beq}{\begin{equation}}
\newcommand{\eeq}{\end{equation}}
\newcommand{\barr}{\begin{array}}
\newcommand{\earr}{\end{array}}
\newcommand{\beqar}{\begin{eqnarray}}
\newcommand{\eeqar}{\end{eqnarray}}

\newtheorem{proposition}{Proposition}
\newtheorem{leftnumbers}{\rule{-0.35em}{0em}}[section]
\newenvironment{themit}[1]{\begin{leftnumbers}{\sc #1.}}{\end{leftnumbers}}

\newcommand{\btheorem}{\begin{themit}{Theorem}\ }
\newcommand{\etheorem}{\end{themit}}
\newcommand{\bproposition}{\begin{themit}{Proposition}\ }
\newcommand{\eproposition}{\end{themit}}
\newcommand{\blemma}{\begin{themit}{Lemma}\ }
\newcommand{\elemma}{\end{themit}}
\newcommand{\bconjecture}{\begin{themit}{Conjecture}\ }
\newcommand{\econjecture}{\end{themit}}
\newcommand{\bproblem}{\begin{themit}{Problem}\ }
\newcommand{\eproblem}{\end{themit}}
\newcommand{\bcorollary}{\begin{themit}{Corollary}\ }
\newcommand{\ecorollary}{\end{themit}}

\newenvironment{remdefexam}[1]{\begin{leftnumbers}{\sc #1.}
 \rm}{\end{leftnumbers}}

\newcommand{\bremark}{\begin{remdefexam}{Remark}\ }
\newcommand{\eremark}{\end{remdefexam}}
\newcommand{\bdefinition}{\begin{remdefexam}{Definition}}
\newcommand{\edefinition}{\end{remdefexam}}
\newcommand{\bexample}{\begin{remdefexam}{Example}\ }
\newcommand{\eexample}{\end{remdefexam}}

\renewcommand\theequation{\thesection.\arabic{leftnumbers}}
\newcommand{\beqn}{\refstepcounter{leftnumbers}
\begin{equation}}
\newcommand{\eeqn}{\end{equation}}

\newtheorem{remit}[proposition]{Remark}

\newtheorem{defineit}[proposition]{Definition}

\newtheorem{examit}[proposition]{Example}

\def\gl{{\frak gl}}

\newcommand{\point}[2]{\mbox{$(#1,\ldots ,#2)$}}

\newcommand{\tp}[2]{\mbox{$#1^{\otimes #2}$}}

\newcommand{\Hom}[2]{\mbox{${\rm Hom}({#1},{#2})$}}

 \newcounter{subeqn}[equation]
\renewcommand{\thesubeqn}{\theequation.\alph{subeqn}}

\newcommand{\calb}{\mbox{$\cal B$}}
\newcommand{\calc}{\mbox{$\cal C$}}
\newcommand{\cald}{\mbox{$\cal D$}}
\newcommand{\cale}{\mbox{$\cal E$}}
\newcommand{\calf}{\mbox{$\cal F$}}

\newcommand{\call}{\mbox{$\cal L$}}
\newcommand{\calm}{\mbox{$\cal M$}}

\newcommand{\bc}{{\bf C}}

\newcommand{\bh}{{\bf H}}
\newcommand{\bn}{{\bf N}}
\newcommand{\br}{{\bf R}}
\newcommand{\bt}{{\bf T}}

\newcommand{\bz}{{\bf Z}}

\newcommand{\frak}{\bf}
\newcommand{\bk}{\mbox{\bf K}}
\newcommand{\bb}{\mbox{\bf B}}
\newcommand{\bq}{\mbox{\bf Q}}
\newcommand{\semiright}{\,\rhd\!<}
\newcommand{\semileft}{\, >\!\lhd}
\newcommand{\sign}{\mbox{\rm sign}}

\newcommand{\emnorm}{||\; ||}
\newcommand{\rnorm}{||\; ||_r}
\newcommand{\qrnorm}[1]{|| #1||_{Q_r}}
\title{An Introduction to $K$-theory and Cyclic Cohomology}
\author{\bf\sf Jacek Brodzki 
  \\
\mbox{}\\
\em  Department of Mathematics,\\
\em University of Exeter\\
\em North Park Road \\
\em Exeter EX4 4QE, U.K.\\
\footnotesize\tt brodzki@maths.ex.ac.uk
}
\date{}
\begin{document}

\maketitle

\eject
\tableofcontents
\vfill\eject

\section{Introduction}
Cyclic cohomology, a major discovery of the eighties, was defined 
by Alain Connes and, independently,  
 by Boris Tsygan. Since its initial construction, various mathematicians contributed 
to the development of this theory, most notably Cuntz, Goodwillie, Karoubi, 
Kassel, Loday, 
Quillen, Wodzicki.
The main impetus for the study of cyclic cohomology was given by 
Connes's noncommutative geometry, in which it 
plays the role similar to de Rham cohomology in the classical differential 
geometry. Noncommutative geometry offers a radically new
approach to geometry, providing new powerful tools  to study
 certain important
spaces that are beyond the  reach of methods of the classical
differential
geometry. Examples of such spaces are: the space of Penrose
tilings of the plane, the space of leaves of a foliation,
spaces of irreducible representations of discrete groups, and
phase spaces arising in quantum field theory.
Additionally, it has already been demonstrated that  noncommutative geometry
is  very useful in addressing  one of the  most fascinating and 
difficult
problems of physics,  the problem  of  unification of fundamental
forces of nature.
Basic ideas of noncommutative
geometry strongly rely on the theory of $C^*$-algebras
and real analysis. For example, the algebra of functions on a topological
space is replaced by a noncommutative $C^*$-algebra. 
Among the many useful properties of cyclic cohomology is its close relation 
to $K$-theory, which I endeavour to explain in these notes.

This is the  extended version of   lectures delivered at the 
Mathematical Institute
of the Adam Mickiewicz University in Pozna\'{n} in April  1995. 
A part of this text is based 
on the graduate lecture course presented by the 
author at the Department of Mathematics of the University of Texas at 
Austin in the Fall Semester of 1991. 
I 
intended to provide the reader with a rapid introduction to 
the basics of $K$-theory and cyclic cohomology, focusing on 
the relations between the two theories and outlining their 
possible applications. Although this exposition is in no way 
complete, I tried to present here some of the more exciting recent 
results. There are now some very good texts on the subjects. 
First, Connes's book \cite{Connes:Book} gives a very full 
and stimulating account of the early successes of noncommutative 
geometry and it is an indispensable reference to the subject. 
Secondly, cyclic homology has been exhaustively described 
in the monograph of Loday \cite{L}. The reader may also wish 
to consult the review articles by the same author. 
Finally, the book by Rosenberg \cite{Rosenberg:book} provides
a very good introduction to the  algebraic $K$-theory, whereas 
Wegge-Olsen's book \cite{wegge} gives an entertaining exposition 
of the basics of the $K$-theory of operator algebras.

\mbox{}

I am very grateful to Professor Micha\l\ Karo\'{n}ski for inviting me to Pozna\'{n},
and for organizing the lectures. 
I wish to express my gratitude to Alain Connes, Joachim Cuntz and Daniel Quillen 
for many discussions and for sharing their insights with me. 
I thank also  Nick Gilbert, Bill Oxbury, Steven Wilson, 
John Parker, Gavin Jones, Tony Scholl,  J\"{u}rg Wildeshaus, Lyndon Woodward, 
Takashi Kimura, Arkady Vaintrob, Jos\'{e} Gracia-Bondia  and 
Peter V\'{a}mos for many interesting conversations. 
I thank especially  Joseph C. V\'{a}rilly for his careful and critical reading of the
early version of my Austin Lecture Notes, and Grzegorz Banaszak 
for reading an earlier draft of these notes and his many comments and useful suggestions.
 Finally, I extend my thanks 
to the almost defunct British Rail for  providing me with 
some interesting working conditions.

\mbox{}

\mbox{}

\begin{flushright}
{\em Exeter-Southampton, April 1996. }
\end{flushright}

\vfill\eject
\section{The $K$-functor in various situations}
\subsection{The Grothendieck group}
We begin with the introduction of a simple yet powerful notion of 
 Grothendieck's $K$-functor. This chapter contains a long list of examples
of various situations in which the $K$-functor appears naturally. 
It will, therefore, be useful to have at our disposal a few 
of its equivalent descriptions. 
\bdefinition
Let $A$ be an abelian semigroup; we assume for simplicity 
that it contains a zero element.  The {\em Grothendieck
group}\index{Grothendieck group} $K(A)$ of $A$ is an
abelian group that  has the following universal property: There is
a canonical semigroup homomorphism $\phi_A:A \ra K(A) $ such that 
for any group $G$ and a semi-group homomorphism $\psi: A \ra G$,
there is a unique
homomorphism $\gamma : K(A) \ra G$ such that $\psi = \gamma\phi_A$.
This means that the following diagram is commutative.
\setlength{\unitlength}{10pt}
\begin{center}
\begin{picture}(10,7)
\put(2,5){$A$}
\put(4,5.5){\vector(1,0){4}}
\put(6,6){$\scriptstyle\phi_A$}
\put(9,5){$K(A)$}
\put(3.5,4.2){\vector(1,-1){2.5}}
\put(3.2,2.5){$\scriptstyle\psi$}
\put(9,4.2){\vector(-1,-1){2.5}}
\put(8,2.5){$\scriptstyle\gamma$}
\put(5.6,0.3){$G$}
\end{picture}
\end{center}
\edefinition
To prove the existence of $K(A)$, and for later reference,  we provide three constructions of 
$K(A)$ (cf. \cite{atiyah}). 
\begin{enumerate}
\item   Let $F(A)$ be the free abelian group generated
by elements of $A$, and let $E(A)$ be the subgroup of $F(A)$ generated by the
elements of the form $a + a' - (a \oplus a ')$, where $ \oplus$ denotes addition in $F(A)$. 
We define $K(A) = F(A)/ E(A)$, with $\phi_A: A \ra K(A)$ being the composition of the inclusion 
$A \ra F(A)$ with the canonical surjection $F(A) \ra K(A)$. 
\item
 Let $\triangle : A\ra A\times A $ be the diagonal
homomorphism of semigroups, and let $K(A)$ be the set of cosets of
$\triangle (A) $ in $A\times A$. {\em A priori} it is only a quotient semigroup, but
it is not difficult to check that  the
interchange of factors in $A\times A$  induces an inverse in $K(A)$
 so that $K(A)$ is in fact 
a group. We then define $\phi_A : A\ra K(A)$ to be the composition
of the map $a\mapsto (a,0) $  with the natural projection $A\times A\ra
K(A)$. 
\item \label{casethree} Let us now consider the following equivalence relation on the 
product $A\times A$.  We put   $(a,b) \sim (a',b')$ when there
exists a  $p\in A$ such that
$$
a + b' + p = a' + b + p
$$
Then by definition $K(A) = A\times A/\sim $. Elements of $K(A)$ will be 
denoted $[(a,b)]$. 
\end{enumerate}

We thus have a way of assigning an abelian group to each abelian 
semigroup. Moreover, one checks that $K$ is in fact  a covariant functor from the category of
abelian semigroups  to the category of abelian groups. This means 
 that the following diagram is commutative, for any homomorphism 
$\gamma: A\ra B$ of semigroups.
\beqn\label{K-universal}
\begin{array}{ccc}
A & \arr{\phi_A} & K(A) \\
\Bdal{\gamma} &&\Bdar{K(\gamma)}\\
B & \arr{\phi_B} & K(B)
\end{array}
\eeqn
If $B$ is a group, the map $\phi_B$ is an isomorphism.

\bexample\label{VectorSpaces}
 Let $\bn $ be the set of natural numbers including zero. Then
$K(\bn)= \bz$. 

To formulate a slightly less trivial variant of this example, 
let us consider the category ${\rm Vec}$ of all finite dimensional 
vector spaces. Since vector spaces of a given dimension are 
all isomorphic, there is a semi-group isomorphism 
$[{\rm Vec}]\ra \bn$. Here the set $[\rm Vec]$ of 
isomorphism classes of vector spaces is a semi-group 
with respect to the operation induced by the direct sum 
of vector spaces. Thus, from the functoriality of the 
$K$-functor, $K([\rm Vec])= \bz$. 
\eexample

\blemma \label{FORMALDIFF}
If $A$ is an abelian semigroup then any element in $K(A)$ can
be written as $[a]- [b]$ for some elements $a, b\in A$ and where
$[a] = [(a,0)]$, in the notation of \ref{casethree}. Moreover,
$$
[a] - [b] = [a'] - [b']
$$
if and only if there exists $c\in A$ such that
$$
a + b' + c =  a' + b + c
$$
\elemma
Proof. According to \ref{casethree}., any element in 
$K(A)$ is the  class of some $(a,b)=(a,0) + (0,b)$. Note 
that 
$$
[(0,b)] = - [(b,0)]
$$
in $K(A)$. \qed

\bcorollary
Let $a,b$ be two elements in $A$. Then $[a] = [b]$ if and only if
there exists a $c\in A$ such that $a\oplus c =  b \oplus c$.
\ecorollary
\subsection{The $K$-theory of spaces}
As a first example of the general construction of the $K$-groups 
we shall discuss the basics of the topological $K$-theory of Atiyah 
and Hirzebruch. Their construction assigns an invariant to any 
compact topological space $X$ by applying the $K$-functor to  
the semigroup of complex vector bundles on $X$. 
We begin our tour with an important result, due to Serre and Swan 
(cf. \cite{atiyah,K:kteoria}).

\subsubsection{The Serre-Swan Theorem}
 \btheorem \label{SWAN} {\rm [Serre-Swan]}
Let $X$ be a compact Hausdorff space. Let $A= C(X)$ be the algebra  of
continuous  functions on $X$. We shall 
denote by  ${\rm Vec}(X)$  the category of complex vector
bundles over $X$, and by ${\bf P}(A)$ the category of finitely generated 
projective modules over $A$. Then the  categories ${\rm Vec}(X)$ and ${\bf P}(A)$ 
are equivalent.
\etheorem
Proof. Let $\Gamma$ be the functor from ${\rm Vec}(X)$ to ${\bf P}(A)$ which
to any vector bundle $E$ associates  the space of its continuous sections $\Gamma(E)$.
It is clear that $\Gamma(E)$ is an $A$-module
with  the $A$ action defined by
$$
(fs)(x) =  f(x)s(x)
$$
for $f\in A$ and $s\in \Gamma(E)$. To check that this functor is 
well  defined, i.e., that $\Gamma(E)$ is in fact a finitely generated 
projective module over $A$, we need to recall the theorem of 
Swan (see, e.g. \cite{K:kteoria}):
\btheorem {\rm [Swan]}
Any vector bundle $E$ over a compact Hausdorff space 
is a direct summand in a trivial bundle. In other words, there 
exists a vector bundle $F$ such that $E\oplus F$ is trivial.
\etheorem
 
We show first that $\Gamma$ maps a trivial rank $n$  bundle over $X$ 
to the free finitely generated $A$-module $A^n$.
Indeed, any continuous section $s$ of $E$
is a continuous map $s\in C(X,\bc ^n) = C(X)^n = A^n$. On the other hand,
any such function is a section of $E$, and so   $E= A^n$. 
If  $E$ is now  an arbitrary vector bundle over the compact space
$X$, then by Swan's theorem it is a direct summand in a trivial
bundle, thus $\Gamma (E)$ is a direct summand in $A^n$, for some
$n$. It follows that $\Gamma(E)$ is finitely generated and projective.

Furthermore, let us consider two trivial bundles $E= X\times \bc^m$ and
$F = X\times \bc^n$.  Any morphism $f$ of two such bundles
is given by  a matrix $M_f(x) = (a_{ij}(x))$, $i = 1, \dots, m$,
$j=1, \dots, n$ with entries in $A$. It is clear that if two bundle
morphisms have the same matrices then they are the same. Moreover,
any such matrix $M$ represents a morphism of trivial bundles. This
shows that the functor $\Gamma $ establishes 
an equivalence between subcategories of trivial bundles and 
free finitely generated $A$-modules.

Before we can treat the general case, we need to introduce the 
notion of the {\em Karoubi envelope} of a category.  
\bdefinition
Let $\calc $ be any category, and let $\calc ' $ be a category whose objects
are pairs $(X,e)$ where $X\in Ob(\calc)$ and $e:X\ra X$ is an idempotent, 
i.e. $e^2 = e$.
Morphisms in the category $\calc' $ are defined by
$$
{\rm Hom} _{\cal C '}((X,e),(X',e')) = \{ f\in {\rm Hom} _{\cal C}(X,X')| fe =
f,
e'f = f\}
$$
\edefinition
The category $\calc '$ associated with $\calc$ is called the Karoubi
envelope of $\calc$.

The reminder of the proof of theorem \ref{SWAN} relies on the
following simple facts.
\blemma
The Karoubi envelope of the category ${\rm Vec}_T(X)$ of {\em trivial} complex
bundles of finite rank on a compact space  $X$ is equivalent to the category 
${\rm Vec}(X)$ of complex vector
bundles of finite rank on $X$.
\elemma
Proof. We need to show that any vector bundle $E$ on $X$ is of the
form $e\tilde{V}$, where $\tilde{V}$ is the trivial complex vector
bundle with fibre $V$ and $e:\tilde{V}\ra \tilde{V}$ is a projection.

If $F$ is a vector bundle and $e$ is a projection on
$F$, then $eF$ and $(1-e)F$ are subbundles of $F$. This follows from
the fact that the rank of a  projection is an upper semicontinuous function
on $X$. Thus in particular, $e\tilde{V}$ is a subbundle of $\tilde{V}$.

To prove the converse, by Swan's theorem there exists
a bundle $E'$ such that $E\oplus E' = \tilde{V}$. But such a direct
sum decomposition of the trivial bundle $\tilde{V}$ defines a projection
$e$ such that $e\tilde{V}= E$. \qed
\blemma
The Karoubi envelope of the category ${\frak F}(A)$ of finitely generated
free modules over $A$ is the category ${\frak P}(A)$ of all finitely generated
projective  $A$-modules.
\elemma
Proof. Let $P$ be a finitely generated
projective $A$-module. Since $P$ is finitely generated,
there is a surjection $\pi: A^n \ra P$ for some $n$. Since $P$ is
projective, $\pi$ has a section, i.e., an $A$-module
map $s:P\ra A^n$ such that $\pi s : P\ra P$ is the identity map.
It is obvious that  the composite $e=s\pi : A^n \ra A^n$ is an idempotent
on $A^n$ such that $eA^n = P$. Thus $P$ is a direct summand in
$A^n$ and the image of  the projection $e$. \qed

\blemma
If $\calc$, $\cald$ are two equivalent categories, then their
Karoubi envelopes $\calc'$, $\cald '$ are also equivalent categories.
\elemma
Proof. The proof is not difficult and is left to the reader. \qed

Theorem \ref{SWAN} is now proved. 
\subsubsection{$K^0(X)$}
\bdefinition
Let $ {\rm Vec }(X)$ be the space of finite rank vector bundles on
a compact Hausdorff space $X$. Let $[{\rm Vec}(X)]$ be the semigroup
of isomorphism classes of vector bundles on $X$ with addition 
induced by the Whitney sum of vector bundles. We define
$K^0(X) = K([{\rm Vec}(X)])$
\edefinition

The following statement is an easy translation of lemma \ref{FORMALDIFF}
to this special case.
\blemma
Every element of $K^0(X)$ can be written as
$[E] - [\theta^n]$ for some $n$ and some vector bundle $E$. Here $[\theta^n]$
denotes the class of a trivial bundle of rank $n$. Moreover,
$$
[E] - [\theta^n] = [F]- [\theta^p]
$$
if and only there exists $q$ such that
\beqn\label{stable}
E\oplus \theta^{p+q} = F\oplus \theta^{n+q}.
\eeqn
\elemma
Proof. Let $[E]- [F]$ be an element of $K^0(X)$, and let $F'$ be
a supplement of $F$ in the trivial bundle of rank $n$, for some
$n$. We have
$$
[E] - [F] = [E] + [F']- ( [F] + [F']) = [E+ F'] - [\theta^n]
$$
The second part of the lemma is now clear. \qed
\bremark
It follows from (\ref{stable}) and (\ref{FORMALDIFF}) that two 
classes $[E]$ and $[F]$ are
equal in $K^0(X)$ if and only if there
exists  $n$
$$
E\oplus \theta^n \simeq F \oplus \theta^n.
$$
Two such  bundles $E$ and $F$ are  called {\em stably equivalent}.
We need to consider stable equivalence instead of just isomorphism
of bundles because there is no cancellation law
in the category of vector bundles.
\eremark
\blemma\label{KPOINT}
If $X=*$ consists of just one point, then $K^0(*)= \bz$.
\elemma
Proof. A vector bundle over a point is just a vector space, and so
this problem reduces to example \ref{VectorSpaces}. \qed
\bremark
This example, although trivial, identifies a class of topological spaces
whose $K^0$-group is trivial. This is thanks to the {\em homotopy 
invariance} of $K^0$, which we shall discuss in the next section. 
We can say here that if a compact space $X$ is contractible to a point, 
then $K^0(X) = K^0(*) = \bz$. This fact may be interpreted geometrically 
by saying that any vector bundle over a contractible compact space 
must be trivial. 
\eremark

Finally, we note that $K^0(X)$ is {\em contravariant} in $X$. This is due to
the fact that any continuous map $f:Y\ra X$ acts by pullback
on vector bundles over $X$ thus inducing a map ${\rm Vec}(X)\ra {\rm Vec}(Y)$
and so a map $K^0(X)\ra K^0(Y)$.

\subsubsection{Reduced $K$-theory}
Let us consider the constant map
$$
p: X\lrar *
$$
 which sends the space $X$ to
a single point $*$. This map induces a map of the corresponding $K$-groups
$$
p^*:K^0(*) \lrar K^0(X)
$$
\bdefinition
The reduced $K$-theory $\tilde{K}^0(X)$ is defined to be the cokernel
of the map $p^*$:
\beqn\label{kdef}
0\lrar \bz \arr{p^*} K^0(X) \arr{} \tilde{K}^0(X) \lrar 0
\eeqn
\edefinition

\bproposition
A choice of the base-point in $X$ defines a canonical splitting of the 
sequence (\ref{kdef})
so that
$$
\tilde{K}^0(X) \simeq \ker [ K^0(X)\lrar K^0(*) \simeq \bz]
$$
and
$$
K^0(X) \simeq \bz \oplus \tilde{K}^0 (X)
$$
\eproposition
Proof. A choice of the base-point in $X$ gives the inclusion  
$i:*\hookrightarrow X$, which  induces a homomorphism 
$i^*:K^0(X)\ra K^0(*)=\bz$ that is a left inverse of the map $p^*$. \qed
\bproposition
If $X$ is a disjoint union  $X= X_1\sqcup\dots \sqcup X_n$ then the
inclusions  of $X_i$ into $X$ induce a decomposition of $K^0(X)$
as a direct product $K^0(X_1) \oplus \dots \oplus K^0(X_n)$.
\eproposition
Proof. 
Any bundle over $X$ is characterized by its restrictions to each $X_i$
so the semigroup  of isomorphism classes of bundles over $X$ is
a product of isomorphism classes of bundles over the components
$X_i$. Hence the result. \qed
\bremark
This proposition is not true for the reduced $\tilde{K}^0$-functor. For
example, if $X$ is the disjoint union of two points $p_1$ and $p_2$ then
$\tilde{K}^0(X)= \bz$ but $\tilde{K}^0(p_i) = 0$ for $i=1,2$.
\eremark
\bdefinition\label{rank-function}
Let $E$ be any finite rank complex vector bundle over 
a topological space $X$, and let us denote by $E_x$ the 
fibre of $E$ over the point $x\in X$. We define the {\em rank function}
$r: X\ra \bn$ by $r(x) = \dim E_x$. Since $E$ is locally  trivial, 
the rank function is also locally constant. 
The space of all locally constant $\bn$-valued
functions on $X$ with pointwise addition is an abelian semigroup, 
which we shall denote by $H^0(X,\bn)$. 
\edefinition
The map $r$ extends in a natural way to a group homomorphism
$$
r: K^0(X) \lrar H^0(X,\bz),
$$
by defining
$$
r([E] - [F]) = r(E) - r(F).
$$
for all $[E]-[F]\in K^0(X)$. The group $H^0(X,\bz)$ 
is the $0$-th \v{C}ech cohomology group of $X$. 
\bproposition
Let $K'(X) = \ker [r:K^0(X) \lrar H^0  (X,\bz )]$. The  short exact sequence
$$
0\ra K'(X) \ra K^0(X) \ra H^0(X,Z) \ra 0
$$
splits canonically. It follows that 
if $X$ is connected, then $K'(X)$ is isomorphic to the reduced 
$K$-group $ \tilde{K}^0(X)$.
\eproposition
Proof. Let $f:X\ra \bn$ be a locally constant function. Since $X$ is
compact, $f$ only takes on a finite number of values $n_1, \dots , n_p$.
Let $X= X_1 \cup \dots \cup X_p $, where $X_i =  f^{-1} (\{ n_i\})$.
Let $E$ be the vector bundle whose restriction to $X_i$ is 
given by  $E|_{X_i} = X_i \times \bc^{n_i}$, $i=1, \dots, p$.
The correspondence $f\mapsto E$ defines a semigroup homomorphism
$$
t: H^0 (X,\bn) \ra [{\rm Vec}(X)]
$$
which induces a map  $H^0(X,\bz)\ra K^0(X)$. This map is a right inverse to 
the map
$r: K^0(X) \ra H^0(X,\bz)$.

When $X$ is connected, $H^0(X,\bz) \simeq \bz$,  and so from the first part 
we have that $K'(X) \simeq \coker[ \bz \ra K(X)]\simeq \tilde{K}^0(X) $.
\qed
\bremark
The group $H^0(X,\bz)$ is the zeroth homology group of $X$ with 
integer coefficients. The rank function  $r$, which takes
values in $H^0(X,\bz)$, is a first example of the Chern 
character.
\eremark

\subsubsection{Homotopy invariance and excision}\label{homotopy-and-excision}
The homotopy invariance and excision are two basic properties of 
the topological $K$-functor which are most useful in computations
of $K$-groups. Also, as they are imitated in various other versions 
of $K$-theory, it is worthwhile to develop a clear geometric picture
for results of this type. The homotopy invariance, which was already 
alluded to in the previous section, describes the behaviour of 
the $K$-functor with respect to continuous deformations of 
the base space $X$.  
The interest in excision, on the other hand, originates from another basic operation 
in algebraic topology. Let $X$ be a compact space and $Y$ a 
closed subspace. We shall define a  compact space $X/Y$, 
in which $Y$ is shrunk to a point, in some sense. We need to extend 
our definition of the $K$-functor so that it behaves well with respect 
to this operation; this will lead to the introduction of the relative 
$K$-theory. This new functor will establish $K$-theory as an \lq exotic'
cohomology theory. For more details the reader may wish to consult 
\cite{atiyah} or \cite{K:kteoria}.

We begin with the homotopy invariance of $K^0$. Let $Y$ be a compact 
Hausdorff space and let $f,g$ be two maps from $Y$ to $X$. We say
that $f$ and $g$ are {\em homotopic} if and only if there 
exists a continuous map $F:Y\times I \ra X$, where $I$ is the unit
interval, such that $F_0=F(\cdot, 0)= f$ and $F_1= F(\cdot , 1) = g$. 
\btheorem\label{K-homotopy} {\em [Homotopy invariance]}
Let $Y$ be a compact Hausdorff space, and let $f, g:Y\ra X$.
Let us assume that  $f$ and $g$ are homotopic maps from $Y$.
For any vector bundle $E$ there is an isomorphism 
$$
f^* E \simeq g^* E
$$
It follows that the maps induced by $f$ and $g$ on $K$-groups
are the same, i.e.:
$$
K^0(f) = K^0(g): K^0(X)\arr{} K^0(Y)
$$
\etheorem 
The proof of this theorem is not difficult and relies on the 
following two results. 
\bproposition{\rm [Tietze's extension theorem]}
Let $X$ be a compact Hausdorff space, $Y$ a closed subspace, 
 and let $E$ be a bundle over $X$.
Any continuous section $s: Y\ra E|_Y$ can be extended to $X$. 
\eproposition
Proof. Take $s\in \Gamma(E|_Y)$. Since, locally, $s$ is a vector-valued
function, we know (from Tietze's theorem for vector valued functions) that
for any $x\in X$ there exists an open neighbourhood $U$ containing
$x$ and $t\in \Gamma(E|_U)$ such that $t$ and $s$ agree on the
overlap $U\cap Y$. Since $X$ is compact, we can  find a finite subcover
$\{U_\alpha\}$ by such sets. Let $t_\alpha \in \Gamma (E|_{U_\alpha})$
be the corresponding section and let $\{p_\alpha\}$  a partition of unity
subordinate to the cover $\{U_\alpha\}$. We  define
$s_\alpha \in \Gamma (E)$ by
$$
s_\alpha = \left\{\begin{array}{cr} p_\alpha(x) t_\alpha(x) ,&  {\rm if} \;
x\in U_\alpha\\
0 & otherwise
\end{array}\right.
$$
Then $\sum_\alpha s_\alpha$ is a section of $E$ whose restriction to $Y$ is
$s$.
\qed
\bproposition
Let $Y$ be a closed subspace of a compact Hausdorff space $X$ and let
$E$, $F$ be two vector bundles over $X$. If $f: E|_Y \ra F|_Y $ is an
isomorphism, then there exists an open set $U$ containing $Y$ and an
extension $f: E|_U\ra F|_U$ of $f$,  which is an isomorphism.
\eproposition
Proof. Since $F$ is a section of $\Hom{E|_Y}{F|_Y}$ it can be extended
to a section of $\Hom{E}{F}$. Let $U$ be a set of those points
for which this map is an isomorphism. Then $U$ is open and contains
$Y$. \qed

We shall now prove Theorem \ref{K-homotopy}. Our goal it to show that 
the isomorphism class of the bundle $F^*_tE$, where $F:Y\times I\ra X$ 
is the homotopy 
between the maps $f$ and $g$, does not depend on $t$. 
Let $\pi : Y\times I \ra Y$ be the projection onto
the first factor. For any $t\in I$, the bundles $F^*E$ and $\pi^*F^*_t E$ are 
identical when restricted to the closed subspace $Y\times \{t\}$ of $Y\times I$. 
Using the fact that $Y$ is a compact space and the previous theorem 
we conclude that the bundles $F^*E$ and $\pi^*F^*_t E$ are isomorphic
on some strip $Y\times\delta t$, where $\delta t$ is a neighbourhood
of $t$ in $I$. Hence the isomorphism class $[F_t^*E]$ of 
the bundle $F_t^*E$ is a locally constant function of $t$. 
Since $I$ is connected, this function is in fact constant and 
this implies that 
$$
f^* E = F_0^*E \simeq F_1^* E= g^*E
$$
as required. The second statement of Theorem \ref{K-homotopy} is now clear. 
\qed 

\bremark 
Using homotopy invariance of the $K$-functor it is now easy to see that 
if $X$ is a contractible space compact space, then every vector bundle over $X$ is 
trivial. 
\eremark
\mbox{}

We now pass to the description of the excision properties of $K^0$. 
We shall keep our assumption that $X$ is a compact Hausdorff space
and that $Y$ is a closed subspace. If we excise $Y$ from $X$, i.e., 
form the space $X-Y$, the resulting space is no longer compact, but
only locally compact. 
\bdefinition
By $X/Y$
we mean $X$ with $Y$ shrunk to a point, which we shall denote by
$\infty$. Thus $X/Y$ is the one-point compactification of the 
locally compact space $X-Y$. 
When $Y$ is the empty set, we take $X^+ =X/\emptyset $ to be
$X$ with a base point. 
\edefinition

In order to understand what the operation $X\ra X/Y$ does to the $K$-theory, 
it will be useful to have some description of vector bundles on 
the quotient space $X/Y$ in terms of bundles on $X$. 
Let us assume that $E$ is a vector bundle over $X/Y$ and that $E_\infty = 
\bc^r$. Since the canonical map $\pi: X\ra X/Y$ identifies $Y$ with $\infty$,
the restriction of the  
the pullback bundle $\pi^*E$  to $Y$ is the trivial bundle $\pi^*E|_Y= Y\times \bc^r$. 
We shall show that there is a one-one correspondence
between vector bundles over the quotient space $X/Y$ and vector 
bundles over $X$ whose restriction to $Y$ is a trivial bundle.

Suppose then that $E$ is a vector bundle over $X$ and that $\alpha :E|_Y \ra
Y\times V$ is an isomorphism. We call $\alpha$ a {\em trivialization}
of $E$ over $Y$. Let us denote by $\pi: Y\times V\ra V$ the projection
onto the second factor. We 
define an equivalence relation on $E|_Y$ by setting 
$e\simeq e'$ if and only if
$$
\pi\alpha (e) = \pi \alpha (e').
$$
This relation identifies the points in the restriction of $E$ to $Y$
which are \lq on the same level' relative to the trivialization
$\alpha$. We then extend this relation trivially to the whole of $E$.

\blemma The set of equivalence classes $E_\alpha$ of this equivalence
relation  is a vector bundle over $X/Y$.
\elemma
Proof. We note that since the trivialization procedure described above does not
affect the bundle away from $Y$, the only problem here is to check that
the bundle $E_\alpha$ is locally trivial near $\infty$.
For this, use Tietze's theorem to extend $\alpha $ to an isomorphism
$\tilde{\alpha} : E|_U \ra U\times V$, where $U$ is some open set containing $Y$.
Then $\tilde{\alpha}$ induces an isomorphism $(E|_U)_{\tilde{\alpha}} \simeq
U \times V$. We remark that $U$ cannot be taken to  be $X$ in general.
\qed

Suppose now that $\alpha_0$ and $\alpha_1$ are homotopic trivializations of
$E$ over $Y$. This means that we have a trivialization $\beta$ of $E\times I$
over $Y\times I\subset X\times I$ inducing $\alpha_0$ and $\alpha_1$ over the
end points of $I$. Let
$$
f: (X/Y) \times I \lrar  (X\times I  ) / ( Y\times I)
$$
 be the natural map. Then $f^*((E\times I)_\beta)$ is a bundle on
$(X/Y) \times I$ whose restriction to $X/Y\times \{ i \}$ is
$E_{\alpha_i}$ for $i= 0,1$. Thus
$$
E_{\alpha_0} \simeq E_{\alpha_1}
$$
and we have established the following.
\bproposition\label{alpha-homotopy}
A trivialization  $\alpha$ of a bundle $E$ over $Y\subset X$ defines a bundle
$E_\alpha$ over $X/Y$. The isomorphism class of $E_\alpha $ depends only on
the homotopy class of $\alpha$.
\eproposition

Now that we know how to describe vector bundles over the quotient space
$X/Y$ we may consider $K^0(X/Y)$, and define the relative $K$-theory
$K^0(X,Y)$. 

\bdefinition
The relative $K$-group $K^0(X,Y)$ is
$$
K^0(X,Y) = \tilde{K}^0(X/Y)
$$
\edefinition
It follows from the properties of $\tilde{K}^0$ that the relative $K$-theory
is a contravariant functor of $(X,Y)$. We also note that
since $K^0(X) = \tilde{K}^0(X^+)$,  we have $K^0(X,\emptyset) = K^0(X)$.

\bproposition
Let $X$ be a compact space and $Y$ a closed subspace. There is 
an exact sequence
$$
K^0(X,Y) \arr{\jmath ^*} K^0(X) \arr{\imath^*} K^0(Y)
$$
where $\imath: Y\ra X$ and $\jmath :(X,\emptyset ) \ra (X,Y)$ are
inclusions.
\eproposition
Proof. Since the composite map $\imath^* \jmath^*$ is induced by the
composition $\jmath\imath : (Y,\emptyset) \ra (X,Y)$, we have that
$\imath^*\jmath^* = 0$. This yields $\ig \jmath^* \subset \ker \imath^*$.

Now
 suppose that $\xi\in \ker \imath^*$. We may represent $\xi$ as the 
formal difference $\xi = [E] -
[\theta^n]$ for some vector
bundle $E$ and a  trivial bundle of rank $n$ over $X$. Since $\imath^*\xi = 0 $
it follows that $[E|_Y] = [\theta^n|_Y]$ in $K^0(Y)$, which means that the
bundles $E|_Y$ and $\theta^n|_Y$
are stably isomorphic.  This implies that 
we have, for some integer $m$,  $(E\oplus \theta^m ) |_Y = \theta^n\oplus
\theta^m $.  In other words, we  have constructed a trivialization
$\alpha$ of the bundle $(E\oplus \theta^m ) |_Y$. This defines a
bundle  $
(E\oplus \theta^m) _\alpha$ on $X/Y$ and so an element
$$\eta = [(E\oplus \theta^m )_\alpha] - [\theta^n\oplus\theta^ m ] \in
\tilde{K}^0(X/Y) = K^0(X,Y)
$$
of the relative $K$-theory. 
Then
$$
\jmath^* (\eta) = [E \oplus\theta^ m ]- [\theta^n \oplus \theta^m ] = [E] -
[\theta^n] = \xi
$$
which proves that $\ker\imath^* \subset\ig\jmath^*$. \qed
\bcorollary
Let $X$ be a compact space, and $Y$ a closed subspace, equipped
with a base point. The sequence
$$
K(X,Y)\lrar \tilde{K}(X) \lrar \tilde{K}(Y)
$$
is exact.
\ecorollary
Proof. Let $y_0$ be the base point of $Y$, and so a base 
point of $X$. The result follows from the fact
$$
\begin{array}{rcl}
K(X) & \simeq & \tilde{K}(X) \oplus K(y_0)\\
K(Y) & \simeq & \tilde{K}(Y) \oplus K(y_0)
\end{array}
$$
and the previous statement. \qed

\subsubsection{More on the Serre-Swan Theorem}\label{more-on-swan}

We shall now proceed to give an algebraic description of the module of
sections of a  vector bundle over $X/Y$. Since $X/Y$ is a compact space, 
the Serre-Swan theorem applies, and so the space $\Gamma(X/Y,\bar{E})$
of continuous  sections of a bundle $\bar{E}$ over $X/Y$ is 
a projective module over $C(X/Y)$. We want to  show this explicitly, 
by constructing 
a free module in which $\Gamma(X/Y, \bar{E})$ is a direct summand.

We have seen that any bundle $\bar{E}$ over $X/Y$ can be described using a
bundle $E$ over $X$ together  with a trivialization $\phi : E|_Y \simeq
Y\times \bc^r$ of the restriction of $E$ to $Y$. Sections of $\bar{E}$ are 
thus sections of $E$ which
are constant over $Y$. Using the isomorphism $\phi$ we can say that
a section of $\bar{E}$ over $X/Y$ is a  section $s$  of $E$ such that
$\phi\circ s :Y \ra \bc^r=\widetilde{\bc^r}$ is constant.  The space of sections
$\Gamma(X/Y, \bar{E}) $ is a module over the algebra of functions
$C(X/Y)$, where functions act on sections by pointwise multiplication.
We note here that a function $f\in C(X/Y)$ is the same as
a function on $X$ which is constant on $Y$.

To focus our attention, let us  consider a trivial bundle $E$ over 
$X$ together with an isomorphism $\phi:E|_Y \stackrel{\sim}{\ra} \widetilde{\bc^r}$ 
which
trivializes $E$ over $Y$.
This   isomorphism $\phi$ may be regarded
as an element of $GL_r(C(Y))$, the group of invertible $r\times r$ matrices
with entries in $C(Y)$. As we have seen before, the bundle $E$ together
with the isomorphism $\phi$ may be regarded as a bundle $\bar{E}$ over
the quotient space $X/Y$. We note that even when the bundle $E$ is 
trivial, the quotient bundle $\bar{E}$ need not be.  We want to 
find maps
$$
\Gamma(X/Y, \bar{E})\ra  C(X/Y) ^N \ra \Gamma (X/Y, \bar{E})
$$
for some $N$, which compose to the identity map.

 We first construct the map on the right. Tietze's extension 
theorem shows that there is a surjection $C(X) \ra C(Y)$. 
We use this
map to lift $\phi$ to $M_r(C(X))$, and we shall denote the resulting matrix by
$p$. This matrix need no longer be invertible, so let us assume
that $\phi^{-1}$ lifts to a matrix $q$. The matrix valued functions
$p$ and $q$ then satisfy $pq|_Y = qp|_Y = 1$, where $1$ denotes the
identity in the group $GL_r(C(Y))$.

For a vector $\xi\in \bc^r$, consider $q\xi\in C(X)^r= \Gamma(X,E)$. Then
clearly $\phi(q\xi)|_Y$ is a constant function on $Y$ and thus $q\xi$ is a
section of $\bar{E}$ over $X/Y$. This is also true if $\xi $ is any
element of $C(X/Y)^r$. This remark defines  a map $C(X/Y)^r \ra \Gamma
(X/Y, \bar{E})$ by sending $\xi\mapsto q\xi$.

Now consider $s\in \Gamma(X/Y, \bar{E})$. Hence  $s\in
\Gamma(X,E)$ is such that $\phi\circ s|_Y $ is constant. 
Since $ps|_Y = \phi \circ s |_Y$ is also constant, this gives a map
$\Gamma(X/Y, \bar{E}) \ra C(X/Y)^r$, $s\mapsto ps$.

Passing to the general case, we want to establish the
following sequence
$$
\Gamma(X/Y, \bar{E}) \arr{\textstyle{\left(\begin{array}{c} p\\
\beta\end{array}\right)}} C(X/Y)^{2r} \arr{(q\;\; \alpha)}  \Gamma(X/Y,
\bar{E} )
$$
where $\alpha, \beta \in M_r(C(X))$ vanish on $Y$. We want to find $\alpha,
\beta, p,q$ such that $p $ and $q$  are lifts of $\phi$ and $\phi^{-1}$ as
before and that the composite map above is the identity, that is
\beqn\label{condition}
qp + \alpha \beta = 1
\eeqn
Let us denote by $I$ the ideal of matrix-valued functions vanishing  on
$Y$. Condition \ref{condition} is equivalent to  $qp\equiv 1\, \mod \, I^2$.
Now it follows from the definition of 
 $p$ and  $q$ that they satisfy  $qp \equiv 1 \, \mod \, I$.
We then write $qp = 1 - (1- qp)$ and upon multiplying this identity 
on the left by  $1+(1-qp)$ we get
$$
(1+(1-qp))qp = 1- (1-qp)^2
$$
Hence if we replace $q$ by $\tilde{q}= (1+(1- qp))q$ then $\tilde{q} p
\equiv 1\, \mod \, I^2$ and $\tilde{q}|_Y = q|_Y$. We can then choose
$\alpha$ and $\beta $ to be any elements of $I$, for instance put
$\alpha = \beta = 1 - qp$.

The example just described corresponds to the following purely algebraic
situation. Let $R$ be a ring and let $I$ be an ideal in $R$. 
Let us denote by $\tilde{I}$ the augmentation of $I$. We establish 
contact with the geometric case described before by means of the 
following identification:
\beqn\label{IDENTIFY}
\begin{array}{rcl}
A & = & C(X)\\
A/I &  = & C(Y)\\
\tilde{I} &  = & C(X/Y)
\end{array}
\eeqn
The objects on the left are related together by the following cartesian square:
$$
\begin{array}{ccc}
\tilde{I}& \arr{}  & \bc \\
\downarrow && \downarrow\\
A &  \arr{\pi} & A/I
\end{array}
$$
Let $\phi\in GL_r(A/I)$. Following Milnor \cite[p.20]{Milnor}, let us define
$M(A^r, \bc^r, \phi)$
to be an $\tilde{I}$-module
consisting  of all $s\in A^r$ such that $\phi(\pi s)$  is a constant function with values in 
$\bc^r$, and hence determines an element of $\bc^r$. We claim
that $M(A^r, \bc^r, \phi)$ is a projective finitely generated
module over $\tilde{I}$. More precisely, if $p\in M_r(A)$ is such that 
$p\cong \phi\; \mod\; I$ and $q\in M_r(A)$ is such that $q\cong 
\phi^{-1}\; \mod\; I$ then we have the following maps
$$
M(A^r, \bc, \phi) \Ar{\scriptstyle{\left(\begin{array}{c} p\\ 1-qp
\end{array}\right)}} (\tilde{I})^{2r} \Ar{(1+(1-qp)q\;\;  (1-qp))} 
M(A^r, \bc^r, \phi)
$$
whose composition is the identity map. We see that the module $M$
is the image of the following projection operator:
$$
M(A^r, \bc^r, \phi) = \ig e = \ig \textstyle{\left(\begin{array}{c} p\\ 1 -
qp\end{array}\right)}((2q - qpq)\;\; (1-qp))
$$
where the matrix $e$ is an idempotent over $\tilde{I}$. We can express
the same fact in yet another way when we notice that
there is the following identity
\beqn\label{MAGIC}
\left(\begin{array}{cc} \phi & 0 \\ 0 & \phi^{-1}\end{array}\right)
= \larray 1 & \phi \\ 0 & 1 \rarray \larray 1 & 0 \\
-\phi^{-1} & 1 \rarray \larray  1 & \phi \\ 0 & 1 \rarray
\larray 0 & - 1 \\ 1 & 0 \rarray
\eeqn
Now  this matrix can be lifted to an invertible matrix
\beqn
\omega =\larray 1 & p \\ 0 & 1 \rarray \larray 1 & 0 \\ -q & 1 \rarray \larray
1 & p \\ 0 & 1 \rarray \larray 0 & -1 \\ 1 & 0 \rarray
\eeqn
Then the projection operator associated
to $\phi$  over $A/I$, whose image is the module $M(A^r, \bc^r, \phi)$
is given by
\beqn\label{projection}
e = \omega \larray 1 & 0 \\ 0 & 0 \rarray \omega^{-1}
\eeqn
We shall shortly define higher $K$-groups. In particular, elements 
of $K^1(Y)$ will be represented by invertible matrices. The above construction
will then be used to provide a connecting homomorphism 
$K^1(Y)\ra K^0(X/Y)$. 

\subsubsection{Higher $K$-groups and Bott periodicity}
If a space $X$ is locally  compact, we can define the $K$-theory
of $X$ using the one-point compactification $X^+$ of $X$ and 
the reduced $K$-theory: 
$$
K^0(X) = \tilde{K}^0(X^+)
$$
For example, in the case when $X$ is a compact space and 
$Y$ is a closed subspace, we have
$$
K^0(X\backslash Y) = \tilde{K}^0((X\backslash \,Y)^+)
= \tilde{K}^0(X/Y)
$$
which makes connection with the convention we have used earlier. 
\bdefinition
For any natural $n$ we define 
$$
K^{-n}(X) = K^0(X\times \br^n) = K^0(S^n(X))
$$
where $S^n(X)= X \times \br^n$ is  called 
the $n$-th reduced suspension of $X$. 
\edefinition

One of the most important properties of these $K$-groups is 
that they have the excision property, which means that 
the following long sequence is exact. 
$$
\ra K^{-n-1}(Y) \ra K^{-n}(X\backslash Y)\ra K^{-n}(X) \ra
K^{-n} (Y) \ra K^{-n}(Y) \ra \dots \ra K^0(Y)
$$
Another quite amazing fact is that the above sequence 
is in fact periodic with period two. This follows from 
the Bott periodicity theorem which implies that 
$$ 
K^0(S^2(X))\simeq K^0(X)
$$
We can therefore say that we have the following exact sequence
of length six
$$
\begin{array}{ccccc}
K^0(X\backslash Y)  & \arr{} & K^0(X) & \ra & K^0 (Y)\\
\uar{} &                &            &        &  \dal{}\\
K^{-1}(Y) & \larr{} & K^{-1}(X) & \larr{} & K^{-1}(X\backslash Y)
\end{array}
$$

\subsection{Algebraic $K$-theory}
\subsubsection{$K_0(\Lambda)$ and projective $\Lambda$-modules}
Let $\Lambda$ be a ring with a unit. We shall now investigate how much
we can say about the $K$-groups of $\Lambda$ proceeding in a purely 
algebraic way. We have already taken some preliminary steps towards
establishing the necessary algebraic formalism in Section \ref{more-on-swan}. 
As we can only hope to sketch a broad outline of this beautiful theory, 
the reader should consult the literature for more details on the
subject. We recommend especially Rosenberg's book \cite{Rosenberg:book}
and a more or less up-to-date overview \cite{LL}. 

Let us denote by ${\bf P}(\Lambda)$ the category of finitely 
generated projective left $\Lambda$-modules. This is a subcategory 
of the category ${}_\Lambda {\rm Mod}$ of all left $\Lambda$-modules. 
\bdefinition
We define $K_0(\Lambda)$ to be the Grothendieck group of $[{\bf P}(\Lambda)]$, 
the semigroup of isomorphism classes of elements of ${\bf P}(\Lambda)$. 
\edefinition

\bexample
Let $\Lambda = k $ be a field, and let $\cal V$ be the category of
vector spaces of finite dimension over $k$, which are the same as 
finitely generated projective modules over $k$. Since two vector 
spaces $V$ and $W$ are isomorphic if and only if their dimensions
are the same, the map $\dim : [{\cal V}] \ra \bn$  establishes 
an isomorphism of abelian semigroups (we assume that $0\in \bn$). 
Then $K_0(k) = K_0(\bn) = \bz$. 
\eexample

If $\phi: \Lambda \ra \Lambda'$ is a ring homomorphism, then 
we may equip $\Lambda'$ with the structure of a right $\Lambda$-module
by $\lambda ' \cdot \lambda = \lambda' \phi(\lambda)$. 
Thus $\phi$ induces a functor $\Phi: {}_{\Lambda}{\rm Mod} \ra 
{}_{\Lambda'}{\rm Mod}$ by $M\mapsto \Lambda'\otimes_\Lambda M$. 
The functor $\Phi$ maps free modules to free modules, preserves
the property of being finitely generated, projective, etc. 
Thus it induces a functor $[{\bf P}(\Lambda)]\ra [{\bf P}(\Lambda')]$
and so a functor $\Phi: K_0(\Lambda) \ra K_0(\Lambda')$. 
Of course, this remark is just an example of the  general 
statement (\ref{K-universal}). (Compare also 
\cite[Theorem 1.1.3]{Rosenberg:book}.)

Finally, when $\Lambda$ is a commutative ring, and $M$, $N$ are 
two modules over $\Lambda$, then the tensor product $M\otimes _\Lambda
N$ has the structure of a $\Lambda$-module $\lambda(m\otimes_\Lambda n)= 
\lambda m \otimes_\Lambda n= m\otimes_\Lambda \lambda n$. 
It follows that the category ${\bf P}(\Lambda)$ is closed under 
the tensor product $\otimes _\Lambda$ and this operation induces
a product in $K_0$: $[P]\cdot [Q] = [P\otimes_\Lambda Q]$ for 
all modules $P$, $Q$ in ${\bf P}(\Lambda)$. 

\bexample
If $\Lambda = \bz$ is the ring of integers, then the map that 
sends a module to its rank induces an isomorphism $K_0(\bz) = \bz$, 
and so all projective $\bz$-modules are free. 

A similar result holds for local rings \cite{Milnor}, and principal 
ideal domains. Recall that an element $\lambda$ of the ring $\Lambda$
is called a {\em unit} if there exists a $\mu\in \Lambda$
such  that $\lambda\mu = \mu\lambda = 1$. The multiplicative
group of units in $\Lambda$ is denoted by $\Lambda^\times$. 
A ring $\Lambda$ is called a {\em local ring} iff the set 
$M= \Lambda - \Lambda^\times$ consisting of all non-units
is a left ideal, and therefore also a right ideal. 
Now if $\Lambda$ is a local ring, then every finitely generated
projective module is free \cite{Milnor,Kaplansky}, and 
so $K_0(\Lambda)$ is the free cyclic group generated by 
$[\Lambda]$. 
\eexample

There is an equivalent description of $K_0(\Lambda)$ in terms 
of idempotents, i.e. maps $p$  such that $p^2 = p$. 
This approach is particularly useful in defining the 
$K$-theory of operator algebras, but it is also quite convenient
to have this description in the present setting. We have already 
touched on this problem in our discussion of the Swan theorem. 
Any projective $\Lambda$-module $P$ is a direct summand 
in a free module $\Lambda^n$, for some natural $n$. Let 
us define a map $p: \Lambda^n \ra \Lambda^n$ by putting
$p$ equal the identity map on $P$ and $0$ on $Q$. Then 
it is clear that $p$ is an idempotent, and that $p$ determines
the $\Lambda$-module $P$ up to isomorphism. We shall denote by $GL_n(\Lambda)$ 
the group of invertible $n\times n$ matrices with entries in $\Lambda$. 
By definition, 
this is the group of units in the algebra of $n\times n $ matrices
$M_n(\Lambda)$. 

We construct the following two direct systems. 
First we note that there is an inclusion  $M_n(\Lambda) \hookrightarrow
M_{n+1}(\Lambda)$ given by $a\mapsto \larray a & 0 \\ 0 & 0\rarray$. 
Let us denote by $M(\Lambda)$ the ring of \lq infinite matrices', 
that is the direct limit of this direct system. $M(\Lambda)$ may 
be considered as the union of all $M_n(\Lambda)$, and so every 
element in $M(\Lambda)$ is in fact a matrix of finite size. 
Secondly, we have the inclusion $GL_n(\Lambda) \hookrightarrow
GL_{n+1}(\Lambda)$ defined by $a \mapsto \larray a & 0\\ 0 & 1
\rarray$. We denote by $GL(\Lambda)$ the direct limit of this
system. Inside $M(\Lambda)$ we single out the set ${\rm Id}(\Lambda)$ 
of idempotent matrices. $GL(\Lambda)$ acts on ${\rm Id}(\Lambda)$
by conjugation. It turns out that the set of orbits of the 
action of $GL(\Lambda)$ on ${\rm Id}(\Lambda)$ may be identified
with the semigroup $[{\bf P}(\Lambda)]$ \cite[p.~8]{Rosenberg:book}. 
\btheorem 
The set $[\rm Id]$ of orbits of the action of $GL(\Lambda)$ on ${\rm Id}(\Lambda)$
is a semigroup, with the semigroup operation induced by 
$(p,q) \mapsto \larray p & 0\\ 0 & q\rarray$. Moreover, 
$K_0(\Lambda) = K([\rm Id])$ . 
\etheorem
To illustrate how this point of view works in practice, let us 
prove the following important property of $K_0$ known as 
the Morita invariance. 
\btheorem 
There is a natural isomorphism $K_0(\Lambda) \ra K_0(M_n(\Lambda))$, 
for any ring $\Lambda $ and a natural number $n$. 
\etheorem
Proof. This result rests on the observation that $M_k(M_n(\Lambda))
= M_{kn}(\Lambda)$. We use this remark to identify 
the group $GL(M_n(\Lambda))$ with $GL(\Lambda)$ and
the set ${\rm Id}(M_n(\Lambda))$ with ${\rm Id}(\Lambda)$. 
Now apply the previous theorem. \qed

We record also the fact that $K_0$ is continuous in the following
sense. Let $\{\Lambda_\alpha\}_{\alpha \in A}$ be a direct 
system of rings indexed by some partially ordered set $A$, and 
let $\Lambda = \displaystyle{\lim_{\arr{}}}\, \Lambda_\alpha$ be 
the direct limit of this system. If we apply the $K_0$-functor
to this direct system, we obtain a direct system of 
abelian groups $\{K_0(\Lambda_\alpha)\}$. 
\btheorem 
$$
K_0(\Lambda) = \lim_{\arr{}} K_0(\Lambda_\alpha)
$$
\etheorem 

\subsubsection{$K_1(\Lambda)$ and the Whitehead Lemma}
There is another interesting subgroup of $GL_n(\Lambda)$ that 
we may single out. 
 Let  $E_n(\Lambda)$ be the subgroup consisting of {\em elementary matrices}
 $e^\lambda_{ij}$ which differ from the identity matrix
only by an  element $\lambda$ in the $ij$-th 
place,  where $i\neq j$.  Then for any $n\geq 3$ the 
group $E_n(\Lambda)$ is perfect,  which means that
$$
[E_n(\Lambda), E_n(\Lambda)] = E_n(\Lambda)
$$
where the group on the left is generated by commutators 
$[x,y] = xyx^{-1}y^{-1}$. 

The direct system used to define $GL(\Lambda)$ restricts to a direct system of 
the groups $E_n(\Lambda)\hookrightarrow E_{n+1}(\Lambda)$, whose direct limit will 
be denoted $E(\Lambda)$. 
The key result linking the properties of $GL(\Lambda)$ and
$E(\Lambda)$ is the following lemma. 
\blemma {\rm [Whitehead]} 
$$
[GL(\Lambda), GL(\Lambda)] = E(\Lambda)
$$
\elemma
Proof. We first need to express the commutator $[x,y]$ of 
two elements of $GL_n(\Lambda)$ by a matrix in 
$GL_{2n}(\Lambda)$: 
\beqn\label{HELP}
\larray
xyx^{-1}y^{-1} & 0 \\
0 & 1 \rarray
= \larray x & 0 \\ 0 & x^{-1} \rarray
\larray y & 0 \\0 &  y ^{-1} \rarray
\larray (yx)^{-1} & 0 \\ 0 & yx \rarray
\eeqn
Then we use the formula (\ref{MAGIC}) to express each 
matrix on the right in terms of elementary matrices in 
$E_2(M_n(\Lambda))$. Finally, use the relations
satisfied by the elementary matrices to show that each of 
the matrices on the right is in fact an element of the
group $E_{2n}(\Lambda)$. (See \cite[Lemma 2.1.2]{Rosenberg:book}
for details).  \qed 

\bdefinition 
$$
K_1(\Lambda) = GL(\Lambda)/ E(\Lambda)
$$
\edefinition

It follows from  Whitehead's lemma  that 
$$
K_1(\Lambda)= GL(\Lambda)_{ab}= 
GL(\Lambda)/[GL(\Lambda), GL(\Lambda)]
$$ 
On the other hand we recall from group homology that 
$GL(\Lambda)_{ab}$ is the first homology group of $GL(\Lambda)$ 
with integer coefficients, 
and so we have the 
isomorphism $K_1(\Lambda) = H_1(GL(\Lambda), \bz)$.
Moreover, a ring homomorphism $\phi: \Lambda \ra \Lambda' $
induces a group homomorphism $GL(\Lambda) \ra GL(\Lambda')$ 
that restricts to a homomorphism $E(\Lambda)\ra E(\Lambda')$.
Hence there  is an induced map
$$
K_1(\phi) : K_1(\Lambda) \arr{} K_1(\Lambda')
$$

\bexample (cf. \cite[p.~27]{LL}\cite{Milnor}) Let us assume that 
$\Lambda$ is a commutative ring. Then taking the determinant 
of a matrix extends to a homomorphism 
$\det : GL(\Lambda) \ra \Lambda^\times $. We shall denote 
by $SL(\Lambda)$ the kernel of this map:
$$
SL(\Lambda) = \ker(\det) = \bigcup _ {n=1}^\infty SL_n(\Lambda)
$$
where, as usual, $SL_n(\Lambda) = \ker\; [\det: GL_n(\Lambda)\ra
\Lambda^\times]
$. The obvious inclusion $E(\Lambda) \subset SL(\Lambda)$ 
implies that $\det $ descends to a group homomorphism 
$$
\det : K_1(\Lambda) \arr{} \Lambda^\times
$$
There is therefore the following short exact sequence
$$
1 \ra SK_1(\Lambda) \arr{} 
GL(\Lambda)/E(\Lambda)\arr{\det} \Lambda^\times
\ra 1
$$
where $SK_1(\Lambda) = SL(\Lambda)/E(\Lambda)$. 
This sequence splits, since $\det $ has a right inverse given by 
the map $\Lambda^\times = GL_1(\Lambda) \hookrightarrow
GL(\Lambda) \ra K_1(\Lambda)$. Thus we see that 
$K_1(\Lambda) = \Lambda^\times \oplus SK_1(\Lambda)$
and so, in the commutative case, the computation of $K_1(\Lambda)$
reduces to finding $SK_1(\Lambda)$. There are some 
special cases when $SK_1(\Lambda)$ is trivial,  making
$\det$ an isomorphism. This happens, for instance, 
when $\Lambda$ is a Dedekind domain, a local ring 
or a finite commutative ring. So in particular, if $\Lambda = \bz$
$K_1(\bz) = \bz^\times = \bz_2$. In the case when 
$\Lambda$ is a field $k$, we have $K_1(k) = k^\times = 
K_1(k[x])$. 
\eexample
\subsubsection{Relative $K$-groups and Milnor's construction}\label{relative-K-section}
We have already encountered the relative $K$-theory in the topological 
setting, where it appeared as a natural object associated with the 
operation of taking the quotient of a topological space. 
In the category of rings it is clearly useful to consider 
short exact sequences of the form
$$
0\ra I \ra R \arr{q} R/I\ra 0
$$
where $I$ is an ideal in $R$. Since the map $q$ is a ring homomorphism, 
it induces a map of $K$-groups
$$
q_*: K_0(R)\ra K_0(R/I)
$$
We shall now define a relative $K$-group $K_0(R,I)$ in such a way 
that it will fit into the following short exact sequence
$$
K_0(R,I) \arr{} K_0(R)\arr{q_*}K_0(R/I)
$$
An important point to be considered is the relation between 
$K_0(R,I)$ and $K_0(I)$. 
\bdefinition\label{relative-K}
If $I$ is an ideal in a ring $R$ then the double of 
$R$ along $I$ is the following subring of the direct product 
$R\times R$:
$$
D(R,I) = \{ (x,y) \in R\times R\mid x-y \in I\}
$$
\edefinition
The double $D(R,I)$ relates to the ring $R$ and the 
ideal $I$ via the following short exact sequence
$$
0 \ra I \ra D(R,I) \arr{p_1} R\ra 0
$$
where $p_1: R\times R \ra R$ is the projection onto the first coordinate.
This sequence splits, and the splitting map is given by 
the diagonal embedding of $R$ into the direct product 
$R\times R$. The kernel of the projection map $p_1$ may 
be identified with the ideal $I$. 
\bdefinition
The relative $K_0$-group of a ring $R$ and an ideal $I$
is defined by 
$$
K_0(R,I) = \ker\{ p_{1*}: K_0(D(R,I)) \ra K_0(R)\}
$$
Similarly, we define
$$
K_1(R,I) = \ker\{p_{1*}:K_1(D(R,I))\ra K_1(R)\}
$$
\edefinition
\bremark \label{nonunital-I}
Milnor \cite{Milnor} uses the above formulae to define the 
$K$-groups of the ideal $I$. This does not lead to any difficulty
in the case of $K_0(I)$, since in that case $K_0(R,I)$ does not
depend on the \lq ambient' ring $R$, and only on the ideal $I$. 
Therefore, if we define 
$$
K_0(I) = \ker \{ K_0(\tilde{I})\ra
K_0(\bz)\}
$$
we can write $K_0(R,I) = K_0(I)$ (cf. \cite[Theorem 1.5.9]{Rosenberg:book}).
This is our first excision statement in algebraic $K$-theory.  
The same statement
is not, in general, true for the group $K_1(R,I)$, which does depend
on the ring $R$. See \cite[Remark 4, p.~34]{Milnor} and \cite{SwanExc}.
\eremark 
In order to understand the general structure of the excision-type results 
for the groups $K_0$ and $K_1$, we shall briefly describe Milnor's patching
theorem \cite[p.~19 ff.]{Milnor} that leads to the Mayer-Vietoris sequence. 
Let us begin with the
following cartesian square of rings
$$
\begin{array}{ccc}
A & \arr{\imath_2} & A_2\\
\dal{\imath_1} &&\dar{\jmath_2}\\
A_1& \arr{\imath_1} & B
\end{array}
$$
Thus, in particular, for any two $a_1\in A_1$ and $a_2\in A_2$ 
such that $\jmath_1(a_1) = \jmath_2(a_2)$ there exists a unique
$a\in A$ such that $\imath_k(a) = a_k$, for $k=1,2$. 
If $P$ is a finitely generated $A$-module, we shall denote by 
$\imath_{1\nt}P $ the induced module 
$\imath_{1\nt}P = P\otimes_{A} A_1$. We note that there is a canonical 
isomorphism 
$$
\jmath_{1\nt}\imath_{1\nt}P = \jmath_{2\nt}\imath_{2\nt}P
$$
which is easy to check when we notice that this equation is the same as
$$
(P\otimes_A A_1)\otimes_{A_1}B = P\otimes_A(A_1\otimes_{A_1}B)
= P\otimes_A B = (P\otimes_A A_2)\otimes_{A_2}B
$$

Let $P_1$ and $P_2$ be finitely generated projective modules over $A_1$ 
and $A_2$, respectively. Let us assume that there is an isomorphism
$h:\jmath_{1\nt} P_1 \simeq \jmath_{2\nt}P_2$ of $B$-modules. 
We define $M= M(P_1,P_2,h) \subset P_1\times P_2$ by 
$$
M= \{ (p_1, p_2) \in P_1\times P_2\mid h\jmath_{1*}p_1 = \jmath_{2*}p_2\}
$$
where $\jmath_{i*}p_i = p_i\otimes 1 \in P_k\otimes_{A_i} B$, $k=1,2$. 
Thus $M$ appears in the following cartesian diagram
\beqn\label{diagram}
\begin{array}{ccc}
M & \arr{} & P_2\\
\dal{} && \dar{\jmath_{2*}}\\
P_1 &\arr{h\jmath_{1*}}  & \jmath_{2\nt}P_2
\end{array}
\eeqn
\btheorem\label{patching} {\em [Milnor's patching theorem]}
If either $\jmath_1$ or $\jmath_2$ is surjective, then $M$ is a finitely 
generated $A$-module such that $\imath_{k\nt}M \simeq P_k$, $k=1,2$. 
\etheorem 
Proof. Assume first that the modules $P_k$ are free:
$P_1 = A_1^r$, $P_2 = A_2^s$. Note that in this case $\jmath_{1\nt} P_1
= B^r$ and $\jmath_{2\nt}P_2= B^s$. Assume that the isomorphism $h: B^r \ra B^s$
may be lifted to an isomorphism $h':A^r_1\ra A^s_2$, so that $\jmath_1 h'=h$.
We form two modules associated with this data: $M(A^r_1, A^s_2, h)$ and 
$M(A^s_1, A^1_s, id)\simeq A^s$ and then use the two corresponding diagrams of the 
form \ref{diagram} to prove that $M(A^r_1, A^s_2, h) \simeq A_1^s$ is therefore
a free module. 

Suppose now that the map $\jmath_2$ is surjective and $P_1$ and $P_2$ are free modules, 
as in the previous step. We claim that the module $M(P_1, P_2, h)$ is projective, 
where $h: (\jmath_{1\nt} P_1 = B^r) \simeq (\jmath_{2\nt} P_2
= B^s)$. Let $h^{-1}: B^s\ra B^r$ be the inverse map. Then, by naturality, 
$$
M(A^r_1, A_2^s, h) \oplus M(A^s_1, A^r_s, h^{-1}) = M(A^{r+s}_1, A^{r+s}_2, 
h\oplus h^{-1})
$$
where  $h\oplus h^{-1}$ is represented by the matrix $\larray h & 0 \\ 0 & h^{-1}\rarray$. 
Now use the formula \ref{MAGIC} to lift the last matrix to an invertible matrix
over $A_1^{r+s}$. This step reduces our problem to the previous case, 
and shows that $M(A^{r+s}_1, A^{r+s}_2, 
h\oplus h^{-1})$ is a free module, and therefore $M(P_1, P_2, h)$ is projective. 

Finally, when $P_1$ and $P_2$ are finitely generated projective modules over 
their respective rings, we can find projective modules $Q_1$ and $Q_2$
such that $P_k\oplus Q_k \simeq A_k^r$, $k=1,2$. Let $h: P_1\otimes _{A_1}B \arr{\sim}
P_2\otimes _{A_2}B$ be an isomorphism postulated by the theorem's hypothesis. 
Since we have
$$
P_k\otimes _{A_k} B\oplus Q_k\otimes _{A_k} B\simeq B^r
$$
we can write:
$$
\begin{array}{rcl}
Q_1\otimes_{A_1} B \oplus B^r & = & Q_1\otimes_{A_1} B \oplus
(P_1\otimes _{A_1} B\oplus Q_2\otimes _{A_2} B)\\
& = & Q_2\otimes_{A_2} B \oplus B^r
\end{array}
$$
This shows that if we replace $Q_k$ by $Q_k \oplus A^r$, we can claim that 
$h$ lifts to an isomorphism
$k: Q_1\otimes_{A_1} B \simeq Q_2\otimes_{A_2}B$. Then the identity
$$
M(P_1,P_2, h) \oplus M(Q_1, Q_2,k) = M(P_1\oplus Q_1, P_2\oplus Q_2, h\oplus k)
$$
shows that $M(P_1,P_2, h)$ is finitely generated and projective by 
the previous case. \qed
\bcorollary\label{Corollary}
Under the same hypotheses as in Milnor's theorem the sequence
$$
K_k(A)\arr{(\imath_{1*}, \imath_{2*})} K_k(A_1)\oplus K_k(A_2) \arr{\jmath_{1*}  -
\jmath_{2*}} K_k(B)
$$
is exact for $k=0,1$. 
\ecorollary
In fact, we can do better than this. Let us define a map $\delta:
K_1(B)\ra K_0(A)$ as follows. Take $[u]\in K_1(B)$, where $u\in GL_n(B)$. 
This matrix determines an isomorphism $h: \j_{1\nt}A_1^n \simeq \j_{2\nt}A_2^n$
of free $B$-modules. Thus there is a projective module $M=M(A_1^n, A_2^n,h)$
over $A$. Define
$$
\delta[u] = [M]-[A^n]\in K_0(A)
$$
We leave it to the reader to check that $\delta$ is a well defined homomorphism. 
Thus combining the short exact sequences of Corollary \ref{Corollary}
with the connecting homomorphism $\delta$ we get a Mayer-Vietoris sequence of 
length six:
$$
K_1(A) \ra K_1(A_1)\oplus K_1(A_2)\ra K_1(B)\arr{\delta}
K_0(A) \ra K_0(A_1)\oplus K_0(A_2)\ra K_0(B)
$$
\subsubsection{Relative $K$-groups and excision}\label{simple-excision}
The Mayer-Vietoris sequence  is the direct analogue of the 
sequence of the same name known from algebraic topology. We shall now demonstrate
a few applications of the techniques introduced in the previous section.
\bproposition\label{I-tilde}
If $I$ is an ideal in a ring $R$ then 
$$
K_0(I) = K_0(R,I)
$$
\eproposition
Proof. Recall that $K_0(I)=\ker\{\epsilon: K_0(\tilde{I})\ra K_0(\bz)\}$
where $\epsilon$ is the augmentation homomorphism:
$$
0\ra I \ra \tilde{I} \arr{\epsilon} \bz \ra 0
$$
and $\tilde{I} = I\semileft \bz$. 
Let us denote $D=D(R,I)$ and let us define a unital homomorphism $\gamma: \tilde{I}\ra
D$ by $\gamma: (x,n)\mapsto 
(n\cdot 1, n\cdot 1 + x)$ for $n\in \bz$ and $x\in I$. Then there 
is the following cartesian diagram
$$
\begin{array}{ccc}
\tilde{I} & \arr{\gamma} & D\\
\dal{\epsilon\;\;} && \dar{p_1}\\
\bz &\arr{\imath} & R
\end{array}
$$
where $\imath :\bz\ra R$ is the inclusion. This diagram satisfies
the conditions of Milnor's theorem, so there is the corresponding Mayer-Vietoris 
sequence:
$$
\ra K_1(R) \arr{\delta} K_0(\tilde{I})\arr{(\gamma_*, \epsilon_*)} K_0(D) \oplus K_0(\bz)
\arr{p_{1*}  - \imath_*} K_0
$$
Let us take $[u]\in \ker \epsilon_* = K_0(I)$. Then 
$$
0= (p_{1*}  - \imath_*)(\gamma_*, \epsilon_*)[u] = p_{1*}\gamma_*[u]
$$
which shows that $\gamma_*[u]\in \ker p_{1*} = K_0(R,I)$. Thus the homomorphism 
$\gamma_*$ restricts to a homomorphism $\gamma_*: K_0(I)\ra K_0(R,I)$. 
We want to show that this map is an isomorphism. First we note that 
$$
\ker \{\gamma_*:K_0(I)\ra K_0(R,I)\} = \ker\gamma_*\cap \ker\epsilon_* 
= \ig \delta
$$
by exactness of the Mayer-Vietoris sequence. Let us then recall the definition
of the connecting homomorphism $\delta$. Let $[u]\in K_1(R)$ be represented
by an invertible matrix $u\in GL_n(R)$. Let us form the module 
$M(\bz^n, D^n, u)$. If we identify 
$D^n$ with the subset of $R^{2n}$ consisting of the elements $(p_1,p_2)$ 
such that $p_k\in R^n$ and $p_1-p_2\in I^n$, then we see that the invertible 
matrix $u$ lifts to an invertible matrix $\larray u & 0\\ 0 & u \rarray \in 
GL_n(D) \subset GL_{2n}(R)$. Thus $M(\bz^n, D^n, u)=\tilde{I}^n$ is a free 
$\tilde{I}$-module. Thus 
$$
\delta[u] = [M] - [\tilde{I}^n] = 0
$$
This implies that $\ker \{\gamma_*:K_0(I)\ra K_0(R,I)\} =0$. 

To show that $\gamma_*$ is surjective, let $ \xi\in K_0(R,I)$. Then 
$p_{1*}\xi = 0$ and 
$$
K_0(D)\oplus K_0(\bz) \ni(\xi,0) 
\stackrel{ p_{1*}  - \imath_*}{\displaystyle\longmapsto} 0
$$
There exists thus a class $\eta\in K_0(\tilde{I})$ such that 
$$
(\gamma_*, \epsilon_*)\eta = (\xi,0)
$$
Since evidently $\epsilon_* \eta =0$, $\eta\in K_0(I)$. This finishes the 
proof. 
\qed

Another important application of the Mayer-Vietoris sequence is given 
by the following prototype of the excision property in the algebraic 
$K$-theory. 
\btheorem \label{excision-prototype}
For any ideal $I\subset R$ we have the exact sequence of $K$-groups
$$
K_1(R,I) \ra K_1(R) \ra K_1(R/I) \ra K_0(R,I) \ra K_0(R) \ra K_0(R/I)
$$
\etheorem
This statement follows in a straightforward way from the Mayer-Vietoris sequence 
applied to the cartesian square
$$
\begin{array}{ccc}
D(R,I) & \arr{p_2} & R \\
\dal{p_1\;}&& \dal{}\\
R &\arr{} & R/I
\end{array}
$$
when we proceed in a similar way to the proof of Proposition \ref{I-tilde}. 

\subsubsection{Quillen's higher $K$-theory}
Let $G$ be a discrete group, and let $EG$ be a contractible
space on which $G$ acts freely. The quotient space 
$BG = EG/ G$ is called the {\em classifying space} of 
the group $G$. The classifying space is used to compute the 
homology of $G$; by definition 
$H_*(G,\bz) = H_*(BG, \bz)$, so that the group homology
$H_*(G,\bz)$ may be computed in terms in homology of 
the topological space $BG$.  The homotopy groups of 
$BG$ are rather simple, since $BG$ is connected, 
$\pi_0(BG) = 0 $, $\pi_1(BG)= G$ and $\pi_n(BG) = 0 $ 
for $n\geq 2$. In other words, $BG$ is the Eilenberg-Maclane
space of type $K(G,1)$. 

\bdefinition
An $H$-space is a space $X$ equipped with maps
$$
\mu : X\times X \ra X, \qquad i : X\ra X
$$
and a point $e$, which satisfy the axioms of a group 
up to homotopy. This means that
\begin{enumerate}
\item The maps $X^3\ra X$ that send
$(x, y , z) $ to $\mu(x, \mu(y,z))$ and $\mu(\mu(x,y), x)$ 
are homotopic. 
\item The maps $X\ra X$ sending $x $ to $\mu(e, x)$ and
$\mu(x, e)$ are homotopic to the identity. 
\item The map $X\ra X$ sending $x$ to $\mu(x, i(x))$
 is homotopic to the constant map $x\mapsto e$, 
and similarly for the map $x\mapsto \mu(i(x), x)$. 
\end{enumerate}
\edefinition

We record here the fact that the  fundamental group of an 
$H$-space is abelian and that the Hurewicz map 
$h : \pi_1(X)\ra H_1(X, \bz)$ is an isomorphism. 

Let us now consider the group $GL(\Lambda)$ defined above, 
and let $BGL(\Lambda)$ be its classifying space. 
Quillen \cite{Quillen-Ktheory} defines a space $BGL(\Lambda)^+$ 
with the following properties. First of all, $BGL(\Lambda)^+$ is 
an $H$-space, and secondly, the canonical inclusion
$$
\phi: BGL(\Lambda) \arr{} BGL(\Lambda)^+
$$
induces an isomorphism in homology
$$
H_*(BGL(\Lambda), \bz)\arr{\simeq}H_*(BGL(\Lambda)^+, \bz)
$$
\bdefinition
For $n > 0 $ we define
$$
K_n(\Lambda) = \pi_n(BGL(\Lambda)^+)
$$
\edefinition
\bexample
Using the Hurewicz map we find that 
$$
K_1(\Lambda) = \pi_1(BGL(\Lambda)^+) \arr{\simeq}
H_1(BGL(\Lambda)^+, \bz) = H_1(GL(\Lambda),\bz)
$$
and $H_1(GL(\Lambda), \bz) = GL(\Lambda)/[GL(\Lambda), 
GL(\Lambda)] = GL(\Lambda)/ E(\Lambda)$ and 
we recover $K_1(\Lambda)$ defined in the previous section. 

Furthermore, it can be shown that $K_2(\Lambda)
= H_2(E(\Lambda), \bz)$, and that $K_3(\Lambda)= 
H_3(St(\Lambda), \bz)$, where $St(\Lambda)$ is the 
so called Steinberg groups associated with 
$\Lambda$ \cite{Milnor}. 
\eexample

Quillen has shown that for any positive $n$ the groups 
$K_n(\Lambda)$ are finitely generated, if $\Lambda$ is the ring of 
integers in a number field. Unfortunately, 
the algebraic $K$-groups are hard to calculate, even 
the $K$-theory of the integers $\bz$ is not known 
in all generality. So far, apart from the examples mentioned
above, it has been shown that $K_2(\bz) = \bz_2$, 
$K_3(\bz) = \bz_{48}$. It has been suggested for some time that $K_4(\bz) = 0$, 
but it seems that a recent proof of this fact was withdrawn. The problem 
of calculating these $K$-groups has deep relations with 
arithmetic geometry and number theory. 

\subsection{$K$-theory of $C^*$-algebras}
\bdefinition
A $C^*$-algebra is a Banach algebra $A$ over $\bc$ equipped with a 
conjugate-linear involution $*: A \ra A$ such that 
$$
(xy)^* = y^*x^*, \;  \norm{x^* x} = \norm{x}^2
$$
for all $x,y$ in $X$. 
\edefinition
It a basic fact of the theory of $C^*$-algebras
that the norm on a $C^*$ algebra is uniquely determined by the 
algebra structure. The spectral radius $R(x)$ of an element $x$
of $X$ is defined by 
$$
R(x) = \sup \{ |\lambda|, \lambda\in \bc
 \mid  x - \lambda \; \mbox{ \rm is not invertible}\}
$$
Then the square of the  norm  of any element $x$ of a $C^*$-algebra is the 
same as the spectral radius of $x^*x$:
$$
\norm{x}^2 = R(x^*x)
$$
The above definition introduces the so called {\em abstract}
$C^*$-algebras. {\em Concrete} $C^*$-algebras are defined 
in terms of their representations in a Hilbert space $\bh$. 
In fact an involutive algebra $A$ is a $C^*$-algebra if and only if 
it admits a $*$-representation $\pi$ on a Hilbert space  $\bh$ such that 
$\pi(x) = 0 $ implies that $x=0$ and the image $\pi(A)$ of $A$ 
is closed in norm in the algebra of bounded operators on $\bh$. 

We remark also that a theorem of Gelfand and Naimark states that 
for any unital commutative algebra $A$ there exists a compact 
topological space $X$ such that $A$ is the algebra of 
continuous complex-valued functions on $X$ (cf. 
\cite[II.1]{Connes:Book}). 
To describe this space, let us define a {\em character} of $A$
to be a homomorphism $\chi: A \ra \bc$ preserving the unit of 
$A$, i.e. $\chi(1) = 1$. The space of all such homomorphisms,
 called the {\em spectrum} of $A$ is compact when equipped 
with the topology of pointwise convergence on $A$. 
The {\em Gelfand transform} is the following map
$$
A \ni x \mapsto (\hat{x}: \mbox{\rm Spec } (A ) \ra \bc)
$$
where $ \hat{x}(\chi) = \chi(x)$. The Gelfand transform
establishes an isomorphism of $A$ onto 
the $C^*$-algebra $C(X)$ of continuous functions on $X$. 

We can therefore say that in the commutative case the 
algebra of continuous functions $C(X)$ on a topological 
space carries the same amount of information as the space 
itself. This observation lies at the basis of noncommutative 
geometry, which has been developed by Alain Connes. 
It allows one to extend most of the classical geometric and analytic 
tools used for the study of a compact space $X$ to the case when 
the corresponding algebra of coordinate functions is noncommutative.

To define $K_0$ of a $C^*$-algebra, let us form the $C^*$-algebra
$M(A)$ of infinite matrices. We define it using the same 
direct system of matrix algebras as 
considered in the purely algebraic case. This time, however, 
this direct system is considered in the category of $C^*$-algebras, 
and the direct limit $M(A)$ is a normed algebra with a $C^*$-algebra
norm, but not necessarily complete. 
Let $p,q$ be two projections, that is self adjoint idempotents in 
$M(A)$. We say that $p$ and $q$ are homotopic in $M(A)$ if and only 
if  
there exists $v\in M(A)$ such that 
$v^*v = p $, $vv^* = q$. 

\bdefinition
If $A$ is a unital $C^*$-algebra then the group 
$K_0(A)$ is defined as the 
Grothendieck group of the abelian semigroup $P(A)$ of the homotopy 
classes of projections in $M(A)$. The addition in $P(A)$
is given by 
$$
[ p ] + [ q ] = \larray p & 0 \\ 0 & q \rarray\in M_{k+m} (A) 
\subset M (A)
$$
for $p,q \in M_{k,m}(A)$. 
\edefinition

We can now follow the pattern established in the previous two sections to 
give the following basic properties of $K_0(A)$.
\begin{enumerate}
\item {\em Homotopy invariance}. Let $A$ and $B$ be 
$C^*$-algebras. We say that two morphisms
$\alpha, \beta: A\ra B$ are {\em homotopic} if there exists a family 
 of morphisms $\gamma_t$, where $t\in [0,1]$ such that
$\gamma_t(a)$ is a norm continuous path in $B$ for every fixed $a\in A$. 
Moreover, $\gamma_0 = \alpha $ and $\gamma_1= \beta$. 
\end{enumerate}
\bproposition
If $\alpha_0, \alpha_1: A\ra B$ are homotopic morphisms, 
the induced maps $K_0(A)\ra K_0(B)$ coincide. 
\eproposition
\begin{enumerate}
\setcounter{enumi}{1}
\item There is a natural isomorphism $K_0(A\oplus B) \simeq
K_0(A)\oplus K_0(B)$. This may be regarded as a special 
case of the following property. 
\item{\em  Half-exactness}.  If $J$ is a nonunital algebra, let us denote 
by $\tilde{J} = \bc \semiright J$ the augmented algebra. Then 
define
$$
K_0(J) = \ker(K_0(\tilde{J}) \ra \bz)
$$
\end{enumerate}
\bproposition
Let 
$$
0 \ra J\ra A\ra B \ra 0
$$
be an exact sequence of $C^*$-algebras. Then the sequence
$$
K_0(J) \ra K_0(A) \ra K_0(B) \ra 0
$$
is exact at $K_0(A)$.
\eproposition
\begin{enumerate}\setcounter{enumi}{3}
\item {\em Morita equivalence}
The map $ x \mapsto \larray x & 0 \\ 0 & 0 \rarray$ 
induces an isomorphism $K_0(A) \ra K_0(M_2(A))$. 
\end{enumerate}
\bdefinition
For an algebra $A$ we define the {\em suspension} $SA$ of the 
algebra $A$
 $$
SA = \{ f : S^1 \ra A\mid f(1) = 0 \}
$$
Then we define $K_1(A) = K_0(SA)$. 
\edefinition
\begin{enumerate}
\setcounter{enumi}{4}
\item{\em  The index map. }
\end{enumerate}
\bproposition
Let us consider a short exact sequence of algebras as in the 
previous Proposition. Then there exists a connecting morphism 
$\delta : K_1(A/J ) \ra K_0(J)$. This map is called the index
map in $K$-theory. 
\eproposition
\begin{enumerate}
\setcounter{enumi}{5}
\item {\em Bott periodicity. }
$$K_i(A)\simeq K_{i + 2}(A)$$
\item  {\em Excision. }
\end{enumerate}
\bproposition 
Let $0 \ra J \ra A\ra A/J \ra 0 $. Then there 
is the following exact sequence of length six. 
$$
\begin{array}{ccccc}
K_0(J)  & \arr{} & K_0(A) & \ra & K_0 (A/J)\\
\uar{} &                &            &        &  \dal{}\\
K_1(A/J) & \larr{} & K_1(A) & \larr{} & K_1(J)
\end{array}
$$
\eproposition
\begin{enumerate}
\setcounter{enumi}{7}
\item {\em Stability} For any algebra $A$ and the algebra of 
compact operators $\bk$ we have $K_0(A\otimes \bk) = K_0(A)$. 
\item {\em Continuity}
\end{enumerate}
\bproposition
Let $A$ be the $C^*$-algebra direct limit of the directed 
system $\{A_i\}_{i\in I}$. Then there is the corresponding
direct system of $K$-groups $\{K_0(A_i)\}_{i\in I}$
and 
$$
K_0(A) = \lim_{\arr{}}K_0(A_i)
$$
\eproposition

We now list some basic examples of the $K$-theory groups \cite[6.5]{wegge}. 
\bexample
\begin{enumerate}
\item $K_0(\bc) = \bz $, $K_1(\bc)  = 0 $. The same result is 
true for the algebra $M_n(\bc)$ of matrices over $\bc$, by the 
Morita equivalence. 
\item Let $\mbox{\bf K } $ be the algebra of compact operators. 
Then $K_0(\bk) = \bz$ and $K_1(\bk) = 0 $. For the algebra 
$\bb$ of bounded operators we have $K_0 = K_1 = 0$. Hence 
from the long exact sequence of $K$-theory groups it follows 
that $K_0(\bb/\bk) = 0 $ and $K_1(\bb/\bk) = \bz$. The algebra 
$\bb/\bk $ is called the Calkin algebra. 
\item Let $\bt_n$ be the $n$-torus. 
Then $K_i(C(\bt_n )) = \bz^{2^{n-1}}$, for $i= 0,1$. 
In  the case of the sphere $S^n$ we have two possibilities. When $n$ is even, 
$K_0(C(S^n)) = \bz/2\bz$ and $K_1(C(S^n))= 0 $, whereas 
for $n$ odd $K_0= K_1 = \bz$. 
\end{enumerate}
\eexample
Cuntz has proved a very interesting characterization of functors of 
$K$-theory type. Namely any functor $F$ from the category of 
$C^*$-algebras to the category of abelian groups which is 
continuous, half-exact, homotopy invariant and stable is the same 
as the $K$-functor. More precisely, if $F(\bc) = \bz$ and $F(S\bc)= 0$, 
then $F = K_0$, whereas in the opposite case when $F(\bc) = 0$ and 
$F(S\bc)= \bz$, $F = K_1$. See \cite{Cuntz:Bott} for more details.
\subsubsection{Algebraic $K$-theory of $C^\star$-algebras}
Since any unital Banach   (or $C^*$)-algebra  $R$ is a ring, it makes sense 
to define both algebraic $K^{\rm alg}(R)$ and \lq topological'  
$K^{\rm top}(R)$ $K$-theories  for $R$.  Here by $K^{\rm top}_i(R)$ we mean the two $K$-groups defined in the 
previous section. 
Comparing the definitions given above we see that $K_0^{\rm top} (R) = 
K_0^{\rm alg}(R)$. In other degrees there is a map 
\beqn\label{K:compare}
K_i^{\rm alg}(R) \arr{} K_i^{\rm top} (R)
\eeqn
which is induced by the (identity) map  $GL(R)^{\rm discrete} \ra GL(R)^{\rm top}$. 
Here $GL(R)^{\rm discrete}$ is the group $GL(R)$ treated as a discrete group, 
whereas in forming $GL(R)^{\rm top}$ we take into account the 
topology of the algebra $R$. The map of $K$-theories is not an isomorphism 
in general, as is expected from the differences in their properties. For 
instance, the Bott periodicity which holds for $C^\star$-algebras (or Banach 
algebras in general) is not shared by the algebraic $K$-theory.

When $i=1$ we have the following descriptions of the $K_1$-groups in both cases: 
$$
K^{\rm alg}_1(R) = GL(R)/E(R), \qquad K^{\rm top}_1(R) =  GL(R)/GL_0(R)
$$
where $GL_0(R)$ is the connected component of $GL(R)$. Since the group 
of elementary matrices $E(R)$ lies in $GL_0(R)$, the natural map 
is surjective. However, it is not even injective  for $\bc$ \cite[p.455]{Rosenberg}, 
since $K_1^{\rm alg}(\bc)$ is isomorphic to the multiplicative group $\bc^\times$  via the 
determinant map, whereas $K_1^{\rm top}(\bc) = 0 $. 

In general, there are results about the comparison of the two $K$-theories only 
for special classes of algebras. For instance, there is a conjecture of Karoubi 
that the map (\ref{K:compare}) is an isomorphism for stable $C^\star$-algebras. 
We say  that a $C^\star$-algebra $A$ is {\em stable} if and only if 
$A\hat{\otimes}_\pi \bk \simeq A$, where $\bk $ is the algebra of compact operators, 
and $\hat{\otimes}_\pi$ denotes the completed projective tensor product. 
We shall comment on Karoubi's conjecture in the context of Suslin-Wodzicki 
excision theorem in $K$-theory. 

\subsection{$K$-homology}
$K$-homology arose from the study of elliptic operators on manifolds and 
it is the dual theory to $K$-theory. The initial ideas are due to Atiyah \cite{Atiyah-ellips}, 
which were later on developed by Brown-Douglas-Fillmore, Kasparov and Connes.  

First we need to explain a relation between $K$-theory and the theory of 
Fredholm operators. Let us then assume that $\bh$ is a separable 
Hilbert  space. A Fredholm operator $T$ on $\bh$ is a bounded linear 
operator on $\bh$ such that $T$ has a finite dimensional kernel and 
cokernel. 
Then the index of $T$ is defined by 
$$
{\rm index}\; T = \dim \ker T   - \dim \coker T
$$
The index of the Fredholm operator is invariant under compact perturbations: 
$$
{\rm index}\; (T+A) = {\rm index}\; T
$$
for any compact operator $A$. Moreover, if $S$ is a bounded operator sufficiently 
near in norm to $T$ then $S$ is also Fredholm and ${\rm index}\; T = 
{\rm index} \, S$. 
The last statement may be reformulated as follows. If we topologise 
the space $\calf$ of Fredholm operators on $\bh$ using  the norm topology, then we 
have defined a continuous map 
$$
{\rm index}: \calf \arr{} \bz
$$
This map is in fact surjective, and so induces a bijection 
$$
\pi_0(\calf) \arr{} \bz
$$
where $\pi_0(\calf)$ is the set of connected components of $\calf$. 

Suppose now that $X$ is a compact space and that $T: X \ra \calf$ is a continuous 
map, so that $T_x$ is a family of Fredholm operators depending continuously 
on the parameter $x\in X$. If $\dim \ker T_x$ is independent of $x$, then the family 
of vector spaces $\ker T_x$ forms a vector bundle $\ker T$ over $X$. A similar 
statement is true in the case of  the cokernel of $T$. 
We can thus define the index of the family of operators $T_x$ by 
$$
{\rm index }\; T = [\ker \; T] - [\coker\; T] \in K^0(X)
$$
If $\ker T_x$ has a variable dimension, one can still reduce to the previous 
case by composing with a projection. Atiyah \cite{atiyah} shows that 
there is a bijection of  the set $[X, \calf]$ of homotopy classes 
of maps $T: X\ra \calf$ onto $K^0(X)$. The composition of operators 
in $\calf $ corresponds to the addition in $K^0(X)$, and the adjoints 
correspond to taking the negative.

Let $P$ be an elliptic  operator of order $m$ on a closed manifold $M$. 
The operator $P$ may be acting on
functions, vector-valued functions, or more generally, on sections 
of vector bundles over $M$. Ellipticity implies that there exists 
 an operator $Q$, called the 
{\em parametrix}, such that $PQ$ and $QP$
are equal to the identity operator, 
modulo operators of lower order. Moreover, for any function $f$, 
the operator $Pf- fP$ is an operator of degree $m-1$. 

If we allow the operator $P$ to be a pseudodifferential operator, 
then as show Atiyah and Singer, it is 
enough to concentrate on operators of order 0. We then say that 
$P$ is a Fredholm operator if and only if there exists a 
pseudodifferential operator $Q$ such that $PQ - 1$ and 
$QP - 1$ are compact operators, and for any function $f$ , 
$Pf - fP$ is a compact operator. 

We can  represent classes in $K$-homology $K_0(X)$ using   homotopy 
classes of such Fredholm operators. To see that this really describes 
an object dual to $K$-theory, let us take a class $[P]\in K_0(X)$ and 
$[e]\in K^0(X)$. Then we have the following pairing
$$
([e], [P]) \arr{} {\rm index}\; (ePe) \in \bz
$$

This idea is extended  by the following important notion \cite{Connes:IHES}. 
\bdefinition
An {\em even} Fredholm module is a quadruple $(A,\bh ,  F, \gamma)$ where
\begin{enumerate}
\item $A$ is a $C^*$-algebra acting on a Hilbert space $\bh$. 
\item $F$ and $\gamma $ are self-adjoint anticommuting 
involutions, i.e., $F^* = F$, $\gamma ^* = \gamma$, 
$F^2 = \gamma^2 = 1$, $F\gamma + \gamma F = 0 $. 
\item For all $a\in A$, $[F,a]$ is a compact operator in $\bh$. 
\end{enumerate}
\edefinition
\bremark
There is a similar definition  of an {\em odd } 
Fredholm module, where we do not assume the existence of the involution 
$\gamma$. In point (3) of the above definition we can replace
the ideal of compact operators by any Schatten ideal 
$\call^p(\bh)$. In that case we call the resulting Fredholm 
module $p$-summable. 
\eremark
\bexample
Let $\bh= \call^2(S^1)$ be the space of $\call^2$-integrable 
functions on the circle. This is a separable Hilbert space with 
a basis given by Fourier modes. It has the  direct sum 
decomposition $\bh = \bh ^+ \oplus \bh^-$ into spaces 
of functions whose Fourier series expansion contains 
positive (negative) frequencies, respectively. Let $\gamma$ 
be the corresponding involution, i.e. we assume that 
$\bh^{\pm}$ are the $\pm 1$ eigenspaces of $\gamma$. 
The involution $F$ is given in this case 
by the Hilbert transform $e^{in\theta} \mapsto e^{-in\theta}$. 
We may describe this transform explicitly using the following 
formula (cf. e.g. \cite[p. 83]{Pressley-Segal})
$$
F(f)(\theta) = \frac{1}{2\pi}P.V. \int^{2\pi}_0 \sum_k\sign\; (k)
e^{ik (\theta -\phi)} f(\phi) d\phi
$$
where $P.V. \int_0^{2\pi} = \lim _{\epsilon \to 0 }\left( 
\int^{\theta - \epsilon }_0 + \int_{\theta + \epsilon} ^{2\pi}\right)
$ is the principal value of the integral. 
A theorem of Fefferman and Sarason \cite{Power-Hankel} then 
states that the operator $[F, f]$ is compact for a function
 $f$ on $S^1$ if and 
only if $f$ is of Vanishing Mean Oscillation. In particular, this 
is satisfied by any continuous function on the circle.  This result 
indicates the possible choice of function algebra on $S^1$ which 
we may consider. By requiring that $f$ belongs to a  
suitable Besov space we may  obtain a $p$-summable Fredholm module. 
\eexample
\bexample \cite[IV.5]{Connes:Book} Let $\Gamma$ be a discrete group. We define 
the reduced $C^\star$ algebra $C^\star_r(\Gamma)$ associated with $\Gamma$ 
as follows. Consider the Hilbert space $l^2(\Gamma)$ (this may be 
understood as the space of $L^2$-functions on $\Gamma$ with the counting measure), 
and let $\{\xi _h\}_{h \in \Gamma} $ be an orthonormal basis. There is 
a natural action of $\Gamma$ on $l^2(\Gamma)$ given by 
$$
g \mapsto \lambda _g; \qquad \lambda_g(\xi_h ) =  \xi _{gh}
$$
This way we define a homomorphism from $\Gamma$ to the group 
of unitary operators in the Hilbert space $l^2(\Gamma)$. 
The algebra $C^*_r(\Gamma)$ is by definition the $C^*$-algebra
generated by these symbols, i.e. we have 
$\lambda^*_g = \lambda_{g^{-1}}$ and $\lambda_g\lambda_h = \lambda_{gh}$. 

Consider now a tree $T$ on which the group acts freely and 
transitively. $T$ is a one-dimensional simplicial complex 
which is connected and simply connected. Take $T^i$ to be 
the set of vertices when $i= 0 $ and edges when $i = 1$. For $p\in T^0$ 
put $\phi : T^0 - \{ p \} \ra T^1$ to be the map which sends 
a vertex $q$ to the unique $1$-simplex containing $q$ which is a 
subset of the path $[p,q]$. Let us define $\bh ^ + = l^2(T^0) $ 
and $\bh ^ - = l^2(T^1)$ which gives an orthogonal decomposition 
of the Hilbert space $\bh$. Furthermore, let us define an 
operator $P$ on $\bh$ by $P\xi_p = 0 $ and $P \xi_q = \xi _{\phi(q)}$. 
Then there is the following involution $F$ on $\bh$ 
anticommuting with the involution $\gamma$ whose eigenspaces
are the Hilbert subspaces $\bh^{\pm}$: 
$$
F = \left( \begin{array}{cc} 0 & P^{-1} \\ P & 0 
\end{array}
\right)
$$
This way we have defined an even Fredholm module 
over the algebra $C^*_r(\Gamma)$. This Fredholm module was constructed 
by Julg and Valette \cite{Julg-Valette} and may be used to 
prove an important theorem by Pimsner and Voiculescu that 
the algebra $C^*_r(\Gamma)$ has no nontrivial idempotents. 
\eexample
\bexample \label{conformal}
Let $M$ be an even dimensional manifold, $\dim M = 2l$, which 
is compact, conformal, oriented. Let 
$\bh = L^2(M, \Lambda^l _\bc T^*)$ be the Hilbert space 
of square-integrable sections of the 
bundle of differential forms in the middle degree, equipped with 
the scalar product 
$$
\langle \omega_1, \omega_2 \rangle = \int_M \omega_1 \wedge * \bar{\omega}_2
$$
Here $*$ denotes the Hodge star-operator, which is determined 
by the conformal structure. Let $\gamma = (-1)^{l(l-1)/2} * $ be 
the corresponding involution. Singer posed the following question: Does 
there exist an involution $F$ anticommuting with 
$\gamma$ which is equal to the identity operator 
when restricted to the space of exact forms? Connes, Sullivan and 
Teleman gave a positive answer to this question using the following 
construction of a Fredholm module. 
There exists a partial isometry $S: \bh^- \ra \bh^+$
such that $\overline{{\rm Im}\; d} $ is the graph of $S$, where 
$$
d: C^\infty(M , \Lambda^{l-1}_\bc T^*)\ra C^\infty(M, \Lambda^l T^*)
$$
and where we denote by $\bh^{\pm}$ the $\pm 1$ eigenspaces 
of the involution $\gamma$. Let $\bh^l$ be the cohomology group
of $M$ in the middle dimension $l$, which we identify with 
the Hilbert space of harmonic forms.  Put 
$$ 
\bh ' = \bh \oplus \bh^l = (\bh \ominus \bh^l) \oplus( \bh ^ l \oplus \bh^l)
$$
Let $F$ be the  involution given by 
$\left(\begin{array}{cc} 0 & S \\ S^* & 0 \end{array}\right)$ 
on the first summand and by $\left( \begin{array}{cc} 0 & 1 \\ 1 & 0 
\end{array}\right)$ on the second. We can now  state the following theorem 
of Connes, Sullivan and Teleman \cite{CST}. 
\btheorem
\begin{enumerate}
\item The triple $( \bh' , F , \gamma)$ is a Fredholm module over 
$C^\infty(M)$ canonically associated with the conformal structure on $M$. 
\item This Fredholm module uniquely determines a conformal structure on $M$. 
\end{enumerate}
\etheorem
This theorem generalizes to quasi-conformal manifolds, where 
the Fredholm module constructed by it is used to derive 
local formulae for Pontryagin classes. Also, we can say that 
this Fredholm module serves as the noncommutative 
analogue of a  manifold equipped with a conformal structure. It is 
then natural to try to find the correct notion of differentiable 
manifold in noncommutative geometry.  This has 
been done by Connes, Ellwood and the author in \cite{BCE}. 
\eexample

\subsection{$KK$-theory}\label{KK}
An important generalization of both $K$-theory and $K$-homology 
has been proposed by Kasparov.  A central notion of his theory 
is the {\em Kasparov module} which we can define as follows. 
\bdefinition
Let $A$ and $B$ be $C^*$-algebras. An {\em even Kasparov module}
is a triple $(\cale, F, \gamma)$ where
\begin{enumerate}
\item $\cale$ is a countably generated Hilbert space over $B$. This means
that it is a complete space and a 
 right $B$-module which admits a scalar product with values in $B$. 
 \item $A$ is represented as a $*$-algebra by bounded operators on $\cale$. 
\item $F$ is a bounded operator on $\cale$ such that for all $a\in A$, 
$(F^2 - 1)a$, $(F^* - F) a $ and $[F,a]$ are compact operators on $\cale$. 
\item $\gamma $ is a self-adjoint bounded involution on $\cale$, anticommuting 
with $F$ and commuting with the action of $A$. This means 
that  $\gamma^2 = 1$, $\gamma^* = \gamma$,  $F\gamma = - \gamma F$,
and for all $a\in A$ $\gamma a = a \gamma$. 
\end{enumerate}
\edefinition
As in the case of Fredholm modules, an {\em odd} Fredholm module
is a pair $(\cale , F)$ satisfying the first three conditions above. 
We define $KK(A, B)$ to be the set of homotopy classes of even Kasparov 
modules. The set of homotopy classes of odd Kasparov modules 
defines $KK_1(A,B)$. 
There is a natural abelian group structure on $KK(A,B)$ which is induced 
by the direct sum of Kasparov modules. The link between $KK$-theory 
and $K$-theory and $K$-homology is provided by the identifications
$$
KK(\bc, A) = K_0(A), \qquad KK(A,\bc) = K^0(A)
$$

A good introduction to $KK$-theory is given in the article by Cuntz \cite{Cuntz:KK}
and in his lecture notes \cite{Cuntz:Elements}; see also the article by 
Higson \cite{Hig}. 

\vfill\eject

\section{Cyclic homology }
Cyclic cohomology was discovered by Connes in 1981. It was introduced 
as an extension of the de Rham cohomology of differentiable manifolds
to the noncommutative setting, and 
serves as a natural target for the Chern character homomorphism from $K$-theory. 
In this chapter we shall concentrate on the formal description of this theory,  
and its relation to $K$-theory will be explained in the next chapter. 
Recent work of Cuntz and Quillen \cite{CQ1,CQ2,CQ3}
has provided  new insights into the foundations of the theories of cyclic
type and, as we shall see later, may be used in the case 
of homology theories associated with topological algebras. 
We begin our discussion of cyclic homology with the 
introduction of the differential graded algebra 
of noncommutative forms $\Omega A$. 

\subsection{Noncommutative differential forms}\label{FORMS}
\subsubsection{The  harmonic decomposition}

Let us assume that $A$ is a unital algebra over $\bc$. The 
algebra of abstract differential forms $\Omega  A$ \cite{Arveson,Karoubi,CQ1} 
over 
$A$ (also called noncommutative differential forms) is
the universal differential graded algebra generated by the algebra $A$ together 
with symbols $da$, $a\in A$ subject to the following relations: 
\begin{enumerate}
\item $da$ is linear in $a$. 
\item Leibniz rule: $d(ab) = d(a) b + adb$. 
\item $d(1)=0$ 
\end{enumerate}
The above conditions imply that $d: A \ra \Omega A$ is a derivation.  In the following we shall 
fix the algebra $A$ and denote $\Omega = \Omega A$. 
We may describe 
$\Omega  $ more explicitly as follows. First, it is a 
graded vector space, which in degree $n$ is isomorphic to 
$\Omega^n = A\otimes \bar{A}^{\otimes n }$, where $\bar{A} = A/\bc$. 
The map that establishes this isomorphism 
is given explicitly by 
$$
a_0da_1\cdots da_n \leftrightarrow (a_0, \dots, a_n)
$$
where $(a_0, \dots, a_n)$ is the image of the simple tensor 
$a_0\otimes \cdots \otimes a_n$ in the space $A \otimes\bar{A}^{\otimes n}$. 
In the basis given by this above isomorphism, the differential 
$d$ is defined by the following formula
$$
\begin{array}{rcl}
d(a_0da_1\cdots da_n) &  =  & da_0\cdots da_n \\
d(a_0, \dots , a_n )  & = & (1, a_0, \dots, a_n)
\end{array}
$$
The graded vector space $\Omega$ becomes a graded algebra when 
introduce the following multiplication of forms
$$
\begin{array}{rcl}
(a_0da_1\cdots da_n) \cdot (a_{n+1}da_{n+2} \cdots da_{k})&   = & 
(-1)^n a_0a_1da_2\cdots da_k \\
&& \mbox{} +  \sum_{i=1}^n(-1)^{n-i} a_0 da_1\cdots d(a_i a_{i+1})
 \cdots da_{k}
\end{array}
$$
Using this definition, we may extend  $d$  to a graded derivation on $\Omega$: 
$$
d(\omega\eta) = d(\omega) \eta + (- 1)^{\deg(\omega)} \omega d\eta
$$
which makes $\Omega $ into a differential graded algebra. 
Finally, $d\Omega^n \simeq \bar{A}^{\otimes n }$. 
\bremark\label{multiple-tensors}
The space of $1$-forms $\Omega^1 $ plays an important role
of its own. First of all, it is the kernel of the multiplication 
map $m: A\otimes A \ra A$, and so it fits into the 
short exact sequence 
$$
0\ra \Omega^1  \arr{i} A\otimes A \arr{m} A
$$
The inclusion map $i$ is given explicitly by 
$$
i(a_0da_1) = a_0(a_1\otimes 1 - 1\otimes a_1)
$$
This inclusion makes the $A$-bimodule structure of $\Omega^1 $ 
quite explicit. The other important property of the 
bimodule $\Omega^1 $ is that, together with the derivation 
$d:A \ra \Omega^1 $  it provides a universal model for derivations
of $A$ with values in an $A$-bimodule. This means that for 
any derivation $D: A\ra M$ of $A$ with values in an $A$-bimodule
$M$, there exists a unique bimodule map $\phi: \Omega^1 \ra 
M$ such that $D = \phi d$. 

It it possible to define the space of $n$-forms 
$\Omega^n$ as the $n$-fold tensor product of the $A$-bimodule
$\Omega^1A$ over $A$:
$$
\Omega^n A = \Omega^1 A \otimes_A \stackrel{n}{\dots} \otimes_A
\Omega^1 A
$$
Let us denote a typical element of this tensor product by 
$\omega_1\cdot \dots \cdot \omega_n$, with $\omega_i\in \Omega^1A$, 
for $i=1, \dots, n$. Then the product on $\Omega A$ may be defined
as the following map $\Omega^n\times \Omega^m \ra \Omega^{n+m}$
$$
(\omega_1\cdots \omega_n, \eta_1\cdots \eta_n) \mapsto
\omega_1\cdots \omega_n \cdot \eta_1 \cdots \eta_n
$$
It is not difficult to check that $\Omega^n$ defined this 
way is isomorphic to $A\otimes\tp{\bar{A}}{n}$. 
\eremark
The algebra $\Omega A$ has the following universal property. 
\bproposition
Let $R$ be a differential graded algebra $R = \bigoplus _{n }R_n$ and let $u: A\ra R_0$
be an algebra homomorphism. Then there exists a unique homomorphism 
$u_*: \Omega \ra R$ of differential graded algebras which extends $u$. 
\eproposition
Proof. It is easy to check that the map 
$$
u_*(a_0da_1\cdots da_n) = ua_0d(ua_1)\cdots d(ua_n)
$$
is a DG algebra homomorphism with the required property. \qed

\bproposition\label{trivial}
The homology of $\Omega $ with respect to the differential $d$ is
trivial in positive degrees and equals $\bc$ in degree $0$. 
\eproposition
Proof. It is clear that there is the following exact sequence
\beqn\label{Omega-short}
0 \ra d\Omega^{n-1} \arr{s} \Omega^n \arr{} d\Omega^{n} \ra 0 
\eeqn
for $n\geq 1$, where $s(\omega) = 1\otimes \omega$. 
In degree zero there is a surjection $d: A \ra dA\simeq 
\bar{A}$, whose kernel is $\bc$.   \qed 

As a first interesting homology theory that can be associated 
with an algebra $A$ we introduce the {\em Hochschild homology} 
$H_\star(A, M)$  of   $A$ with coefficients in an  $A$-bimodule
$M$. The quickest way to define $H_n(A,M)$ is to say that the $n$-th 
homology group is the $n$-th left derived functor of the commutator 
quotient space functor 
$$
M_\natural  =  M/[A,M]
$$
Following Quillen \cite{Q:Extensions} it is sometimes convenient to use 
the following suggestive notation to denote the commutator quotient space: 
$M\otimes_A = M/[A,M]$. Let us now explain what this means. 

To calculate Hochschild homology one uses the following standard 
resolution of the algebra $A$ by  free $A$-bimodules: 
$$
\arr{b'} A\otimes \bar{A}^{\otimes 2}\otimes A 
\arr{b'} A \otimes \bar{A} \otimes A \arr{b'} A\otimes A \arr{b'} A \ra 0
$$
where
$$
b'(a_0, \dots , a_{n+1}) = \sum_{i=0 } ^n (-1)^i(a_0, \dots, a_ia_{i+1}, 
\dots , a_{n+1})
$$
One checks that $(b')^2 = 0 $. 
Alternatively, one may write the same resolution in terms of 
$\Omega$ as follows. Using the isomorphism $\Omega ^n = A\otimes \bar{A}^n$, 
we have
\beqn\label{standard}
\arr{b'} \Omega^2 \otimes A \arr{b'} \Omega^{1}\otimes A \arr{b'}
A\otimes A \arr{b'} A \ra 0 
\eeqn
where the differential $b'$ is now given by the following formula
\beqn\label{b'-differential}
b'(\omega da \otimes a' ) = (-1)^{\deg(\omega)}(\omega a\otimes a'
- \omega \otimes aa')
\eeqn
and  in degree zero we have $b'(a\otimes a') = aa'$. 
The above complex is contractible, the contracting homotopy being 
provided by the map $s$ introduced in Proposition \ref{trivial}.

If we  apply the commutator quotient space functor  to the 
complex $\Omega\otimes A$ of (\ref{standard}) and use the 
canonical isomorphism 
$$
(\Omega\otimes A)\otimes_A \simeq \Omega; \qquad \omega \otimes a
\mapsto a\omega
$$
(this isomorphism holds for any $A$-bimodule), we obtain the complex
$\Omega = \bigoplus_n\Omega^n$ equipped with the differential
$$
\begin{array}{rcl}
b(a_0, \dots, a_n) & = & 
\displaystyle{\sum_{j=0}^{n-1} } (a_1, \dots, a_ia_{i+1}, \dots, a_n)\\
 && \mbox{} + (-1)^n(a_na_1, \dots , 
a_{n-1})
\end{array}
$$
Equivalently, this differential may be defined using the following 
compact formula
$$
b(\omega da) = (-1)^{\deg(\omega)}[\omega , a]$$
and $b(a) = 0 $. 
It is now easy to check that $b^2 = 0 $: 
$$\begin{array}{rcl}
b^2(\omega dada') &=& (-1)^{|\omega| +1} b (\omega da a' - a' \omega da) \\
&= & (-1)^{|\omega|+1} b(\omega d(aa') - \omega a da'- a'\omega da)\\
& = &\mbox{} - (\omega a a' - aa'\omega - \omega aa'+ a'\omega a - a'\omega a
+ aa' \omega)\\
&=& 0
\end{array}
$$
\bdefinition\label{Hochschild-homology}
The Hochschild homology $HH(A)$  of the algebra $A$ is the homology of
$\Omega $ with respect to the differential $b$: 
$$
H_n(\Omega A, b) = HH_n(A)
$$
\edefinition

We shall now define the {\em Karoubi operator} and state its basic properties. 
 \cite[Lemme 2.12]{Karoubi}\cite{CQ2,Kastler}.
\bdefinition
 If $\omega\in \Omega A$ we put
$$
\kappa (\omega da) = (-1)^{\deg(\omega)}da \omega
$$
\edefinition
Equivalently, we may write the above formula in the following form. 
$$
\kappa(\astr{0}{n})= (-1)^n(a_n, \astr{0}{n-1})+ (-1)^{n-1}(1, a_n
\astr{0}{n-1})
$$
\bremark\label{Karoubi-as-cyclic}
Using the definition of $\Omega^n A$ given in Remark \ref{multiple-tensors}
we can define the Karoubi operator as the generator of the action 
of the cyclic group on $\Omega^n$:
$$
\kappa: \omega_1\cdots \omega_n\mapsto (-1)^n \omega_n\cdot\omega_1\cdots 
\omega_{n-1}
$$
\eremark
The Karoubi operator is linked to the differentials $b$ and $d$ by 
the following formula. 
\blemma\label{kappa-homotopy}
$$
bd +db = 1-\kappa
$$
\elemma
Proof. Using the definition of $b$ and $\kappa$ we may write
$$
\begin{array}{rcl}
(bd + db)(\omega\, da) & = & b (d\omega da) + (-1)^{|\omega|}d [\omega, a]\\
& = & (-1)^{|\omega | +1}[ d\omega, a ] + (-1)^{|\omega|} d [\omega , a]\\
& = & [\omega , da] = \omega da - (-1)^{|\omega |}da \omega
\end{array}
$$
which gives the required formula. \qed

It now follows easily that $\kappa$ commutes with both $b$ and $d$, 
e.g.
$$
b\kappa  =  b(1 - (bd + db)) = b - bdb = \kappa b
$$
Let us define
 $$
B= \sum_{j=0}^n \kappa ^j d
$$ as an operator $B:\Omega ^n \ra \Omega^{n+1}$. 
It  follows that $dB= Bd = B^2= 0$ and that
$B$ commutes with $\kappa$ 
$$
B\kappa = \kappa B= B
$$
Moreover, it is not difficult to see that  $bB + Bb = 0 $; this 
fact derives from 
the following  lemma \cite{CQ2}.
\blemma\label{kappa-n}
The following identities hold in the space $\Omega^n$:
$$
\begin{array}{rcl}
\kappa ^n & = & 1+b\kappa^{-1}d\\
\kappa^{n+1}&=& 1-db\\
\kappa^{n(n+1)} & = & 1-Bb= 1+bB
\end{array}
$$
\elemma
Proof. We note that for any $i$, $0\leq i \leq n $
$$
\kappa^i(a_0 da_1\dots da_n) = (-1)^{i(n-1)}da_{n-i+1}\dots da_n a_0
da_1\dots da_{n-i}
$$
which gives in particular 
$$\kappa ^n (a_0 da_1\dots da_n) = 
da_1\cdots da_n a_0.
$$
This we can rewrite as follows
$$\begin{array}{rcl}
\kappa^n(a_0da_1\dots a_n)& = & a_0 da_1\dots da_n + [da_1\dots da_n, a_0]\\
& = & a_0 da_1\dots da_n + (-1)^n b(da_1\dots da_n da_0)
\end{array}
$$
which yields  the first formula. The second formula follows from
$$
\kappa^{n+1} = \kappa (1+b\kappa ^{-1}d)= \kappa +  bd= 1-db
$$
 when we use
the first identity of the lemma. Finally we have
$$
\kappa ^{n(n+1)}-1 = \sum_{i=0}^n \kappa^{ni} (\kappa^ n - 1) =
\sum_{i=0}^ n b\kappa ^{in-1}d = bB
$$
where we use the fact that $\kappa$ is an operator of order $n +1$ on
$d\Omega^n$.
Then we note that we can write the
same expression in the following way
$$
\kappa^{n(n+1)} -1 = \sum _{i=0}^{n-1}\kappa^{(n+1)i}(\kappa^{n+1}-1) =
- \sum_{i=0}^{n-1} \kappa^{(n+1)i} db = - Bb
$$
We thus get the third formula of the lemma. \qed

We note that the last formula in Lemma \ref{kappa-n} implies 
that $bB + Bb = 0 $ which indicates that $(\Omega , b, B)$ may 
be organized into the following double complex. 

\beqn\label{double}
\begin{array}{ccccc}
 \Bdar{b}&& \Bdar{b}&& \Bdar{b}\\
 \Omega^2 &\Lla{B}& \Omega^1 &\Lla{B}& \Omega^0\\
 \Bdar{b}&&\Bdar{b}\\
\Omega^1 &\Lla{B}& \Omega^0\\
 \Bdar{b}\\
\Omega^0 
\end{array}
\eeqn
where $\Omega^0 = A$, and where in bidegree $(p,q)$ we find 
$\Omega^{p-q}$. 

Cuntz and Quillen noticed \cite{CQ2} that in 
 analogy with  the de Rham complex on a Riemannian manifold, 
the two differentials $b,d$, of degree $-1$ ($+1$, respectively), 
on the graded algebra  $\Omega $ 
may be regarded as formal adjoints of each other, 
with $1 - \kappa$ playing the role of the Laplacian. It is 
therefore natural to study the spectral decomposition of $\Omega$ into 
eigenspaces of the Karoubi operator $\kappa$, paying special attention
to the subspace of \lq harmonic' forms, i.e. the space where
$\kappa = 1$. 

\blemma 
$$
(\kappa^n - 1 )(\kappa^{n+1} - 1) = 0 
$$
on $\Omega^n$. Consequently, there is a direct sum decomposition
$$
\Omega = \ker(\kappa - 1)^2 \oplus \bigoplus_{\xi} \ker(\kappa - \xi)
$$
where $\xi$ denotes all roots of unity different from $1$. 
\elemma
Proof. The first result  follows directly from identities in Lemma 
\ref{kappa-n}. Since $(x - 1)^2$ is the largest power of $(x-1)$ 
dividing $(x^n  - 1)(x^{n+1} - 1)$, the generalized eigenspace for $\kappa$ with 
eigenvalue $1$ is the kernel of $(\kappa - 1)^2$. 
Moreover, it follows from the first formula of this lemma that there 
is an eigenspace $\ker(\kappa - \xi)$ associated to 
each root of unity $\xi \neq  1$. 
\qed

Let us denote by $P$ the operator of spectral projection onto 
the generalized eigenspace where $\kappa$ has an eigenvalue $1$. 
We then have the following {\em harmonic decomposition}
$$
\Omega = P\Omega \oplus (1- P ) \Omega
$$
of the space of differential forms $\Omega$. 
Note that $1 - \kappa$ is invertible on $(1 - P)\Omega$. 
This decomposition is preserved by any operator commuting with $\kappa$. 
Thus, for example, it it invariant under the action of the operators 
$b$, $d$, $B$.
Moreover, one can  introduce the following analogue of the \lq Green's operator'
$G$ for $1 - \kappa$, which is defined to be the inverse of $1 - \kappa$
on $(1 - P)\Omega$, and $0$ on $P\Omega$. Then $PG = GP = 0$, 
and $G$ commutes with the operators commuting with $\kappa$. 

Let us denote $P^\perp = 1 - P$. We have the following \cite[2.1]{CQ2}. 
\bproposition 
The complex $P^\perp\Omega$ is contractible with respect to both $b$ and $d$. 
Consequently, 
the complexes $\Omega $ and $P\Omega$ have the same homology 
with respect to differentials $b$ and $d$. 
\eproposition
Proof. On $P^\perp\Omega$ one has the identity $1 = G(bd + db)$. Since $G$ commutes
with $b$ and $d$ we have that $(Gdb)d = G(bd + db)d = d$. Hence $Gbd$ is 
a projector with image $d P^\perp\Omega$, and similarly, $Gbd$ projects onto 
$b P^\perp\Omega$. 
Since the sum of these projectors is $1$, we have the decomposition 
$P^\perp\Omega = dP^\perp\Omega \oplus bP^\perp\Omega$. 
Furthermore, we have that  $(Gd)b = Gdb$ and $b(Gd) = Gbd$ which indicates that 
$b$ and $Gd$ are inverse isomorphisms between $dP^\perp\Omega$ and $bP^\perp\Omega$. 
In the same way we show that $d$ maps $bP^\perp\Omega$ isomorphically onto 
$dP^\perp\Omega$ with inverse $Gb$.  The statement of the proposition is 
now clear. 
\qed

\subsubsection{Cyclic homology}
We shall now devote some attention to the properties of the total complex 
$(\calb(\Omega), b+B)$ of the double complex (\ref{double}).
\bdefinition
The cyclic homology of the algebra $A$ is the homology of the total complex 
$\calb \Omega$. Hence, in degree $n$, we have:
$$
HC_n(A) = H_n(\calb(\Omega A), b + B)
$$
\edefinition
\bexample
The simplest example of this situation is the case when  $A=\bc$ is 
the algebra of complex numbers. Then $\Omega^0 \bc = \bc$  and 
$\Omega^n \bc = 0 $ for $n \geq  1$. Thus the double complex 
(\ref{double}) has a particularly simple form
$$
\begin{array}{ccccc}
 && \Bdar{0}&& \Bdar{0}\\
 0& \Lla{0}& 0  &\Lla{0}& \bc\\
 \Bdar{0}&&\Bdar{0}\\
 0  &\Lla{0}& \bc\\
 \Bdar{0}\\
\bc
\end{array}
$$
It follows that  the total complex $\calb(\Omega\bc)$ is 
$$
\bc \larr{0} 0 \larr{0} \bc \larr{0} 0 \larr{} \dots 
$$
with $\bc $ in even degrees and $0$ in odd degrees. Thus
$HC_{2k}(\bc) = \bc  $ and $HC_{2k+1} (\bc) = 0$.   
\eexample

The total complex $\calb(\Omega)$ is equipped with the following 
{\em Connes periodicity operator} $S$. Elements of $\calb(\Omega)$ 
may be described as polynomials  $\sum_{p\geq 0 } u^p\omega_{n-2p}$, 
where $\omega _{n- 2p} \in 
\Omega^{n-2p}$. Then we put $S(u^0\omega_k) = 0 $ and
$S(u^p\omega_k) = u^{p-1}\omega_k$ for $p \geq 1$, which is illustrated 
by the following diagram

\beqn\label{S-complex}
\begin{array}{ccccc}
 \Bdar{b}&
\begin{picture}(30,20)(0,0)
\put(35,20){\vector(-3,-2){40}}
\put(11,10){$\scriptstyle{S}$}
\end{picture} 
& \Bdar{b}&
\begin{picture}(30,20)(0,0)
\put(35,20){\vector(-3,-2){40}}
\put(11,10){$\scriptstyle{S}$}
\end{picture} & \Bdar{b}\\
 \Omega^2 &\Lla{B}& \Omega^1 &\Lla{B}& \Omega^0\\
 \Bdar{b}& 
\begin{picture}(30,20)(0,0)
\put(35,20){\vector(-3,-2){40}}
\put(11,10){$\scriptstyle{S}$}
\end{picture} 
&\Bdar{b}& \begin{picture}(30,20)(0,0)
\put(35,20){\vector(-3,-2){40}}
\put(11,10){$\scriptstyle{S}$}
\end{picture} \\
\Omega^1 &\Lla{B}& \Omega^0\\
 \Bdar{b}& \begin{picture}(30,20)(0,0)
\put(35,20){\vector(-3,-2){40}}
\put(11,10){$\scriptstyle{S}$}
\end{picture} \\
\Omega^0 
\end{array}
\eeqn
It is clear from this diagram, that there is the following 
exact sequence of complexes
\beqn\label{S-seqn}
0 \ra \Omega \arr{I} \calb(\Omega) \arr{S} \calb(\Omega)[2] \ra 0
\eeqn
where the map $I$ is the inclusion of the column $p=0$.
Furthermore, $\calb(\Omega)[2]$ denotes
the {\em suspension} of the  complex $\calb(\Omega)$. Recall
that for any complex $X = \bigoplus_{n\in \bz}X_n$ equipped with 
differential $d$ of degree $-1$ the $k$-th suspension is the
complex $X[k]$ defined by 
$$
X[k]_n = X_{n-k}, \qquad d_{X[k]} = (-1)^{k}d_X
$$
With these conventions in  mind we see that the sequence 
(\ref{S-seqn}) leads to the following long exact sequence 
of homology groups first derived by Connes in the dual case 
of cohomology \cite{Connes:IHES,LQ}. 
$$
\ra HH_{n+2}(A) \arr{I} HC_{n+2}(A) \arr{S}
HC_n(A) \arr{B} HH_{n+1}(A) \ra 
$$
Here the connecting map $B$ is induced by 
$\sum_{p\geq 0}u^p\omega_{n-2p}\mapsto B\omega_n $. 

Finally, we define the {\em periodic cyclic homology} to be the 
homology of the complex 
$$
HP_\nu(A) = H_\nu(\prod_n\Omega_n, b + B)
$$
where $\nu \in \bz_2$. Thus the periodic cyclic homology 
is a $\bz_2$-graded homology theory associated with the 
following periodic complex
$$
\arr{} \hat{\Omega}^{\rm even} \arr{b + B} \hat{\Omega}^{\rm odd}\arr{b+B}
\hat{\Omega}^{\rm even} \arr{}
$$
where $\hat{\Omega}^{\rm even} = \prod_{k\geq 0 }\Omega^{2k}$, and 
similarly in the odd case. 
\bexample
Take the algebra of polynomials in one variable $P$ and consider
also the algebra $Q$ of Laurent polynomials in one variable $z$. 
Then in both cases, the bimodule of forms $\Omega ^1$ is of dimension
one over $P$, that is it is generated by one element $dz$. 
We want to calculate its periodic cyclic cohomology. First, the 
Hochschild-Kostant-Rosenberg theorem states that there is 
a quasi-isomorphism of mixed complexes
$$
(\Omega A, b, B)\ra (\Omega_A, 0, d)
$$
where $A= P,Q$, $\Omega_A$ is the de Rham complex of the algebra $A$, 
and $d$ is the usual de Rham differential. In our case, $\Omega_A$ 
is zero in degrees higher than one.

To calculate the periodic cyclic cohomology we note that the total 
complex of the 
mixed complex $(\hat{\Omega}_A, 0, d)$ may be arranged as the 
following double complex
$$
\arr{0} P \arr{d} P\arr{0} P \arr{d} P \arr{}
$$
where the domain of $d$ is in the even degrees and its range in the odd. 
Here we have used that $\Omega^1 P \simeq P$. Then $\ker d: P\ra P$ 
is isomorphic to $\bc$, whereas its image coincides with $P$, and integration
gives the inverse map. Thus
\blemma
$$
HP^*(P) = \left\{\begin{array}{cl}
\bc & \mbox{\rm when}\; * = 2\nu\\
0& \mbox{\rm otherwise}
\end{array}
\right.
$$
\elemma
In the case of the algebra of Laurent polynomials $Q$, the map $d:Q\ra Q$
from even degrees to odd is not surjective: the one-dimensional 
space of monomials  $a_{-1}z^{-1}dz$
is not in its image. Therefore in this case we have a different result:
\blemma
$$HP^*(Q) = \left\{\begin{array}{cl}
\bc & \mbox{\rm when}\; * = 2\nu\\
\bc & \mbox{\rm otherwise}
\end{array}
\right.
$$
\elemma
It is worth noting that the isomorphism is induced by the evaluation
map $\epsilon_0$ in the even degrees  and by the residue map 
in the odd degrees. 
\eexample
\bremark\label{b'-complex}
The double complex (\ref{double}) has the following origin. 
We define an action of the cyclic group of order $n$ on the space
$\tp{A}{n}$ by letting the generator act as
$$
\lambda(\astr{1}{n})= (-1)^{n-1}(a_{n}, \astr{1}{n-1})
$$
We denote by $N$ the corresponding  cyclic norm operator
$N=
\sum_{i=0}^{n-1}\lambda^i$. The four operators $b$, $b'$, $\lambda$ and $N$
are related as follows \cite{T}.
\blemma
One has
\beqn\label{tsygan}
b(\lambda -1)= (\lambda-1) b', \qquad b'N= Nb
\eeqn
\elemma
Proof. We check
that
$$
b'= \sum_{j=1}^n\lambda^{j-1}c\lambda ^{-j}\qquad
b=\sum_{j=1}^{n+1}\lambda^{j-1} c \lambda ^{-j}
$$
where
$$c(\astr{1}{n})= (-1)^{n-1} (a_na_1, \dots,
a_{n-1})
$$
is the crossover term in the formula for $b$.These identities
imply  (\ref{tsygan}) \qed

It follows from (\ref{tsygan}) that we can organize the two
complexes with differentials $b$ and $b'$, respectively, into the
following double complex which we shall denote by $\bc_{\bullet\bullet}(A)$.

\beqn\label{doublecx}
\begin{array}{cccccc}
\Bdar{b}&&\Bdar{-b'}&&\Bdar{b}\\
\tp{A}{2}& \Lla{1-\lambda} & \tp{A}{2}& \Lla{N} & \tp{A}{2}&
\Lla{1-\lambda}\\
\Bdar{b}&&\Bdar{-b'}&&\Bdar{b}\\
 A& \Lla{1-\lambda} & A& \Lla{N} & A &
\Lla{1-\lambda}
\end{array}
\eeqn
Columns in  the even degrees are copies of the  Hochschild
complex 
$$
C(A)= (\bigoplus_{n\geq 0 } \tp{A}{n+1}, b)
$$
and the odd degree columns are (when $A$ is unital) the acyclic Hochschild
complexes $B(A) = (\bigoplus_{n\geq 0 } \tp{A}{n+1}, b')$. 
The rows in the above complex are exact.
(This is true if $A$ is an algebra over a field $k$ of characteristic zero).
Loday and Quillen have demonstrated that  the total complexes 
of ${\bf C}_{\bullet\bullet}(A)$ and of $\calb(\Omega )$ are quasi-isomorphic, 
by which we mean the following. There is a canonical surjection
${\rm Tot}{\bf C}_{\bullet\bullet}(A) \ra \calb(\Omega )$ 
induced by the projection $\tp{A}{n}\ra \Omega^n$ described previously. 
This  surjection  induces an isomorphism $H_\star({\rm Tot}{\bf C}_{\bullet\bullet}(A))
\simeq H_\star(\calb(\Omega))$. This statement is known as the 
{\em simplicial normalization theorem} \cite{LQ,Kassel:Crelle}

\bdefinition
The cokernel of the map $1-\lambda$ from the $b'$-complex to the $b$-complex
in the zeroth column of the above complex
is the cyclic complex $C^\lambda_{\bullet}(A)$. It is equipped with the 
differential $b$. 
\edefinition

Thus the cyclic complex $C^\lambda(A)$ may be defined as the 
cokernel of the map $(1-\lambda)$ in the following augmented complex
$$%\label{augmented}
\begin{array}{cccccccccc}
&&\Bdar{b}&& \Bdar{b}&&\Bdar{-b'}&&\Bdar{b}\\
0& \larr{}& C_1^\lambda(A)&\larr{}& \tp{A}{2}& \Lla{1-\lambda} & \tp{A}{2}& \Lla{N} &
 \tp{A}{2}& \Lla{1-\lambda}\\
&& \Bdar{b}&& \Bdar{b}&&\Bdar{-b'}&&\Bdar{b}\\
0 &\larr{} &C_0^\lambda(A)&\larr{}& A& \Lla{1-\lambda} & A& \Lla{N} & A &
\Lla{1-\lambda}
\end{array}
$$

Using the periodicity of the double complex 
$\bc_{\bullet\bullet}(A)$ we check that the
following short sequences of complexes are exact
$$
\begin{array}{c}
0\larr{} C^\lambda(A)\larr{} C(A)\larr{} \ig (1-\lambda)\larr{} 0\\
0\larr{} \ig(1-\lambda) \larr{} B(A)\larr{} C^\lambda(A)\larr{} 0
\end{array}
$$
From the first sequence we have the following long exact sequence of homology groups
$$
\ra  HC_{n+1}(A) \ra H_n(\ig (1-\lambda))
\ra HH_n(A) \ra HC_n(A) \ra
$$
whereas, when $A$ is unital, the second sequence yields
$$
H_n(\ig(1-\lambda))= HC_{n-1}(A)
$$
since in this case the $b'$-complex is acyclic. The two facts together
give once again Connes's long exact homology sequence. 
$$
\lrar HH_n(A)\lrar HC_n(A)\arr{S}HC_{n-2}(A) \lrar HH_{n-1}(A)
\lrar
$$
The connecting homomorphism in the sequence is the Connes $S$-operator.
\eremark
\subsubsection{The nonunital case}
When the algebra $A$ does not have a unit, a slight 
modification of the formalism is required, pretty much following
the same lines as when we defined $K$-groups for nonunital rings.
Let $\tilde{A}$ be the augmented algebra $\tilde{A}  = \bc \semiright A$, which 
is characterized by the existence of the {\em augmentation} map, i.e.
an algebra homomorphism $\epsilon: \tilde{A} \ra \bc$. 
This map induces an augmentation $\Omega\tilde{A} \ra \Omega\bc = \bc$, which 
makes $\Omega \tilde{A}$ into an augmented differential graded algebra
\beqn\label{Omega-augment}
\Omega\tilde{A} = \bc \oplus \bar{\Omega}\tilde{A}
\eeqn
The augmentation ideal $\bar{\Omega}\tilde{A}$ in this algebra is the 
universal nonunital differential graded algebra generated by the nonunital 
algebra $A$. This is precisely the universal algebra used by Connes \cite{Connes:IHES,Connes:Book} (see also \cite{Arveson}); this algebra is 
also discussed under the name of differential envelope by Kastler \cite{Kastler}, 
and by Karoubi \cite{Karoubi}. We remark that $\bar{\Omega}\tilde{A}$ 
is closed under all the operators defined on $\Omega\tilde{A}$. 

To describe the algebras $\Omega\tilde{A}$ and $\bar{\Omega}\tilde{A}$ in 
more detail, let us note that in degree $n$ we have that $\Omega^n\tilde{A} = 
\tilde{A}\otimes \tp{A}{n}$
which gives the following isomorphism 
$$
C_n(A) \oplus C_{n-1}(A) \arr{\sim} \bar{\Omega}^n\tilde{A}
$$
Here $C_n(A)= \tp{A}{n+1}$ denotes as in \cite[\S 5]{CQ2} the space of unnormalized 
$n$-chains.  
This isomorphism allows us to represent elements of $\bar{\Omega}\tilde{A}$ 
as column vectors with components in $C$ and operators on $\bar{\Omega}
\tilde{A}$ are described by $2\times 2$-matrices of operators on $C$. It is 
not difficult to check that the operators $\tilde{d}$, $\tilde{b}$, $\tilde{B}$, $\tilde{\kappa}$
are given by the following expressions
$$
\begin{array}{ccc}
\tilde{d} = \left( \begin{array}{cc} 0 & 0 \\ 1 & 0 \end{array}
\right),  && \tilde{b} = \left(\begin{array}{cc} b & 1 - \lambda \\
0 & -b'\end{array}\right),  \\
\rule{0mm}{10mm}
\tilde{\kappa} = \left(\begin{array}{cc}
\lambda & 0 \\ b' - b & \lambda
\end{array}\right), &&
\tilde{B} = \left( \begin{array}{cc}  0 & 0 \\ N_\lambda & 0 \end{array}\right),
\end{array}
$$
where  $N_\lambda = \sum_{i = 1}^{n-1} \lambda ^i$ in degree $n$. 
Furthermore, the identities $\tilde{b}^2 = 0 $ and $\tilde{b}\tilde{B} 
+ \tilde{B}\tilde{b} = 0 $ on $\bar{\Omega}\tilde{A}$ are equivalent
to the Connes-Tsygan identities (\ref{tsygan}). It is now obvious that 
$\bar{\Omega}\tilde{A}$ corresponds to the first two columns of the 
double complex \ref{doublecx}. 

\subsection{Deformations of the algebra $\Omega A$}
\subsubsection{Universal algebra extensions}\label{RA}
The algebra $\Omega$ of noncommutative differential forms 
is just one example of a universal  algebra 
that can be associated  with a given algebra $A$. We shall now 
introduce two important deformations of $\Omega$. First of all
we define  the algebra $RA$ which is important in 
the study of algebra extensions. Secondly, in the next section, 
we shall discuss the Cuntz algebra $QA$, which is sometimes 
called the algebra of quantized differential forms over $A$. 

\bdefinition \label{based-linear}
Let $A$ and $B$  be  unital algebras. A {\em based linear\/} map is a linear map 
$\rho : A \ra B$ that sends the unit of $A$ to the unit of $B$. 
\edefinition

Let $\rho$ be a based linear map. 
We define its {\em curvature} to be the bilinear map
\beqn\label{curvature-of-rho}
\omega(a_1, a_2) = \rho (a_1a_2)  - \rho(a_1)\rho(a_2)
\eeqn
In other words, $\omega$ measures the deviation of $\rho$ from an algebra 
homomorphism. Since $\omega$ vanishes when  either of its arguments is $1$, 
it can be viewed as a linear map $\omega : (\tp{\bar{A}}{2}= d\Omega^1 A) \ra B$. 

\bdefinition 
Let $T(A) = \bigoplus_{n\geq 0} \tp{A}{n}$ be the tensor algebra of the 
underlying vector space of $A$. Then we put \cite[\S 1]{CQ1}
$$
RA = T(A) /T(A) ( 1_T - 1_A)T(A)
$$
where $1_T$ is the identity of $T(A)$, and $1_A$ is the identity of $A$, regarded as 
a tensor of degree one. 
\edefinition
Let $\hat{\rho} : A \ra RA $ be the map defined as the composition of  the inclusion 
$A\ra T(A)$ with the canonical surjection $T(A) \ra RA$. It is clear that
 $\hat{\rho}$
is a based linear map. The algebra $RA$ has the following universal 
property. Let $R$ be a unital algebra and let $\rho : A \ra R$ be a based linear map. 
Then there exists a unique algebra homomorphism $ RA \ra R$ which sends 
$\hat{\rho}$ to $\rho$. 

If we apply the universal property to the identity map $A \ra A$ we see that there 
is a unique homomorphism $RA \ra A$. Let $IA$ be its kernel. This way 
we obtain an algebra extension $A = RA/ IA$ for which $\hat{\rho} $ is a 
based linear lifting, i.e. a based linear map which is a section of the canonical 
surjection $RA \ra A$.  $RA$ is called the universal extension of $A$. 

To describe the algebraic structure of $RA$ in a way that reflects the product
on $A$, we need first to introduce a certain deformation of the product in the 
algebra $\Omega A$ of differential forms. This product has been discovered 
by Fedosov \cite{F} in his work on index theorems.  If $\omega $ and $\eta $ 
are homogeneous elements of $\Omega A$ then 
\beqn\label{Fedosov}
\omega\circ \eta = \omega\eta - (-1)^{|\omega|} d\omega d\eta
\eeqn
where we denote by $|\omega| = \deg(\omega)$. This product is associative
and compatible with the even-odd grading on $\Omega A$, and so makes 
$\Omega A$ into a superalgebra. We note that this product is not commutative
in general. The following proposition establishes a link between the 
universal extension $RA$ and the algebra $\Omega^{ev}A$ of 
even differential forms equipped with the Fedosov product \cite[Propn. 1.2]{CQ1}.
\bproposition
There is a canonical isomorphism $RA \simeq \Omega^{ev}A$, where the 
algebra $\Omega^{ev}A$ of even forms is equipped with the Fedosov product. 
Explicitly, this isomorphism is given by the map 
$$
\hat{\rho}(a_0)\hat{\omega}(a_1,a_2)\cdots \hat{\omega}(a_{2n-1}, a_{2n}) 
\leftrightarrow a_0 da_1 \cdots da_{2n}
$$
Under this isomorphism,  the ideal $IA^n$ corresponds to 
$\bigoplus _{k\geq n} \Omega^{2k} A$. 
The associated graded algebra ${\rm gr} _{IA}RA = \bigoplus IA^n / IA^{n+1}$ is isomorphic
to the algebra $\Omega^{ev} A$ of even forms with the ordinary product. 
\eproposition
Proof. (Sketch) Let us denote by $\theta : A \ra \Omega^{ev}A$ the 
canonical inclusion, which is a based linear map. Note that $\theta$ is not 
an algebra homomorphism. The curvature of $\theta$ is given by 
$$
\theta(a_1a_2) - \theta(a_1)\circ \theta (a_2) = a_1a_2 - a_1\circ a_2= da_1 da_2
$$
Let $\Psi: RA \ra \Omega^{ev}A$ be the corresponding morphism given by 
the universal property of $RA$. Then $\Psi(\hat{\rho}(a))= \theta(a) = a$, 
and 
$$
\Phi(\hat{\omega}(a_1,a_2)) = \Psi(\hat{\rho}(a_1a_2) - \hat{\rho}(a_1)\hat{\rho}(a_2) ) 
= da_1 da_2
$$
from which it follows that 
$$
\Phi(\hat{\rho}(a_0)\hat{\omega}(a_1,a_2)\cdots \hat{\omega}(a_{2n-1}, a_{2n})) = 
a_0da_1 \cdots da_{2n}
$$
where we use that the Fedosov product of closed forms coincides with the usual 
product. On the other hand there is a well defined map $\Phi: \Omega^{ev}A \ra RA$
$$
\Phi: a_0da_1\cdots da_{2n} \mapsto \hat{\rho}(a_0)\hat{\omega}(a_1,a_2)\cdots \hat{\omega}(a_{2n-1}, a_{2n})
$$
We check that $\Phi$ is surjective and that $\Phi$ and $\Psi$ are inverses
of each other. 
\qed

If we identify $RA$ with the algebra of even forms as in the above proposition, 
then the universal property of $RA$ may be stated as follows. For any based linear 
map $\rho: A\ra B$ there exists a unique algebra homomorphism $\rho_\star: RA\ra B$
given by 
$$
\rho_\star(a_0da_1\cdots da_{2n}) = 
\rho(a_0)\omega(a_1,a_2)\cdots \omega(a_{2n-1},a_{2n})
$$
where $\omega $ is the curvature of $\rho$.  Moreover, 
using the same identification
we see that $IA = \bigoplus_{k\geq 1}\Omega^{2k}A$. 

\mbox{}

Let $R$ be an algebra. We now want to describe briefly the following 
interesting supercomplex \cite[\S 4]{CQ3} which we shall denote 
$X(R)$: 
$$
X(R): \qquad \Arr{\natural d} R \Arr{\beta} 
\Omega^1 R_\natural \Arr{\natural d} R \Arr{\beta} 
$$
Here $\Omega ^1R_\natural = \Omega^1R/ [\Omega^1 R, R] = 
\Omega^1 R/b\Omega^2 R$ and $\natural : \Omega^1R \ra \Omega^1R _\natural$
denotes the canonical projection. Moreover, 
$ \beta(\natural(xdy)) = b(xdy) = [x,y]$. The $X$-complex has been 
defined and studied in detail in the case of free algebras and coalgebras by Quillen
\cite{Q:Cochains} and by Cuntz and Quillen in general \cite{CQ1,CQ3}. The importance 
of this complex lies in its relation to the periodic cyclic homology of an algebra. 
It may be regarded as the first order approximation of the  standard double complex 
\ref{double}. There are, however, important cases when the $X$-complex contains 
sufficient cohomological information. This complex may be regarded as a 
noncommutative analogue of the ordinary de Rham complex $\Omega_R$ for commutative
algebras. If $R$ is a smooth commutative algebra, then an extension of the 
Hochschild-Kostant-Rosenberg theorem \cite{HKR,LQ} states that this complex 
calculates the cyclic homology of the algebra. Moreover, there is a class 
of algebras called {\em quasi-free} for which the $X$-complex is sufficient 
to calculate their homology. An algebra is quasi-free if and only if its dimension 
(relative to Hochschild cohomology) is not greater than one. 

There is another important situation, when the $X$-complex comes into its own, 
namely when $R = RA$, and our immediate task will be to find an explicit 
description of $X(RA)$. We have already seen that we may identify $RA$ with the
algebra of even differential forms $\Omega^{\rm even}A$ equipped with the 
Fedosov product, and so we now need to find a similar description for the 
commutator quotient space $\Omega^1(RA)_\natural$.  To avoid confusion with the differential $d$ on $\Omega A$, 
we shall denote by $\delta$ the canonical derivation from $RA$ to $\Omega^1 (RA)$.
First we note that there is a canonical isomorphism of $RA$-bimodules
$$
RA\otimes \bar{A}\otimes RA \arr{\sim} \Omega^1(RA); \qquad \qquad x\otimes a\otimes y 
\mapsto x\delta a y
$$
from which it follows that there is an isomorphism  of commutator quotient spaces 
$$
RA \otimes \bar{A} \arr{\sim} \Omega^1(RA)_\natural; \qquad\qquad x\otimes a \mapsto
\natural(x\delta a)
$$
It is clear that we can write the last isomorphism as 
$$
\Omega^{od} A \arr{\sim} \Omega^1(RA) _\natural; \qquad\qquad xda \mapsto \natural (x\delta a)
$$
Since we have seen before that $RA\simeq \Omega^{ev}A$, we have that the $X$-complex 
of the algebra $RA$ may be identified with $\Omega A$ as a $\bz_2$-graded vector
space. It is not very difficult to identify the differentials $\beta$ and $\natural \delta$ in 
terms of this isomorphism \cite[\S 5, Eqn. (9), Lemma 5.3]{CQ3}. 
From our general definition of the $X$-complex follows that 
$$
\beta(\natural(x\delta y )) = x\circ y - y\circ x = xy-yx -dxdy+dydx = b(xdy)- (1+\kappa)d(xdy)
$$
which gives 
$$
\beta = b - (1+\kappa)d: \Omega^{od}A \arr{} \Omega^{ev}A
$$
The calculation of the other differential is slightly more involved, we get
in this case
$$
\natural \delta = - N_{\kappa^2}b + B : \Omega^{ev} A \arr{} \Omega^{od}A
$$
where $N_{\kappa^2} = \sum_{j=0}^{n-1}\kappa^{2j}b$ on $\Omega^{2n}A$. 

\mbox{}

Let $c$ be the scaling operator which is multiplication by $c_q$ on $\Omega^q A$, 
where $c_{2n} = c_{2n+1} = (-1)^n n!$.  Cuntz and Quillen proved 
the following important fact \cite[Thm. 6.2]{CQ3}. 
\btheorem
The maps 
$$
cP: X(RA) \arr{} \Omega A, \qquad\qquad c^{-1}P : \Omega A \arr{}X(RA),
$$
where $P$ is the projection onto the space of harmonic forms, 
are inverse modulo homotopy, so that $(X(RA), \beta \oplus\natural \delta)$ and 
$(\Omega A, b + B)$ are homotopy equivalent supercomplexes. 
\etheorem

The importance of this theorem stems from the fact that the $X$-complex
has additional functorial properties, which are not shared by the 
standard $b+B$-complex $\Omega A$, as we shall see in the following sections. 

\subsubsection{The Cuntz algebra $QA$}\label{Cuntz-algebra}

We now pass to the description of the Cuntz algebra $QA$ \cite{Cu2}
associated with an
algebra $A$. If $A$ is a unital algebra, then
the Cuntz algebra $QA$ is defined as  the  free
product $QA = A * A$ in the category of unital algebras, i.e., with
amalgamation over the identity of $A$. There are two canonical homomorphisms
$\iota, \iota^\gamma$ from $A$ to $QA$ and a canonical automorphism of
$QA$ of order two:  $\omega \mapsto \omega^\gamma$
which interchanges  $\iota $ and $\iota^\gamma$. The algebra $QA$
is a superalgebra which means that $QA$ is $\bz_2$-graded
and the grading is compatible with multiplication.

For $a\in A$, let $p(a)$ and $q(a)$ denote the even and odd components
of $\iota(a)$ with respect to the $\bz_2$-grading of $QA$. If we put
$a=\iota(a)$ and $a^\gamma = \iota^\gamma(a)$, then
$p(a) = (a + a^\gamma)/2$ and $q(a) = (a - a^\gamma)/2$.
The following relations  determine the product in $QA$.
$$
\begin{array}{rcl}
p(a_1 a_2) & = & p(a_1) p(a_2) + q(a_1) q(a_2) \\
q(a_1a_2) & = & p(a_1) q(a_2) + q(a_1) p(a_2)
\end{array}
$$
We deduce from these relations that any element in $QA$ may be
represented as a linear combination of elements of the form
$p(a_0) q(a_1) \cdots q(a_n)$.
Note finally that the unit of the algebra $A$ is also the unit of $QA$,
and that  we have $p(1) = 1$, $q(1) = 0$.

There is a canonical  map $QA \ra A$ which identifies the two copies
of the algebra $A$ inside $QA$. This map is  called the folding map, and its kernel is
denoted  by $qA$.   There
is a natural completion of this algebra  $\hat{Q}A = \displaystyle{\lim_{\llar}} \; QA /{qA}^n$ with respect 
to the $qA$-adic topology induced by the ideal $qA$. 
This completion is again a $\bz_2$-graded algebra.

There is a linear  isomorphism $\mu$ between the underlying graded
vector space of the Cuntz algebra $QA$ and the graded vector space of
forms $\Omega A$ which on generators is defined by
\beqn\label{mu-map}
\mu : p(a_0) q(a_1) \cdots q(a_n) \mapsto a_0 da_1 \cdots da_n
\eeqn
In particular, if $a\in A$ as a subalgebra of $QA$ then
$\mu(a) = \mu(p(a)) + \mu (q(a)) = a + da\in \Omega A$.

The  map $\mu$  becomes an isomorphism of  ${\bf Z}_2$-graded  algebras
when the algebra  $\Omega A$ is equipped with the Fedosov product (\ref{Fedosov}).
We may therefore think of the Cuntz algebra as a  deformation of
the algebra of
differential forms $\Omega A$. 
In particular, under this isomorphism, the ideal $J$ of forms of degree one or higher
 in $\Omega A$
is identified with the folding ideal $qA$ of the Cuntz algebra.

Using this isomorphism  it is easy to see that if $\tilde{A}$
is an augmented algebra then $Q\tilde{A}$ is also augmented
$$
Q\tilde{A} = \bc \oplus \bar{Q}\tilde{A}
$$
where $\bar{Q}\tilde{A}$, the reduced part of the Cuntz algebra,
is the Cuntz algebra associated with $A$ in the nonunital category.

If $A$ is in fact a unital algebra then the unit $1$ of $A$ has two images
$u$ and $u^\gamma$ in $Q\tilde{A}$ which are idempotents. The corresponding
involutions $f = 2u - 1$ and $f^\gamma $ generate inside the algebra
$\hat{Q}\tilde{A}$ an algebra which is isomorphic to $\hat{Q}\tilde{\bc}$, the group
algebra of the infinite dihedral group. We shall record its basic properties
for future use \cite{Cu2}.

The dihedral group generated by the two involutions
$f$ and $f^\gamma$ may be presented  using the following
two generators. First, let
$$
L= \log(ff^\gamma)= - \sum_{n\geq 1} (2fq(f))^n/n
$$
Since $(ff^\gamma)^\gamma = f^\gamma f = (ff^\gamma)^{-1}$,
we have on the one hand
$(e^L)^\gamma = e^{-L} $ while on the other $(e^L)^\gamma =
e^{(L^\gamma)}$
from which follows that $L^\gamma = - L$, i.e. $L$ is odd with respect
to the canonical automorphism $\gamma$.

We  put $W = e^{L/2}$. Since  $ff^\gamma = \exp L$ we
have
$$
Wf^\gamma = W f e^{L} = e^{L/2} e^{-L} f = e^{-L/2}f = fe^{L/2}
$$
so that
\beqn\label{conjugation}
f^\gamma = W^{-1} f W
\eeqn

It is now a good moment to state explicit
formulae for $L$ and $W$ that we shall need later \cite{ja23}.
\bproposition \label{L}
Let $L = \log(ff^\gamma) = -\sum_{n\geq 1} (2fq(f))^n /n$.
Then
$$
\log(ff^\gamma) = - \sum_{i\geq 1} 2^{2i-1}\frac{((i-1)!)^2}{(2i-1)!}p(f)q(f)^{2i-1}.
$$
\eproposition
Let us now define $W_t = \exp(tL/2)$ for $t\in [0, 1]$. \bproposition\label{Wt}
For all $t\in [0,1]$ we have
$$
W_t = 1 + \sum_{n\geq 1} (-1)^{n} 2^{2n-1}{t/2 + n - 1 \choose
2n-1} ((t/2n)q(f)^{2n} + p(f)q(f)^{2n-1})
$$
Moreover, we have that   $W_{-t} = W^\gamma_t$. Thus  $p(W_{-t})  = p(W_t)$
and $q(W_{-t}) = -q(W_t)$.
\eproposition
\subsubsection{The Cuntz algebra and $KK$-theory}
The original motivation for introducing the Cuntz algebra $QA$ 
was its relation to the $KK$-theory, which we now explain very briefly. 

Let $(\cale, F, \gamma)$ be an even Kasparov module as in Section \ref{KK}. 
We shall assume 
that $F$ is in fact a self-adjoint involution, as $F$ can always be modified
in this way \cite{Connes:IHES}. 
It is not difficult to check that such a Kasparov module induces 
a homomorphism $\alpha$ from $qA$ to $\bk$, the ideal of 
compact operators on $\cale$. 
This homomorphism is defined by setting $\alpha(\omega) = \frac{1+\gamma}{2}\omega$, 

$ \alpha^\gamma(\omega) = \frac{1+\gamma}{2} F\omega F$. 

First, we want to describe $K_0(A)$ in terms of the ideal $qA$. 
Note that there are canonical maps $\pi_0, \pi_1: qA \ra A$, which are 
restrictions of the natural maps ${\rm id}* 0, \; 0*{\rm id} : QA \ra A$. 
In particular, there is an induced map 
\beqn\label{K(qA)}
K_{0}(\pi_0)  : K_0(qA) \arr{} K_0(A)
\eeqn
It turns out that this map is in fact an isomorphism \cite[Propn. 3.2]{Cuntz:Elements}. 
\bproposition
The map \ref{K(qA)} is an isomorphism with the inverse 
$$
K_0(\imath) - K_0(\imath^\gamma) : K_0(A) \arr{} K_0(qA)\subset K_0(QA)
$$
Here $\imath, \imath^\gamma: A \ra QA$.
\eproposition

This proposition has the following interesting corollary. We shall 
denote by $[A,B]$ the set of homotopy classes of homomorphisms 
between two $C^*$-algebras $A$ and $B$. For two such algebras, 
we endow $[A, \bk \otimes B]$ with the following semigroup
structure: 
$$
[\phi] + [\psi] = \left[\larray \phi & 0 \\ 0 & \psi \rarray\right]
$$
where 
$$
\larray \phi & 0 \\ 0 & \psi \rarray : A \arr{} M_2(\bk \otimes B) \simeq \bk \otimes B
$$
is the homomorphism that maps $\omega $ to $\larray \phi(\omega) & 0 \\ 
0 & \psi(\omega) \rarray$. 
\bproposition
There is an isomorphism 
$$
[q\bc , \bk \otimes A] \arr{\sim} \bk \otimes A
$$
To define this map, let $e$ be the canonical  generator of 
$K_0(q\bc) \simeq K_0(\bc) \simeq \bz$. Then to  any $[\phi]\in [q\bc, \bk \otimes A]$
we assign the element $K_0(\phi)(e) \in K_0(\bk \otimes A)\simeq K_0(A)$. 
\eproposition
Moreover, it turns out that $[q\bc , \bk \otimes A]$ is an abelian group. 
If $\gamma$ is the restriction to $q\bc$ of the canonical involution 
on $Q\bc$ that exchanges the two copies of $\bc $ inside $Q\bc$, then 
the inverse element to $[\phi] $ is $[\phi \gamma]$. 

With these results at hand we may define $KK(A,B)$ as follows. 
\bdefinition
Given two $C^\star$-algebras $A$ and $B$ we put 
$$
KK(A,B) = [qA, \bk \otimes B]
$$
\edefinition

One can use this definition to demonstrate the existence of the 
Kasparov product
$$
KK(A_1, A_2) \times KK(A_2, A_3) \ra KK(A_1, A_3)
$$
for any $C^*$-algebras $A_1$, $A_2$ and $A_3$.

\subsection{Cochains and characters}
In this section we develop basic features of  the theory of cochains on an 
algebra in both $\bz$ and $\bz_2$-graded cases. 
A thorough discussion of cochain complexes in cyclic cohomology in terms of the 
bar construction on an algebra is
given in Quillen's paper \cite{Q:Cochains} and in the $\bz_2$-graded case in 
 the paper of  Cuntz and Quillen \cite{CQ2}. 

First, let us define the Hochschild cochain complex. 
\bdefinition
Let $\calm$ be a
bimodule over the algebra $A$. An $n$-cochain  $f$ on $A$ with 
values in $\calm$ is a linear map $f: \tp{A}{n}\ra \calm$. The space 
of such cochains is denoted by $C^n(A,\calm)$. The graded space 
$C^\star(A,\calm)$ is made into a differential complex using the following 
differential 
$$ 
\begin{array}{rcl}
(bf)(a_1, \dots , a_{n+1}) & = & a_1f(a_2, \dots , a_{n+1}) +\displaystyle{\sum_{i=1}^n}
(-1)^if(a_1, \dots, a_ia_{i+1}, \dots , a_{n+1})\\
\rule{0mm}{8mm}
& & \mbox{} + (-1)^{n+1} f(a_1, \dots , a_n)a_{n+1}
\end{array}
$$
\edefinition
Of particular interest are cochains on $A$ with values in the dual space $A^*$ which 
is an $A$-bimodule with the usual action of $A$: 
$$
\begin{array}{rcl}
(af)(a')& =& f(a'a)\\
(fa)(a')&= &  f(aa')
\end{array}
$$
It is clear that one can identify the space of linear maps $\tp{A}{n}\ra A^*$ 
with the space of linear functionals $\tp{A}{n+1}\ra \bc$. In this case, 
we obtain the Hochschild complex which calculates the Hochschild
cohomology $HH^\star(A)$, the dual theory to Hochschild homology 
defined in \ref{Hochschild-homology}.

There is a simplicial normalization theorem that states that 
Hochschild cohomology may be calculated using the {\em reduced} (or
simplicially normalized) cochain complex $(\bigoplus_{n} (\Omega^n A)^*, b)$, 
where $b$ is the transpose of the differential defined on  noncommutative
forms. 

\bdefinition
The cyclic cochain complex is the complex $(C^\star_\lambda(A), b)$ of 
cochains $f:\tp{A}{n+1} \ra \bc$ which are invariant under the action of 
cyclic group: $\lambda f = f$, 
where
$(\lambda f)(a_0, a_1, \dots , a_n) = (-1)^nf(a_n, a_0, \dots , a_{n-1})$. 
The differential $b$ is the restriction of the differential on the space of 
Hochschild cochains to the space of cyclic cochains:
$$
\begin{array}{rcl}
(bf)(a_0, \dots , a_{n+1}) & = & \displaystyle{\sum_{i=0}^n }(-1)^i f(a_0, \dots, a_i a_{i+1}, 
\dots , a_{n+1} ) \\
\rule{0mm}{8mm} && \mbox{} + (-1)^{n+1}f(a_{n+1}a_0, \dots , a_n)
\end{array}
$$
\edefinition
This complex computes the cyclic cohomology of the algebra $A$. 
\bexample
Let $M$ be a compact oriented manifold, and let $A$
be the algebra of smooth functions on $M$. Put
$$
\phi(\astr{0}{n}) = \int_M a_0 da_1\dots da_n
$$
for $a_0, \dots, a_n\in A$. Then   $\phi$ is a cyclic cocycle.
\eexample

This example may be generalized in the following way. Recall that
by
a $k$-dimensional current on the manifold $M$ we mean a continuous linear
functional   $C : \Omega^k(M)\ra \bc$. For example, if $k=0$, $C $
is a distribution on $M$. The value of $C $ on any $k$-form
$\omega$ is denoted by
$$
\omega \mapsto \int_C \omega
$$
A $k$-current  $C$ is called closed, if
$$
\int_C d\eta = 0
$$
for all $k-1$-forms $\eta$. One then has the following.
\blemma
If $C$ is a closed $k$-current, then
$$
\phi(\astr{0}{k}) = \int_{C} a_0 da_1\dots da_k
$$
is a cyclic cocycle.
\elemma

\bexample
Let $S$ be a Clifford bundle over a compact manifold
$M$ and let $D$ be the corresponding Dirac operator. Let $V$ be a vector
bundle over $M$, represented as the range of a projection-valued function
$e$. Then the Dirac operator $D_V$
on $S$ with coefficients in $V$ is given by $D_V = e(D\otimes 1)e$. Assume
that $\ker D=0$; then $D^{-1}$ is a bounded operator. Define
$$
\tau (f_0, f_1, \dots, f_{2p}) = {\rm Tr}_s
(D^{-1}[D,f_0] D^{-1}[D,f_1]\dots
D^{-1}[D,f_{2p}])
$$
where $f_i$ are (possibly matrix valued) smooth functions on $M$. Here
$Tr_s$ denotes a supertrace. (More on supertraces in the 
next section.) We note that for large enough $p$, the
expression in brackets is a trace class operator.
\blemma
$\tau$ is a cyclic cocycle.
\elemma

The importance of cyclic cocycles defined by the above formula is
illustrated by the following statement.
\blemma
If $e$ is a projection-valued function corresponding to a vector bundle
$V$, then
$$
{\rm Index} \, D_V= \tau(e, \dots, e)
$$
\elemma
\eexample
\bexample
Let $(A, \bh, F,\gamma)$ be an even $p$-summable Fredholm module. 
Then 
$$
\tau_{2n } = {\rm Tr}\;\gamma a_0[F,a_1]\dots [F, a_{2n}]
$$
is a cyclic cocycle called the {\em character} of the Fredholm module. 
\eexample
\subsection{Supertraces}\label{Supertraces}
We finish the chapter with  a discussion of 
$\bz_2$-graded cochain complexes introduced by Connes in \cite{Connes:Entire}
which will be needed in our discussion of entire cyclic cohomology of Banach
algebras. The material in this section is taken from \cite[\S 8]{CQ2}. 

A cochain $f$ is now an element of $(\Omega A)^*$, the dual space to $\Omega A$. 
Hence $f$ is the same as a sequence $(f_n)$, where $f_n$ is a linear 
functional on $\Omega^n A$. 

Let us identify the Cuntz algebra $QA$ with the algebra $\Omega A$
equipped with the even-odd grading and the Fedosov product \ref{Fedosov}
$$
\omega\circ \eta = \omega\eta - (-1)^{\deg (\omega)}d\omega d\eta
$$

\bdefinition
A cochain $f$ on $\Omega A$ is a {\em supertrace} on $QA$ if and 
only if for any two homogeneous elements $\omega$, $\eta $ of $QA$
$$
f(\omega\circ \eta) = (-1)^{\deg(\omega)\deg(\eta)} f(\eta\circ \omega)
$$
In other words, a supertrace is a cochain vanishing on supercommutators
with respect to the Fedosov product. 
\edefinition
We have the following result \cite[Propn. 8.1]{CQ2}. 
\bproposition
The following conditions are equivalent for $\tau\in (\Omega A)^*$: 
\begin{enumerate} 
\item $\tau $ is a supertrace on $QA$. 
\item Connes' infinite cycle identity holds {\em \cite{Connes:Entire}}
$$
\tau([\omega, \eta]) = (-1)^{\deg(\omega)}2\tau(d\omega d\eta)
$$
Here 
$[\omega, \eta] = \omega \eta - ( - 1)^{\deg(\omega)\deg(\eta)}\eta\omega$
is the supercommutator with respect to the ordinary product on $\Omega A$. 
\item $\tau $ is $\kappa$-invariant, i.e. $\tau\kappa = \tau$, and
$\tau_{n-1}b = 2\tau_{n+1}d$ for $n>0$. 
\end{enumerate}
\eproposition
Proof. Assume that $\tau$ satisfies Connes' cycle identity. Then 
$\tau([d\omega, d\eta]) = 0 $.  Thus if we apply $\tau $ to the 
Fedosov supercommutator $[\omega, \eta]_s$ 
$$
[\omega, \eta]_s = [\omega, \eta] -(-1)^{\deg(\omega)}2 d\omega d\eta
+ ( -1)^{\deg(\omega)}[d\omega , d\eta]
$$
we see that $\tau([\omega , \eta]_s) = 0 $, hence $\tau$ is a supertrace on 
$QA$. Conversely, if $\tau$ is a supertrace on $QA$, then $\tau$ vanishes 
on terms of the form $[d\omega, d\eta]_s = [d\omega , d\eta]$, which 
implies that $\tau $ satisfies Connes' identity. 

To see that a supertrace $\tau$ satisfies the identities in point $3)$ above, 
we note that as the Cuntz algebra $QA$ is generated by $a$ and $da$, 
$\tau$ is a supertrace if and only if it vanishes on the supercommutators
$$
[\omega, a]_s = (-1)^{\deg(\omega)}(b - ( 1 + \kappa)d)(\omega da)
$$
and 
$$
[\omega, da]_s = (1 - \kappa)(\omega da)
$$
\qed

Let us now define the rescaling operator $z$ on $\Omega A$ which multiplies 
by $z_n$ on $\Omega^n A$, where we define
$$
z_{2m} = \frac{(-1)^m}{m!}, \qquad \qquad z_{2m+1}= \frac{(-1)^m2^m}{(2m+1)!!}
$$
where $(2m+1)!! = 1 \cdot 3\cdots (2m+1)$. This rescaling is precisely what is 
needed to translate the last condition of the previous statement into 
a cocycle condition in the $b + B$-complex $((\Omega A)^*, b + B ) $. 
\bproposition\label{supertrace-condition}
$\tau$ is a supertrace on $QA$ if and only if the rescaled cochain $\tau^z = \tau \circ z$ 
is $\kappa$-invariant and such that $\tau^z (b + B) = 0$ except in degree zero. 
\eproposition

In degree zero a slight modification is needed, but it only affects odd supertraces
on $QA$ (supertraces whose nonzero components are in odd degrees). 
In even degrees we have a bijection between even, $\kappa$-invariant 
$b+B$-cocycles and even supertraces on $QA$. 
In odd degrees $\kappa$-invariant odd $b+B$-cocycles correspond 
to odd supertraces $\tau$ such that $\tau(da)= 0 $ for any $a\in A$. 

For future use we record the following fact \cite[Propn. 8.5]{CQ2}. 
\bproposition\label{fbB}
For a cochain $f \in (\Omega A)^*$, the following are equivalent: 
\begin{enumerate}
\item $f$ is $\kappa$-invariant, i.e. $f(1-\kappa) = 0 $. 
\item $fbB= 0 $ and $fP= f$. 
\end{enumerate}
\eproposition
It follows in particular that harmonic cocycles, i.e. cocycles $f$ such 
that $fP = f$ are automatically $\kappa$-invariant.

\vfill\eject

\section{Chern characters}
\subsection{$K$-theory and Chern classes}
One of the features of the topological  $K$-theory which makes
it so useful in a variety of applications is the 
existence of the Chern character homomorphism. It is also one of 
the key properties of cyclic cohomology, as we shall see presently. 
To set the necessary background for our discussion, we shall review 
briefly the construction of the Chern character in the topological 
$K$-theory. 

Let $E$ be a complex vector bundle over a compact topological space
$X$. We can associate with it  certain cohomology classes $c_i(E) \in H^{2i}(X, \bz)$, 
which are called the {\em characteristic classes of $E$}. These classes 
are used to classify bundles in a certain way and their 
basic properties are as follows: 
\begin{enumerate}
\item $c_0(E) = 1 \in H^0(X, \bz)$. 
\item For all $n\geq 0 $, $c_n(E\oplus F) = \sum_{i+j= n } 
c_i(E) \cup c_j(E)$. 
\item If $f: Y\ra X$ is a continuous map, then 
$c_n(f^*E) = f^*c_n(E)$. 
\end{enumerate}
The second of the above axioms, the Whitney sum formula,  implies 
that $c(E) = \sum_{i >0} c_i(E)\in H^\star(X, \bz)$  depends 
only on the class of the bundle $E$  in $K^0(X)$. Thus 
$c_i$ extend to functions $c_i: K^0(X) \ra H^{2i}(X, \bz)$ and
we can define the Chern character ${\rm ch} :  K^0(X) \ra H^*(X, \bq)$, 
which satisfies the following properties: 
\begin{enumerate}
\item ${\rm ch}_0(E)$ equals the {\em rank} ${\rm r}(E)$ of $E$, 
where ${\rm r}(E)\in H^0(X, \bz) \simeq \bz$. 
\item ${\rm ch} (E \oplus F) = {\rm ch}(E) + {\rm ch}(F)$. 
\item ${\rm ch}(f^*(E))= f^*({\rm ch}(E))$, for a continuous map 
$f: Y\ra X$. 
\end{enumerate}
We remark here that the first of these axioms makes an explicit connection 
with the rank function defined in (\ref{rank-function}. It turns out that 
the rank of a vector bundle gives rise to a characteristic class in 
degree zero. 
In the case when $X$ is a smooth manifold, there is a useful explicit 
description of the Chern character. We assume that  $E$ is 
a smooth vector bundle equipped
with connection $\nabla_E$, whose curvature is $\nabla_E^2$. 
Then the Chern character ${\rm ch}(E)\in H^*(X,\br)$
is represented by the closed differential form 
$$
{\rm ch}(E) = \tr(\exp(\nabla^2_E/2\pi i))
$$
To obtain numerical invariants of $X$,  we can
associate with any closed de Rham current $C$,   a map $I_C$
from $K^*$ to $\bc$ defined by 
$$
I_C(E) = \langle C, {\rm ch}(E)\rangle
$$

\mbox{}

This geometric construction can be reformulated in algebraic 
terms, when we replace the topological $K$-theory by 
the $K$-theory of a suitable ring or a $C^\star$-algebra. 
Connes's great discovery was that in such algebraic case
it is possible to construct a Chern character which takes its values
in cyclic cohomology (or more precisely, the periodic 
cyclic cohomology). We shall indicate the main points of 
this translation in the following two sections.

\subsection{Higher traces}\label{Higher-Traces}
We set out to discuss the algebraic version of the Chern-Weil 
theory of characteristic classes in a vector bundle. Let us first 
consider two unital rings 
$A$ and $B$ and a {\em based linear} map $\rho: A \ra B$ with curvature
$\omega$ (cf. (\ref{based-linear}):
$$
\omega(a_1,a_2) = \rho(a_1a_2) - \rho(a_1)\rho(a_2)
$$
If we assume that the linear map $\rho$
is to play the role of a linear connection in a vector bundle, then 
it is reasonable to treat the bilinear form $\omega$ as the 
formal analogue of the curvature of $\rho$. Since $\omega$ vanishes
if either of its arguments is the identity of $A$, it is 
in fact a two-cochain $\omega: d\Omega^1A\subset \Omega^2 A \ra \bc$. 
It turns out that one can reproduce practically the complete formal 
structure of the Chern-Weil theory starting with the \lq connection'
$\rho$ and its curvature $\omega$, and this was done in detail in 
Quillen's paper \cite{Q:Cochains}. For example, 
there is the Bianchi identity
$$
b'\omega = - [\rho, \omega]
$$
which is not difficult to check directly. There is a natural cyclic
cocycle that may be regarded as the algebraic analogue of the 
Chern character: $\tau(\omega^n/n!)N$, 
for any trace $\tau: B\ra \bc$. We shall use this 
cocycle to construct Connes classes in cyclic cohomology associated 
with so higher traces. 

Let $R$ be an algebra and let $I$ be an ideal in $R$.
Furthermore, let $u: A\ra R/I$ be an algebra isomorphism.
\bdefinition
Given an algebra $A$, a pair consisting of an algebra extension
$A = R/I$ and a trace $\tau : R/I^{m+1} \ra \bc$ on the quotient algebra
$R/I^{m+1}$, for $m\geq 0$, is called an {\em even higher trace} on the
algebra $A$.

An {\em odd higher trace}  on an algebra $A$ is an extension $A= R/I$
together  with a linear map $\tau: I^m\ra \bc$ satisfying $\tau [ R,I^m] = 0 $,
where $m>0$.
\edefinition
 
We recall here that if $B$ is an algebra then a trace $\tau : B \ra \bc$
is a linear map that vanishes on commutators 
of elements in $B$, i.e., it annihilates the space $[B,B]$. Similarly, if $M$ is 
a $B$-bimodule, then $\tau: M\ra \bc$ is a trace on $M$ iff it vanishes 
on $[B,M]$. This notion may be generalized slightly by considering 
traces with values in a vector space $V$. 

In his article \cite{Q:Extensions}, Quillen constructs certain cyclic 
cocycles associated with higher traces. 
In the case of an even higher trace $\tau$ we
choose a linear lifting $\rho: A\ra R$ and consider its curvature $\omega$
$$
\omega(a_1, a_2) = \rho(a_1a_2)- \rho(a_1) \rho(a_2)
$$

Since $\rho$ is an homomorphism of algebras modulo the ideal $I$, the
curvature $\omega$ has values in $I$. This data can be used to construct the
following Chern-Simons form,
\beqn\label{Chern-Simons}
{\rm cs}_{2n+1}= \int^1_0 \tau(\rho(tb'\rho +
t^2\rho^2)^n/n!)N dt\in C^{2n}_\lambda(A)
\eeqn
It is a cyclic $(2n)$-cocycle for $n\geq m$, hence the name even
higher trace for the trace $\tau$. This cocycle determines a class in
cyclic cohomology which we call a {\em  Connes class} and denote
${\rm cn}_{2n}(\tau)$. This is a cyclic cohomology class of degree
$2n$ associated with the higher trace $\tau$.

In the case of an odd higher trace, we can construct the Chern character
form
\beqn\label{Chern-character}
{\rm ch}_{2n}(\tau ,\rho) = \tau (\omega^n/n!)N \in C^{2n-1}_\lambda(A)
\eeqn
which is a cyclic $(2n-1)$-cocycle for $n\geq m$. Its cohomology class
will be denoted by ${\rm cn}_{2n-1}(\tau)$ and called the Connes
class associated to the odd higher trace $\tau$.
\bremark
In the case, when the trace $\tau $ takes values in a vector space
$V$, the cyclic cocycles obtained this way are also $V$-valued. 
Since we assume that $V$ is a trivial $A$-module, the notion 
of a cyclic cocycle extends in a straightforward way to include
$V$-valued cyclic cocycles.
\eremark
To summarise, we have defined two types of cyclic cohomology classes
associated with higher traces. 
In the even case, the Connes class is represented by an algebraic
analogue of the Chern-Simons form, whereas in the odd case, 
the representative of a Connes class is the analogue of the Chern 
character. 
Connes classes are related by the action of the $S$-operator, in a way which
resembles Bott periodicity in the standard $K$-theory.
\btheorem
The Connes classes associated to a higher trace are independent of the
choice of the lifting $\rho$ and they satisfy
$$
S\,{\rm cn}_j(\tau) = {\rm cn}_{j+2}(\tau)
$$
\etheorem
This theorem was first proved in the odd case by Connes in \cite{Connes:IHES}
and in both cases by Quillen in \cite{Q:Extensions}.

\subsection{Even higher traces and $K_0(A)$}
Higher traces and Connes classes are key elements in establishing a relation 
between cyclic cohomology and algebraic $K$-theory. We begin with a discussion of 
the even case, following closely the paper
of Quillen \cite[ II.2]{Q:Extensions}.

Recall that any idempotent matrix
$e\in M_r(A)$ represents a class in $K_0(A)$, when $A$ is a unital algebra.
Any  trace $\tau:A\ra V$ on $A$ induces a map of abelian groups
$$
\tau_*: K_0(A)\ra V
$$
which is defined by
$$
\tau_* :[e]\mapsto \tilde{\tau}(e),
$$
where $\tilde{\tau}$ is the extension of the trace $\tau$ to the matrix
algebra $M_r(A) = A\otimes M_r(\bc)$
$$
\tilde{\tau} = \tau \otimes \tr: A\otimes M_r(\bc) \ra V.
$$
We use here $\tau$ to denote the usual matrix trace on $M_r(\bc)$. 

Let us now assume that we are given an even higher trace on the algebra
$A$, i.e. an extension $A= R/I$ and a trace $\tau:R/I^{m+1} \ra V$.
We can define an index map associated with this higher trace as follows.
We put ${\rm ind}_\tau$ equal  to the following composition
$$
K_0(A) \arr{\simeq} K_0(R/I^{m+1}) \arr{\tau_*} V.
$$
The first map is induced by the surjection $R/I^{m+1} \ra R/I$.
Since the ideal $I$ is nilpotent in $R/I^{m+1}$, the kernel of this
map is nilpotent and the induced map on $K$-theory is an isomorphism.
The reason for that is  basically that any idempotent  matrix over $A$ can be
lifted to   an idempotent matrix over $R/I^{m+1}$. To describe this lifting process 
in more detail, let $e$ be an idempotent matrix in $M_r(A)$
 and let $f= \rho(e)$ be its lift in $M_r(R)$. The matrix $f$ is an idempotent
modulo the ideal $M_r(I)$ which means that the matrix $f - f^2$ has entries
in $I$.  We want to modify this matrix so that it becomes an idempotent matrix.

To do this, let us consider the polynomial ring $k[x]/(x(1-x))^{n+1}$. From
the Chinese Remainder Theorem  we can decompose
$$
k[x]/(x(1-x))^{n+1} = k[x]/x^{n+1} + k[x]/(1-x)^{n+1}
$$
Let $f_n(x)$ be the unique polynomial of degree not greater than $2n+1$ whose class in
$k[x]/(x(1-x))^{n+1} $ is an idempotent and $f_n(0) = 0 $ and $f_n(1) = 1$.
Thus $f_n(x) \cong x \, \mod\, x(1-x)$.

We note that for $n\geq m$ we have  a homomorphism
$$
k[x]/((x-x^2)^{n+1}) \ra M_r(R/I^{m+1}).
$$
which is defined by sending $x$ to $f$. By the Chinese remainder theorem,
there is a unique idempotent $p$ in $k[x]/((x-x^2)^{n+1})$ which reduces to
$x$ modulo $(x-x^2)$. The idempotent lifting of $e$ is then the image of $p$
under this homomorphism.

To represent $p$ (uniquely) by a polynomial of degree not greater that
$2n-1$ we introduce the following function

$$
\begin{array}{rcl}
f_n(x) & = & \frac{\displaystyle\int_0^x t^n(1-t)^ndt}{\displaystyle\int^1_0
t^n(1-t)^ndt}\\
\rule{0mm}{8mm}
& =  & \displaystyle{\frac{(2n + 1)!}{(n!)^2}\int^x_0}(t-t^2)^ndt \\
\rule{0mm}{8mm}
& = & \displaystyle{\frac{(2n+1)!}{(n!)^2}\int^1_0} x(tx-t^2x^2)^ndt
\end{array}
$$
Note that $f_n(1)=1$.

Let $\rho: A\ra R$ be a linear lifting of $A$ into $R$ and let $\tilde{\rho}$
denote the extension of $\rho$ to matrices. From the  above construction
follows that for $n\geq m$ the image of $f_n(\tilde{\rho}(e))$ in the algebra
$M_r(R/I^{m+1})$ is an idempotent matrix reducing to $e$. Thus we have
$$
\mbox{ind} _\tau[e] =
\frac{(2n+1)!}{(n!)^2}\int^1_0\tilde{\tau}\tilde{\rho}(e)(t\tilde{\rho} (e) -
t^2\tilde {\rho}(e)^2)^n dt
$$
Let us denote by $(e)_{2n+1}$ the image of $e^{\otimes 2n+1}$ in
$C^\lambda_{2n}(M_r(A))$. This is a  cyclic $2n$-cycle in the cyclic complex
of $M_r(A)$.
\bdefinition
The Chern character homomorphism is a map $\mbox{\rm Ch}_{(2n)} : K_0(A) \ra
HC_{2n} (A) $ defined by
$$
\mbox{\rm Ch} _{(2n)} [e] = \frac{(2n)!}{n!} [\tr(e)_{2n+1}]
$$
\edefinition
This homomorphism was defined by Connes in \cite{Connes:IHES}.

\btheorem
The index map associated to the higher trace $\tau$ is given by pairing the
associated Connes cyclic cohomology class with the Chern character:
$$
\mbox{ind}_\tau [e] = \langle \mbox{\rm cn}_{(2n)}(\tau),
\mbox{\rm Ch}_{(2n)}[e]\rangle
$$
\etheorem
Proof. We use the definition of Connes classes to calculate the following
$$
\begin{array}{rcl}
\langle \mbox{cs}_{2n+1}(\tau, \rho) \frac{(2n)!}{n!} [\tr(e)_{2n+1}]\rangle
& = & \langle \mbox{cs}_{2n+1}(\tilde{\tau}, \tilde{\rho})\displaystyle{
\frac{(2n)!}{n!}} [\tr(e)_{2n+1}]\rangle \\
\rule{0mm}{8mm}
& = &
\displaystyle{
\frac{(2n)!}{(n!)^2}}\displaystyle{\int^1_0}\tilde{\tau}(\tilde{\rho}(e)
(tb'\tilde{\rho}(e,e) + t^2\tilde {\rho}^2(e,e))^n N (e)_{2n+1})dt
\end{array}
$$
Now use that $N(e)_{2n+1} = (2n+1) (e)_{2n+1}$ (this follows from
the fact that $(e)_{2n+1}$ is a cyclic cycle), and the cochain formulae
$$
\begin{array}{rcl}
b'\tilde{\rho}(e,e)& = &\tilde{\rho}(e\cdot e) = \tilde{\rho}(e)\\
\tilde{\rho}^2(e,e) & = & -\tilde{\rho}(e)\tilde{\rho}(e)
\end{array}
$$
to check that this expression is the same as the formula that we have derived
for the index. \qed

\subsection{Toeplitz operators}
As an example of a concrete situation corresponding to the 
theory explained above we shall discuss Toeplitz operators 
and a certain index theorem. 

Let $\bh = L^2 (S^1, d\theta/2\pi)$ be the space of square-integrable
functions on the circle $S^1$. This space is a module over the
algebra of smooth functions $A= C^\infty (S^1)$, where the module
structure is given by pointwise multiplication. The family of functions
$\{z^n\}$, for $z\in S^1$ and $n\in \bz$ is a basis for $\bh$. By definition,
the {\em Hardy space} is the closure of  the  linear span of
$\{ z^n \}$ for $n\geq 0$. Let us denote by $e$ the 
projection from $\bh$ onto the Hardy space. 

We define a Toeplitz operator for any function $f\in A$ by
$$
T_f(g) = e(fg)
$$
for any $g$ in the Hardy space. Hence we can write that $T_f = efe : \bh \ra \bh$.
The Toeplitz operator  has the following matrix form with respect to
the decomposition  $\bh= e\bh \oplus (1-e) \bh $ given by the projection $e$.
$$
T_f= \left(
\begin{array}{cc}
efe & e f (1-e) \\
(1-e) f e & (1-e) f (1-e) \end{array}\right)
$$
\blemma
If $g\in A$ then $[e,g]$ is an operator with a smooth Schwartz kernel,
in particular it is of trace class. Also, for $f,g \in A$
$$
\tr _\bh (f[e,g]) =  1/2\pi i \int _{S^1} fdg
$$
\elemma
Proof.  Use the Fourier transform of the function $g$ to write
$$
g(\theta) = \sum_n e^{in\theta} \int e^{-i n \theta'} g (\theta') d\theta' /
2 \pi.
$$
Since
\begin{eqnarray*}
(eg) (\theta) &=& \sum _{n\geq 0 } e^{in\theta} \int e ^{-in\theta'}
g(\theta') d\theta '/2\pi \\
& = & \int \sum _{n\geq 0} e^{i n (\theta - \theta')} g (\theta') d \theta'/2
\pi\\
& = & \int \frac {1}{1 - e^{i(\theta - \theta')}}g ( \theta') d\theta '/2\pi
\end{eqnarray*}
the kernel function of the operator $e$ is
$$
e(\theta, \theta') = \frac{1}{1 - e^{i(\theta - \theta')}}.
$$
A very similar calculation yields the kernel function of the operator $[e, g]$:
$$
[e,g] ( \theta - \theta') = \frac{g(\theta - \theta')}{ e^{i(\theta -
\theta')}}
$$
which is a smooth function over the torus. The trace of the operator
$[e,g]$ is then given by
$$
\tr _\bh (f[e,g]) = \int f(\theta) ([e,g](\theta, \theta)) d\theta /2 \pi=
1/2\pi \int f(\theta) g'(\theta) d \theta
$$
\qed

We may state the result of the last lemma in a slightly different way. 
Let $R=L(\bh)$ be the algebra of bounded operators
on $\bh$ and let $I = L^1 ( \bh)$ be the ideal of trace class operators on
$\bh$. If $e$ defines the projection onto the Hardy space as before, then
let us define a linear map $\rho : A \ra R $ by
$$
\rho (f) = efe = T_f
$$
It follows from the last lemma that $\rho $ is a homomorphism of algebras modulo
trace class operators.
\bremark
Note that
$$
\phi (f, g ) = \tr_\bh (\rho(f)\rho (g) - \rho (fg)) - \tr _\bh (\rho(g)\rho
(f) - \rho(gf))
$$
is a cyclic cocycle.
\eremark
\blemma
In the notations introduced above
$$
\phi(f,g) = \tr_{\bf H} (f[e,g])
$$
\elemma
Proof. Take a Toeplitz operator $\bt\in I$. Then
$$
\tr_{\bf H} (e\bt)= \tr_{\bf H} (e\bt e ) = \tr_{e\bf H}(e\bt e )
$$
which follows if we use the matrix representation of $\bt $. 

Now we note that if $D$ is a derivation, and $e$ is a projection, then $(De)e\,
+ eDe = De$
and $eDe = De\, (1-e) $ , $ (1-e) De = (De) e $.
Thus $ e[e,g] = [e,g](1-e)$, etc. Using these identities we calculate
\begin{eqnarray*}
\tr _\bh (\rho(f)\rho(g) - \rho(fg)) & = & \tr_{e\bh} ( ef[e,g]e) \\
& = & ( ef[e,g])
\end{eqnarray*}
and
\begin{eqnarray*}
\tr_{e\bh}( \rho(g) \rho(f)) - \rho(gf) & = & \tr_{e\bh} ( e f ege - egfe) \\
&=& \tr_{e\bh}(eg(e-1) fe ) \\
& = & \tr_{\bh} (eg(e-1) f ) \\
 & = & \tr_{\bh}( [e,g](e-1) f )\\
& = & \tr _{\bh}((e-1)f[e,g]).
\end{eqnarray*}
If we subtract the two expressions, we get the required formula. \qed

If $f$ is an invertible function on the circle, then the Toeplitz operator
$\rho(f) : e\bh \ra \bh$ has a parametrix $\rho(f^{-1})$, i.e.,
$\rho (f)$ is a Fredholm operator. The index of such an operator is defined
by
$$ {\rm index} \; \rho (f) = \tr_{e\bh} ( 1 - \rho(f^{-1}) \rho(f) ) -
\tr_{e\bh}(1-\rho(f)\rho(f^{-1})).
$$

It follows from our calculations that
$$
{\rm index} \; \rho(f) = \phi(f^{-1}, f ) = 1/2 \pi i\int_{S^1} f^{-1} df
$$
Thus the index of a Toeplitz operator associated to a function
$f$ is the  winding number of $f$.

\subsection{Odd higher traces and $K_1(A)$}
We want to show how to construct an index map associated to an odd higher trace. 
For this we shall study again the connecting homomorphism in the 
$K$-theory exact sequence
$$
K_1(I)\ra K_1(R)\ra K_1(R/I) \arr{\delta} K_0(I) \ra K_0(R) \ra K_0(R/I).
$$
which we encountered in \ref{excision-prototype}. 
Let us recall the main idea of Milnor's patching construction (see \ref{patching}). 
  For any $u\in GL_r(R/I)$
we construct  a finitely generated projective module $P_u$ over $\tilde{I}$.
One proves that $P_u$ is projective by lifting $u\oplus u^{-1}$ to an
invertible matrix $T\in M_{2r}(R)$. This gives $P_u$ as the image
of the projection operator
$$
e = T e_0 T^{-1}, \qquad e_0 = \left(\begin{array}{cc} 1 & 0 \\ 0 & 0
\end{array}\right)
$$
acting on the space $(\tilde{I})^{2r}$.
Then the image $\delta(u)$ is represented by the idempotent $e$. We shall now
describe this procedure in some detail \cite{Q:Extensions}.

Let $p$ and $q$ be the lifts to $M_r(R)$ of $u$ and $u^{-1}$ respectively.
We define $x$ and $y$ to be such that $qp = 1 - x$ and $ pq = 1- y$. It follows 
from these definitions that the matrices $x$ and $y$ have entries in
$I$ and that $xq = qy$. Using this remark and applying induction we show that 
$x^i q = q y^i$. Let $n\geq 1$. Our goal is to construct an invertible
matrix $T$, 
which is the lift of $u\oplus u^{-1}$ to $M_{2r}(R)$.
For this we need to modify the parametrix $q$ in the following way. Let
$$
q_n = (1+ x + \dots + x^{2n-1} ) q = q ( 1 + y \dots + y^{2n-1}).
$$
Then it is not difficult to check that $q_np = 1- x^{2n}$ and $pq_n = 1 -
y^{2n}$. We  have
$$
\left( q_n\;\; x^n\right) \left(\begin{array}{c} p \\ x^n\end{array}\right) = 1
$$
and from this follows that the matrix
$$
e_n = \left(\begin{array}{c} p\\x^n \end{array}\right)\left( q_n\;\; x^n\right) =
\left(\begin{array}{cc} 1 - y^{2n} & px^n \\
x^nq_n & x^{2n}
\end{array}\right)
$$
is an idempotent with entries in ${\tilde{I}}$. It is now clear that 
the matrices
$$
T= \left(\begin{array}{cc} p & - y^n \\ x^n & q_n \end{array}\right),
\qquad T^{-1} = \left( \begin{array}{cc} q_n & x^n \\ - y^n & p \end{array}
\right)
$$
are inverses of each other over $R$. Moreover,  $T$ reduces to
$u \oplus u^{-1} $ modulo the ideal $I$ and  $T e_0 T^{-1} = e_n$.
Thus the idempotent  $e_n$ represents the class $\delta[u]$ for any $n\geq 1$.

\mbox{}

We proceed with the construction of the index map on $K_1(R/I)$. 
Let $\tau$ be a trace defined on $I$, regarded as an algebra, which
means that $\tau[I,I] = 0 $. There is an induced map $\tau_* : K_0(I)
\ra V$ which  can be described as follows. We first
extend $\tau$ to a trace on ${\tilde{I}}$  by putting it to
be zero on $\bc$. Then we extend it to matrices over  ${\tilde{I}}$ in the usual
way, using the formula $\tilde{\tau} = \tau\otimes \tr$, where $\tr$ is the standard
matrix trace. Next we define $\tau_*$ on
$K_0(I)$ by applying $\tilde{\tau}$ to idempotent matrices over $\tilde{I}$.
We find that  
\beqn\label{we-find-that}
\begin{array}{rcl}
\tau_*(\delta[u]) & = & \tilde{\tau}(e_n)=  (\tilde{\tau}(x^{2n} - y^{2n}))\\
& = & \tilde{\tau}\left((1 - qp)^{2n} -  ( 1 - pq)^{2n}\right)
\end{array}
\eeqn

\blemma
Assume that $R/I$ is nilpotent, i.e. $R^m\subset I$ for some $m$. Then
$K_0(I)$ maps onto $K_0(R)$ and the kernel is killed by $\tau_*$, so that
$\tau_* $ descends to a map defined on $K_0(R)$.
\elemma
Proof. The assumption of nilpotence implies that  $K_0(R/I)=0$. From the $K$-theory exact
sequence then follows that the map $K_0(I) \ra K_0(R)$ is surjective. The
kernel of this map is the image of $\delta$, so we need to show that
$\tau_*\delta=0$. By nilpotence, any matrix $u = 1- w$ with $w \in M_r(R/I)$
is invertible with inverse
$$
u^{-1} = 1 + w + \cdots + w^n
$$
for some $n$. Lifting $w$ to a matrix $z$ over $R$, we can take $p = 1-z $
and
$$
q = 1 + z + \cdots + z^n.
$$
Then $p$ and $q$ commute, and
$\tau_*\delta[u]=0$ from (\ref{we-find-that}). \qed

We are now ready to define the index map associated to an odd higher trace on. 
Let us then assume that the higher trace is given by
an extension $A= R/I$ and a linear map $\tau : I^m \ra V$ vanishing on
the commutator space $[R,I^m]$. If we apply the above lemma in the case
of the algebra $I$ and the ideal $I^{m}$, we see that the trace map
$\tau_*: K_0(I^m)\ra V$ descends to a map $\tau_* : K_0(I) \ra V$.
\bdefinition
The index homomorphism $\ind_\tau$ associated to the higher trace $\tau$
is the composition
$$
K_1(A)= K_1(R/I) \arr{\delta} K_0(I) \arr{\tau_*} V.
$$
\edefinition
This definition can be illustrated with the following diagram.
$$
\begin{array}{cccl}
\arr{} & K_1(R/I) &\arr{\delta} & K_0(I) \\
&\Buar{}& & \Buar{} \\
& K_1(R/I^m) & \arr{\delta} & K_0(I^m)  \arr{\tau_*} V
\end{array}
$$
Here the first horizontal row is a part of the $K$-theory sequence
associated to the algebra extension $A= R/I$, and the vertical arrows
are surjective maps, and the one on the right is  explained in the previous lemma.
The lemma  states that one can calculate  the index by taking
a class $[u]$ to $\delta [u]$ using the connecting homomorphism $\delta$,
 lifting it to $K_0(I^m)$ and then applying  the trace $\tau_*$.
We use this fact in the following statement.
\bproposition \label{index-formula}
Let $u$ be an invertible matrix over $A$ and let $p$ and $q$ be
matrices over $R$ reducing modulo $I$ to $u$ and $u^{-1}$ respectively.
Then
$$
\ind_\tau[u] = \tilde{\tau} \{ (1 - qp)^n - (1-pq)^n\}
$$
for $n\geq m$.
\eproposition
Proof.  
In our computation of the connection morphism 
$\delta$  we have seen that the class $\delta [u]$ 
is represented  by the idempotent $e_n$ over $(\tilde{I})^n$, where $n$ could be
chosen as large as we like. Hence taking $n\geq m$ we find that
$$
\ind_\tau (e) = \tilde{\tau} (x^{2n} - y^{2n})
$$
with $x= 1- qp $ and $y= 1- pq$. As $\tau[R,I^m] = 0$,
we have for $i\geq m$ that
$$
\begin{array}{rcl}
\tilde{\tau} (x^{i+1} - y^{i+1} ) & = & \tilde{\tau} (x^i - qpx^i - y^i +
pqy^i) \\
& = & \tilde{\tau}(x^i - y^i - qpx^i + px^iq) \\
& = & \tilde{\tau}(x^i - y^i)
\end{array}
$$
\qed

Let $u\in GL_r(\tilde{A})$ be such that $u-1 \in M_r(A)$, so that also
$u^{-1} - 1 \in M_r(A)$. We write
$$
(u^{-1}, u)_{2n} = (u^{-1} - 1, u-1, \dots , u^{-1} -1, u - 1) \in
M_r(A)_\lambda^{\otimes 2n}= C^\lambda_{2n-1}(M_r(A))
$$
This is a cycle in the cyclic complex of the matrix algebra $M_r(A)$. To
check this we identify the cyclic complex of $M_r(A)$ with the reduced
cyclic complex  of the augmented algebra $\widetilde{(M_r(A))}$.
Then the matrices $u$ and $u^{-1}$ can be viewed as elements of
$\widetilde{(M_r(A))}$,  and that $(u^{-1}, u , \dots, u^{-1}, u )$ is obviously
a cycle in this reduced cyclic complex. Thus $(u^{-1} - 1 , u- 1)_{2n}$,
which is the corresponding element of the cyclic complex of the algebra
$M_r(A)$ is a cycle.

\bdefinition
The Chern character homomorphism
$$
\mbox{\rm Ch}_{(2n-1)} : K_1(A)\lrar HC_{2n-1}(A)
$$
is given  by
$$
\mbox{\rm Ch}_{(2n-1)}[u] = (n-1)![\tr (u^{-1}-1, u-1)_{2n}]
$$
\edefinition
The numerical factor in the above definition is needed to make sure that
the Chern homomorphism is compatible with the $S$ operation.

We finish this section proving another index theorem-type statement. Recall
that  an odd higher traces give rise to Connes classes 
$\mbox{\rm cn}(\tau) \in HC^{2n-1}(A,V)$, for $n\geq m$,
which are represented by cyclic cocycles represented by the 
Chern character \ref{Chern-character}. 
\bproposition
The index map  associated to  the higher trace $\tau$ is given by pairing
the associated Connes cyclic cohomology class with the Chern
character
$$
\ind_\tau[u] = \langle \mbox{\rm cn}_{(2n-1)}(\tau) , \mbox{\rm
Ch}_{(2n-1)}[u]\rangle
$$
\eproposition
Proof. Let $\rho: A \ra R$ be a linear lifting of $A$ into $R$. We extend it
to a map $\tilde{A} \ra \tilde{R}$ so as to be the identity on $k$, and we denote by
$\tilde{\rho}$ the obvious extension to matrices. The right side
is then
$$
\begin{array}{rcl}
\langle \mbox{\rm ch}(\rho, \tau), (n-1)! \tr (u^{-1} - 1, u - 1)
_{2n}\rangle& = & \langle \mbox{\rm ch}_{2n}(\tilde{\rho}, \tilde{\tau}), (n-1)!
(u^{-1} - 1, u - 1) _{2n}\rangle  	\\
& = & n^{-1} \tilde{\tau}(\tilde{\omega}^n )N ( u^{-1} - 1 , u - 1 )_{2n}
\end{array}
$$
where we have used that
$$
\mbox{\rm ch}_{2n}(\tilde{\rho}, \tilde{\tau}) = \mbox{\rm ch}_{2n}(\rho,
\tau)\tr
$$
Now we use the following
$$
\begin{array}{rcl}
N(u^{-1} - 1, u - 1)_{2n} & = & n \{ ( u^{-1} -1, u-1)_{2n} - (u-1, u^{-1} -
1 )_{2n}\}\\
\tilde{\omega}(u^{-1} - 1 , u-1) & = & ( b'\tilde{\rho} + \tilde{\rho}^2)
(u^{-1} - 1 , u - 1 ) = 1 - \tilde{\rho}(u^{-1}) \tilde{\rho}(u)
\end{array}
$$
and similarly with $u^{-1}$ and $u$ interchanged. We now see that the right
hand side of the index formula is
$$
\tilde{\tau}\{(1-\tilde{\rho}(u^{-1})\tilde{\rho}(u))^n - (1 -
\tilde{\rho}(u)\tilde{\rho}(u^{-1}))^n\}
$$
which is the same as the formula for the index in Proposition 
\ref{index-formula}.
\qed

\vfill\eject

\section{Excision in cyclic homology and algebraic $K$-theory}
This chapter brings three important results in cyclic homology and 
algebraic $K$-theory. The first is the Loday-Quillen theorem on the relation 
between cyclic homology and the Lie algebra homology. This theorem 
establishes cyclic homology as as sort of additive $K$-theory, 
as will be explained below. Then we briefly describe the theorem of Wodzicki 
about excision in cyclic homology and the theorem of Suslin and 
Wodzicki on excision in the rational algebraic $K$-theory. The proof
of the last theorem relies on the result of Loday and Quillen and 
completes the discussion of excision begun in \ref{simple-excision}.
\subsection{The Loday-Quillen theorem}
If ${\frak g}$ is any Lie algebra over a field $k$, then its homology with coefficients
in $k$  is defined to be $H_n({\frak g}) = \Tor{U({\frak g})}{n}{k}{k}$,
where $U({\frak g})$ is the  universal enveloping algebra of ${\frak g}$.
There is a Koszul complex which  computes this homology. 
The space of $n$-chains $C_n(\frak g)$ in this complex is the 
$n$-th exterior power of $\frak g$, $\Lambda^n {\frak g}$. This complex 
is equipped with a degree $-1$ differential, which on the space 
of $n$-chains is defined by the following formula
$$
d(x_1\wedge \dots \wedge x_n) = \sum_{i=1}^n (-1)^{i+j}[x_i, x_j] \wedge
\dots \wedge \hat{x}_i \wedge \dots \wedge \hat{x}_j \wedge \dots \wedge
x_n
$$
We shall be mostly interested in Lie algebras associated with various 
matrix algebras. Let $A$ be an algebra over a field $k$. 
The algebra of matrices $M_r(A)$ becomes a Lie algebra if we equip it
with the obvious bracket $[x,y] = xy - yx$. This Lie algebra
will be denoted ${\frak gl}_r(A)$. Then for any $r$ there is an inclusion
${\frak gl}_r(A) \ra {\frak gl}_{r+1}(A)$, which is induced by 
the canonical inclusion defined for matrix algebras in previous sections. 
We put ${\frak gl}(A) $ to be the union of  all the Lie algebras
${\frak gl}_r(A)$.

When $A$  is a  unital algebra, then the Lie algebra
${\frak gl}(A)$ of matrices over $A$ contains the Lie algebra of scalar
matrices ${\frak gl}(k)$. Let us consider matrices of a
fixed  dimension $r$. The Lie algebra ${\frak gl}_r(k)$ acts on the chain
complex of  ${\frak gl}_r(A)$ by means of the Lie derivative
$$
X(v_1\wcw v_k) = \sum_{i=1}^k v_1\wcw [X,v_i] \wcw v_k
$$
for $X$ in ${\frak gl}_r(k)$ and $v_i$ in ${\frak gl}_r(A)$.
 The quotient of the Koszul complex $C_\bullet({\frak gl}(A))$ by the
above action is called   the complex of {\em coinvariants} and denoted
by $C_\bullet({\frak gl}(A))_{{\frak gl}(k)}$.

The complex $C_\bullet({\frak gl}(A))$ of the Lie algebra chains has a natural
structure of a coalgebra. The diagonal map
$$
\triangle :{\frak gl}_n(A) \lrar {\frak gl}_n(A) \oplus {\frak gl}_n(A)
$$
given by $X\mapsto X\oplus X$ induces the coproduct
$$
C_\bullet ({\frak gl}(A)) \lrar C_\bullet ({\frak gl}(A)) \otimes
C_\bullet({\frak gl}(A)).
$$

This coproduct induces a coproduct on homology thus making
$H_\bullet({\frak gl}(A))$ into a coalgebra. On the other hand there is the
following map
$$
{\frak gl}_n(A) \oplus {\frak gl}_n(A) \ra {\frak gl}_{2n}(A) \ra {\frak gl}(A)
$$
where the first arrow is given by the direct sum of matrices. This map
induced the map of  complexes
$$
C_\bullet ({\frak gl}_n(A) \oplus {\frak gl}_n(A) ) = C_\bullet({\frak
gl}_n(A)) \otimes  C_\bullet({\frak gl}_n(A)) \lrar C_\bullet ({\frak
gl}_{2n}(A))
$$
Taking the direct limits of both sides we obtain a product
$$
C_\bullet({\frak gl}(A)) \otimes C_\bullet({\frak gl}(A)) \lrar
C_\bullet({\frak gl}(A)).
$$
This product is neither  associative nor commutative, because in general
the matrix $(x\oplus y ) \oplus z$ is only conjugate to the
matrix $x\oplus (y\oplus z)$. In fact, these two matrices are related by a
permutation.  The same is true for the matrices $x\oplus y$
and $y\oplus  x$. However, on passing to the complex of coinvariants
$C_\bullet({\frak gl}(A))_{{\frak gl}(k)}$, one obtains a commutative,
associative product. The proof of this statement uses the fact that 
the complex of coinvariants $C_\bullet({\frak gl}(A))_{{\frak gl}(k)}$
is the same as the complex $C_\bullet({\frak gl}(A))_{GL(k)}$ together 
with the fact that matrices representing permutations are elements 
of $GL(k)$. 
\bproposition
The complex of coinvariants $C_\bullet({\frak gl}(A))_{{\frak gl}(k)}$
is a commutative and cocommutative graded Hopf algebra. 
\eproposition
This is not the only Hopf algebra which results from the two operations
defined here. It turns out that the Lie algebra homology $H_\bullet({\frak gl}(A))$
inherits the structure of a Hopf algebra. 

Recall that an element $x$ of a general Hopf algebra $H$ is called primitive
if and only if
$$
\triangle x = x\otimes 1 + 1 \otimes x
$$
where $\triangle$ is the coproduct in $H$. In the commutative and cocommutative 
case, $H$ is generated as a Hopf algebra by its {\em primitive part}
${\rm Prim}H$, which by definition consists of all its primitive elements. 
It is therefore a natural problem to  determine the primitive part of the Hopf algebra
$H_\bullet({\frak gl}(A))$. The answer turns out to be very interesting---the
primitive part of  $H_\bullet({\frak gl}(A))$ is  the cyclic complex of the
algebra $A$. This is the statement of the theorem of Loday and Quillen. 
A first step in the proof of this result is the following Proposition.  
\bproposition
The Lie algebra homology $H_\bullet({\frak gl}(A))$ is a commutative and
cocommutative graded  Hopf algebra, whose primitive part is the homology
of the complex ${\rm Prim}\,(C_\bullet({\frak gl}(A))_{{\frak gl}(k)})$.
\eproposition
This proposition simply states that the operation of taking the primitive
part commutes with the homology functor. It also indicates that we should
take a closer look at the complex of coinvariants 
$C_\bullet({\frak gl}(A))_{{\frak gl}(k)}$. For this we shall need some 
elementary facts from the invariant theory. A classic reference for this topic is the
book by Herman Weyl \cite{Weyl}, but it is known for  its unreadability.
There is a very nice proof of the basic theorems of the theory of invariants
given in the paper by Atiyah, Bott and Patodi \cite{ABP}, and a thorough 
discussion may  be found in  
Loday's monograph \cite{L}. We use
Loday's exposition in the  summary presented here. A brief description of the 
invariant theory
is also given in the lecture notes by Husemoller and Kassel \cite{HK} 
and in the paper by Feigin and Tsygan \cite{ft-LNM}.

We begin with the invariant theory for ${\frak gl}_r(k)$. A key observation
here is that if $M$ is a ${\frak gl}_r(k)$-module  with the usual 
Lie algebra action, and at the same time it is a $GL_r(k)$-module with the 
adjoint action, then there is an isomorphism $M_{{\frak gl}_r(k)}\simeq
M_{GL_r(k)}$. 

Let  $V$ be an $r$-dimensional vector space over $k$. The group
of automorphisms of $V$ is then the group $GL_r(k)$.
The group of permutations $\sigma_n$ of $n$ letters acts on the $n$-fold
tensor power $\tp{V}{n}$ of $V$ on the left as follows
$$
\sigma(\vstr{1}{n}) = (\vstr{\sigma^{-1}(1)}{\sigma^{-1}(n)}).
$$
We can extend this action by linearity   to the action of the group algebra
$k[\Sigma_n]$ of the symmetric group $\Sigma_n$. This way we obtain
a ring homomorphism $\mu: k[\Sigma_n] \ra \End(\tp{V}{n})$.
On the other hand, the group $GL_r(V)$ acts diagonally on the space
$\tp{V}{n}$ and by conjugation on the algebra $\End(V)$. The image of the
map $\mu$ falls into the part of $\End(V)$ which is invariant under the
action of $GL_r(V)$. A first fundamental theorem of invariant theory
states that the map
$$
\mu: k[\Sigma_n] \lrar \End(\tp{V}{n})^{GL_r(V)}
$$
is surjective. The space on the right hand side is the invariant subspace
of $\End(\tp{V}{n})$, that is the space of all  $\phi\in \End(\tp{V}{n})$ such that
$g^{-1} \phi g = \phi$ for any $g\in GL_r(V)$. This claim is then improved
to the following statement. (We assume that $k$ is a field of characteristic 
zero.)
\bproposition
The map $\mu $ is an isomorphism for $r\geq n$. Otherwise, its
kernel is a two sided ideal $I_{n,r}$ in the ring $k[\Sigma_n]$ generated by
$\epsilon_{r+1} = \Sigma_{\sigma\in \Sigma_{r+1}}\sgn(\sigma)\sigma$.
\eproposition

Let $V^*$ be the dual space of $V$. Then
$\End(V) = V\otimes V^*$. Moreover, $\End(V)\simeq \End(V)^*$ 
and $\End(\tp{V}{n}) = \tp{\End(V)}{n}$.
It is also not difficult to see that, if $G$ is a group, then 
for any $G$-module $M$ we have
$$
(M_G)^* = (M^*)^G.
$$
We can now  define a map
$$
T: k[\Sigma_n]\lrar \End(\tp{V}{n}_{GL_r(V)})^*
$$
as the composition of the map $\mu$ with the isomorphism
$$\End(\tp{V}{n})^{GL_r(V)} \simeq (\End(\tp{V}{n})_{GL_r(V)})^*$$.
\bproposition\label{mapT}
When $n \geq r = \dim \,V$, the  map $T$ is an isomorphism.
\eproposition

Let us denote the Lie algebra
associated with the algebra $\End(V)$ by ${\frak g}={\frak gl}_r(k)$.
For any permutation $\sigma\in \Sigma_n$ we define a map 
$$
\tr_\sigma: \tp{\frak g}{n} \arr{} k
$$
as follows. Let $(X^1, \dots , X^n) \in \tp{\frak g}{n}$. Suppose 
that the permutation $\sigma $ has the following decomposition 
as a product of cycles: $\sigma = \sigma_1\cdots \sigma_k$. For
a  cycle $\sigma_l = (i_{n_1}, \dots, i_{n_l})$, $l=1, \dots, k$ put
$X^{(\sigma_l)}$ to be the product $X^{(\sigma_l)}= 
X^{i_{n_1}}\cdots X^{i_{n_l}}$. The functional $\tr_\sigma $ is defined 
by
$$
\tr_\sigma(X^1, \dots, X^n) = \prod_{i=1}^k \tr X^{(\sigma_i)}
$$
\blemma
For any $\sigma\in \Sigma_n$, the functional $\tr_\sigma$ is ${\frak g}$-invariant. 
\elemma
Proof. If $\sigma$ is a cyclic permutation,
then for any $X\in {\frak g}$ we have
$$
\begin{array}{rcl}
-(X\tr_\sigma) & = & \displaystyle{\sum_{i=1}^n} \tr_\sigma (X^1, \dots,
[X,X^i], \dots,  X^n) \\
\rule{0mm}{5mm}
& = & \displaystyle{\sum _i }\tr(X^{\sigma(1)}, \dots, X^{\sigma(n)})\\
\rule{0mm}{5mm}
& = & \tr [X, X^{\sigma(1)}\cdots X^{\sigma(n)}] = 0
\end{array}
$$
which shows that $\tr_\sigma$ is an invariant functional in this case. 
In general, if
the permutation $\sigma$ is a product of two or more cycles, we note that
the action of ${\frak g}$ affects each cycle separately. Thus using the above
calculation we  show that $\tr_\sigma$ is ${\frak g}$-invariant for any
$\sigma\in \Sigma_n$. \qed

It turns out that the maps $\tr_\sigma$ are good models for invariant 
$n$-cochains on ${\rm gl}_r(k)$ in view of the following result. 
\bproposition\label{basis}
The family $\{ \tr_\sigma\}$, with $\sigma\in \Sigma_n$, spans the space of
invariant $n$-functionals on the Lie algebra ${\frak gl}_r(k)$. The maps
$\tr_\sigma $ are linearly independent if and only if $r\leq n$.
\eproposition
We remark that the isomorphism of proposition \ref{mapT}
is  given explicitly by
$$
T: \sigma \mapsto \tr _\sigma
$$

We now set out to find the primitive part of the Lie algebra homology
of ${\frak gl}(A)$. We do this in a few steps following the original paper of
Loday and Quillen \cite{LQ}. First we shall use the fact that the Lie algebra
${\frak gl}_r(k)$ is reductive. This will allow us to replace the standard
complex calculating the Lie algebra homology by a quasi-isomorphic complex
of coinvariants. Then we apply the invariant theory to calculate the complex
of coinvariant elements explicitly. In the end we identify the complex of
coinvariants with the cyclic complex.

First we recall some facts from the theory of Lie algebras.
\bdefinition
A representation of a Lie algebra is called  {\em simple} or {\em reducible}
if it  does not contain a nontrivial subrepresentations. A representation
of a Lie algebra is called {\em semisimple} if  it is a direct sum
of simple representations. Finally, a Lie algebra is called reductive if
its adjoint representation is semisimple. 
\edefinition
\bexample\label{reductive}
The Lie algebra ${\frak gl}_r(k)$ is reductive for any finite value of $r$.
\eexample
\bproposition\label{reductive2}
Let ${\frak g}$ be a Lie algebra, and let $\frak h$ be a Lie subalgebra. Let
$V$ be  a ${\frak g}$-module, and assume that ${\frak g}$ and $V$ are semisimple
$\frak h$-modules.  Then the surjective map $C_\bullet({\frak g}, V) \ra
C_\bullet({\frak g}, V)_{\frak h}$  is a quasi-isomorphism, i.e. it induces
an isomorphism
$$
H_\bullet({\frak g}, V) \cong H_\bullet(C_\bullet({\frak g}, V)_{\frak h}, d)
$$
\eproposition
Proof. The assumption of semisimplicity implies that $C_n({\frak g}, V)$
is a direct sum of representations of $\frak h$. The trivial representation
corresponds to the module  of coinvariants $C_n({\frak g}, V)_{\frak h}$.
Let $L_n$ be the sum of all other components.

Since $\frak h$ is a Lie subalgebra, its action is compatible with the
differential $d$ on the standard complex which means that one has the
following direct sum decomposition of complexes
$$
C_\bullet({\frak g},V) = C_\bullet({\frak g}, V)_{\frak h} \oplus L_\bullet
$$
We need to show that $L_\bullet$ is an acyclic complex. To prove this
we consider the action of the  Lie algebra $\frak h$ on its  homology.
From our assumptions follows that $L_\bullet$ is a direct sum of nontrivial
${\frak h}$-modules which implies that the corresponding components in $H_\bullet(L_\bullet)$
are  also nontrivial ${\frak h}$-modules. But on the other hand one knows that 
$H_\bullet({\frak
g},V)$  is a trivial ${\frak g}$-module, and thus a trivial ${\frak
h}$-module. This implies that $H_\bullet(L_\bullet)$ has to be zero.
\qed
\bcorollary
If $A$ is a unital algebra then the  complex $C_\bullet({\frak gl}_r(A))$ is
quasi-isomorphic to the complex $C_\bullet({\frak gl}_r(A))_{{\frak gl}_r(k)}$.
\ecorollary
Proof. Apply Proposition \ref{reductive2} to the case when ${\frak g}={\frak
gl}_r(A)$, ${\frak h}= {\frak gl}_r(k)$. \qed

With theses results at hand, we can compute the coinvariant complex
explicitly. We first note that there is the following isomorphism
$$
\tp{{\frak gl}_r(A)}{n} = \tp{{\frak gl}_r(k)}{n}\otimes \tp{A}{n}.
$$
which implies that
$$
(\tp{{\frak gl}_r(A)}{n} )_{{\frak g}} = \left(\tp{{\frak gl}_r(k)}{n}
\right)_{{\frak
g}}\otimes \tp{A}{n}.
$$
This is turn gives
$$
(\tp{{\frak gl}_r(A)}{n} )_{{\frak g}} = k[\Sigma_n]\otimes \tp{A}{n}.
$$
where  we use Proposition \ref{basis} to show that  
$\left( \tp{{\frak gl}_r(k)}{n}\right)_{{\frak g}} =
k[\Sigma_n]$, for $r\geq n$. 
We now apply these considerations to the standard Lie complex
$C_\bullet({\frak gl}_r(A))$. We write
$$
C_n({\frak gl}_r(A))_{{\frak g}} = ((\tp{{\frak
gl}_r(A)}{n})_{\Sigma_n})_{{\frak g}} =
(k[\Sigma_n]\otimes\tp{A}{n})_{\Sigma_n}.
$$
where we use that the actions of ${\frak g}$ and the permutation group
$\Sigma_n$  commute. Thus we have the following statement.
\blemma
One has the following isomorphism:
$$
C_n({\frak gl}_r(A))_{{\frak g}} = (k[\Sigma_n]\otimes\tp{A}{n})_{\Sigma_n}.
$$
for any $r\leq n$.
\elemma
There is the following description of the coproduct on 
$C_\bullet({\frak gl}_r(A))_{{\frak g}}$ 
in terms of this isomorphism. 
Let $\sigma\in \Sigma_n$ and let $(I,J)$ be an ordered partition 
of the set $\{1, \dots, n\}$ such that $\sigma(I)\subset I$ 
and $\sigma(J)\subset J$. Let $\sigma_I$ ($\sigma_J$, respectively), 
denote the restriction of $\sigma$ to $I$ ($J$). 
Then for any $(a_1, \dots, a_n)\in \tp{A}{n}$ we have
$$
\triangle(\sigma \otimes (\astr{1}{n})) = \sum _{I,J}(\sigma_I\otimes
a_I)\otimes (\sigma_J\otimes a_J)
$$
This coproduct makes the above isomorphism an isomorphism of 
coalgebras. A primitive element with respect to this coproduct 
will correspond to a permutation $\sigma$ whose only invariant subset
of $\{1, \dots, n\}$ is the whole set, i.e., a cyclic permutation. 
 Let $U_n$ be the conjugacy class of the cycle $\tau = (1,
\dots, n)$.
\bproposition {\em [Loday-Quillen]}
The primitive part of the graded vector space $C_\bullet({\frak
gl}(A))_{{\frak g}}$ is the graded space 
$C^\lambda_\bullet(A)$ of cyclic chains on $A$.
\eproposition
Proof. From the remark preceding this lemma follows that the primitive 
part of $C_\bullet({\frak g}_r(A))_{\frak gl}$ is
$$
{\rm Prim}\, C_\bullet({\frak gl}_r(A))_{\frak g} = \bigoplus_{n\geq 1}(k[U_n]\otimes
\tp{A}{n})_{\Sigma_n}
$$
Let us assume that the cyclic group $\bz_n$ is embedded in the group of
permutations $\Sigma_n$ by sending the canonical generator to the
cyclic permutation $\tau$. There is a bijection between $U_n$, regarded as a
$\Sigma_n$-space and the quotient space $\Sigma_n/\bz_n$. The bijection is
given by  $\sigma\tau\sigma^{-1} \mapsto [\sigma]\in \Sigma_n/\bz_n$.
It follows that the $\Sigma_n$-module $k[U_n]$ is isomorphic to a module
induced from the trivial $\bz_n$-module $k$ by the inclusion $\bz_n \ra
\Sigma_n$. This means that we have the following isomorphism
$$
\begin{array}{rcl}
(k[U_n] \otimes \tp{A}{n})_{\Sigma_n} & \simeq & k\otimes_{k[\Sigma_n]}
k[\Sigma_n/\bz_n]\otimes \tp{A}{n}\\
& = & k\otimes _{k[{\bf Z}_n]}(k \otimes \tp{A}{n})\\
& = & (\tp{A}{n})_{{\bf Z}_n}
\end{array}
$$
We claim that the space $(\tp{A}{n})_{{\bf Z}_n}$ is the space of 
cyclic chains. To check this, 
notice that the generator of $\bz_n$ acts on $(\tp{A}{n})_{{\bf Z}_n}$
by $\sgn(\tau)\tau$. This is the same as the action of the operator
$\lambda$ which we used to define the cyclic complex. This gives
$$
(\tp{A}{n})_{{\bf Z}_n} = \tp{A}{n}/(1-\lambda) = C^\lambda_n(A)
$$
\qed

We have thus  proved
that, as a graded vector space, the space of primitive coinvariant 
chains on ${\frak gl}_r(A)$ is the same as the space of cyclic chains. Now we
need to take care of the differentials to make sure that this identification
is a map of complexes and therefore induces an isomorphism on the
level of homology. This comes down to the following simple lemma.
\blemma
The isomorphism of graded vector spaces
$$
(\tp{A}{n})_{{\bf Z}_n} = C^\lambda_n(A)
$$
is a chain map.
\elemma
Proof. Calculation. \qed
\btheorem {\em [Loday-Quillen]}
The primitive part of the Lie algebra homology $H_\bullet (\gl(A))$
is the cyclic homology  $HC_\bullet(A)$  of the algebra $A$. More precisely,
one has the following isomorphism
$$
{\rm Prim}\, H_\bullet (\gl(A)) \cong HC_{\bullet -1} (A)
$$
\etheorem

\subsection{Excision in cyclic homology}\label{Section-excision}

In this section we prove the excision property of cyclic homology, which 
was first demonstrated by Wodzicki \cite{Wodzicki}. 
We now want to study a relation between an  excision property
associated with an algebra extension $A=R/I$ and the notion of
$H$-unitality. We recall that the ideal $I$ in the
algebra extension $A=R/I$, regarded as a nonunital algebra, has the
excision property if there exists the following long exact sequence
in cyclic homology.
$$
\ra HC_{q+1}(A) \ra HC_q(I) \ra HC_q(R) \ra HC_q(A) \ra
$$
Wodzicki shows that not all ideals $I$ have this property, but only 
those that are $H$-unital (homologically unital), by which we mean the
following. 
\bdefinition
A non-unital algebra $I$ is called {\em H-unital} (homologically unital) if 
and only if the complex $(B(I), b')$ (see Remark \ref{b'-complex}) is acyclic. 
\edefinition
\btheorem\label{excisionthm}{\em [Wodzicki]}
Let $0\ra I \ra R \ra A\ra 0$ be a short exact sequence of algebras, 
where the kernel $I$ is an $H$-unital algebra.  Then there is the following long 
exact sequence in cyclic homology
\beqn
\ra HC_{q+1}(A) \ra HC_q(I) \ra HC_q(R) \ra HC_q(A) \ra .
\eeqn
\etheorem
The idea of proof that we present here is quite simple \cite{ja1}.   
To the short exact sequence
$$
0\ra I \arr{\imath} R \arr{\pi}  A\ra 0
$$
we can associate a sequence of cyclic complexes
$$
C^\lambda_\bullet (I) \arr{\imath_*} C^\lambda_\bullet(R) \arr{\pi_*}
 C^\lambda_\bullet (A).
$$
This sequence is not exact. However, we shall show that, 
from the point of view of homology, this sequence behaves as if it were
exact. A key ingredient in this approach is the theorem by Goodwillie, 
which describes the cyclic complex of the semidirect product algebra \cite{Good2,ft-LNM}.

Let  $A$ be an algebra and let   $M$  be a bimodule over  $A$.
We form the semidirect product algebra $A\oplus M$ in which $M$ is assumed to
be  an ideal of  square zero. Thus the  product in this algebra is given 
by the following formula
$$
(a, m)(a', m') = (aa', am' + ma')
$$
There is a  natural grading on $A\oplus M$  which places  $A$ in degree zero
and $M$ in degree one.  This grading
induces a grading on the cyclic complex $C^\lambda_\bullet(A\oplus  M)$
$$
 C^\lambda_\bullet(A\oplus M) = C^\lambda_\bullet(A) \oplus \bigoplus_{p\geq 1}
 C^\lambda_\bullet
(A\oplus M)(p)
$$
where the summands on the right are subcomplexes and where
$C^\lambda_\bullet(A\oplus M)(p)$ denotes the part of the cyclic complex
\mbox{$C^\lambda_\bullet(A\oplus M)$} which is of degree $p$ in $M$.   In degree
 $n$ the
space $C^\lambda_n(A\oplus M)(p)$ is spanned by the symbols $(a_0, \dots , a_n)$
where each $a_i$ is homogeneous in $A\oplus M$ and $\sum\deg\; a_i=p$.

Let us define the $n$-fold cyclic tensor product 
$(M\otimes_A)^n$ by 
$$
(M\otimes _A)^n = M\otimes_A\stackrel{n}{\cdots} \otimes _A M\otimes _A = 
(M\otimes_A\stackrel{n}{\cdots} \otimes _A M)\otimes _{A\otimes A^{op}}A
$$
where $A^{op}$ is the \lq opposite' algebra.
This space carries two natural actions 
of the  symmetric group of order $n$.
 First we let
$$
\sigma\point{m_1}{m_n}= (m_n, m_1, \dots, m_{n-1})
$$
and the second action, which we used above in the definition of the 
cyclic complex,  is obtained by letting the generator act as
$\lambda =(-1)^{n-1}\sigma $.
We shall denote the corresponding quotient spaces  of these actions by
$(M\otimes_A)^n_{\sigma}$  and $(M\otimes_A)^n_{\lambda}$, respectively.

\mbox{}

Let us now assume that $A$ is a nonunital algebra and 
let us denote by $\tilde{A}$ the augmented algebra $\tilde{A} = \bc\semiright A$.
We recall that the algebra of differential forms $\Omega \tilde{A}$ 
is an augmented DG algebra:
$$
\Omega \tilde{A} = \bc \oplus\bar{\Omega}\tilde{A}
$$
The difference between 
$\Omega\tilde{A}$ and $\bar{\Omega}{\tilde{A}}$ is in degree zero, where
$\Omega^0\tilde{A} = \bc \semiright A$ and $\bar{\Omega}^0\tilde{A} = A$. 

Let $\calc(A)= \bar{\Omega}\tilde{A}\otimes \tilde{A}$, so that in degree 
$n$ we find 
$$
\calc_n(A) = \tilde{A} \otimes \tp{A}{n}\otimes \tilde{A}
$$
We shall consider $\calc(A)$ as a complex equipped with the differential $b'$
\ref{b'-differential}.

Let $M$ be an $A$-bimodule. Then the tensor product $M\otimes _{\tilde{A}}
\calc(A)$ is a complex of $A$-bimodules. We  extend  the definition
 of the cyclic tensor product to form the complex
$$
(M\otimes_{\tilde{A}} \calc(A) \otimes_{\tilde{A}})^p
$$
whose typical generator in degree $p+n$ is  of the form
$$
(m_1, a_1, \dots, a_i, m_2, a_{i+1}, \dots, a_j, m_3,a_{j+1} \dots,m_p ,\dots
,a_n)
$$
for $m_i\in M$ and $a_j\in A$.
Again,
there are two  actions of  the cyclic group  order $p$
that one can introduce in this complex.
The first is
given by an extension of the action of the generator
$\sigma$ to complexes, and the other action
is defined with the use
of the generator $\lambda$. The  action of the latter operator
in degree $n$  can be described explicitly by
$$
\lambda(\underbrace{m_1,\dots, a_i, \dots,}_{n-l-1\;{\rm terms}}
 m_p, s_1, \dots, s_l )=
(-1)^{(n-l-1)(l+1)} (m_p, s_1, \dots, s_l, m_1, \dots \;\;)
$$
for $r_i, s_j\in A$ and
$m_k\in
M$.

Let us now consider the space $( M[1] \otimes_{\tilde{A}} \calc (A) \otimes _{\tilde{A}})
_\lambda^p$  of
coinvariants with respect to the action of the cyclic group of order $p$.
Here we treat the $A$-module $M$ as a complex concentrated in degree zero,
so that $ M[1]$ is of degree 1. There is  a canonical map which sends
a generator
$$
(m_1, x_1, \dots ,x_{i_1}, m_2, x_{i_1+1}, \dots, m_p, x_{i_{p-1}+1}, \dots,
x_l)\in (M[1]\otimes_{\tilde{A}} \calc (A)\otimes_{\tilde{A}})^p
$$
in degree $p+l$ to itself, regarded as a cyclic chain of degree $p+l-1$.
This way one defines a map
$$
(M[1]\otimes _{\tilde{A}} \calc(A) \otimes _{\tilde{A}} )^p_{\sigma}\lrar
 C^\lambda_\bullet(A\oplus M)(p)[1].
$$
We note here that since we work with complexes, the action of $\sigma$
 involves signs. Since we have
$$
(M[1]\otimes _{\tilde{A}} \calc(A) \otimes _{\tilde{A}})^p 
=(M \otimes _{\tilde{A}}
\calc(A))[p]
$$
and the cyclic group acts on $\tp{\bc[1]}{p}$ via the sign of
a permutation,
we have in fact defined a map
$$
( M\otimes_{\tilde{A}}\calc(A)\otimes _{\tilde{A}})^p_{\lambda}[p-1]
 \lrar
C^\lambda_\bullet(A\oplus M)(p).
$$
The reason that $[p-1]$ rather than $[p]$ appears on the left
is the degree convention in the cyclic complex.
We can now state the following theorem of Goodwillie \cite{Good2,Q:Extensions}.
\btheorem\label{Goodtheorem}
The canonical map
$$
( M\otimes_{\tilde{A}}\calc(A)\otimes _{\tilde{A}})^p_{\lambda} [p-1]\lrar
C^\lambda_\bullet(A\oplus M)(p)
$$
is an isomorphism of complexes for $p\geq 1$.
\etheorem
\bremark
If $R$ is a unital algebra, we can replace in the above formulae
$\tilde{A}$ by $R$, and  the nonunital algebra $A$ by $\bar{R}$ to get analogous
statements in the unital case. Both cases have been discussed in \cite{Q:Extensions}.
\eremark

\mbox{}

In the case we are interested in, we consider an $H$-unital ideal $I$  in 
an algebra $R$. We want to  treat $R$ as a filtered algebra, and 
so we shall digress a little to recall what we  need from the 
theory of filtered algebras. 

First, we say that  a vector space  $V$ is equipped with an increasing
filtration if there is a family of subspaces $F_pV$, $p\in \bz$  such that
$$
\cdots\subset F_{p-1} V\subset F_p V\subset\cdots   V.
$$
Associated to this data there is a graded vector space which we shall
denote by \mbox{$gr\, V=\bigoplus _p gr_p V$} where
$$
gr_p  V = F_pV/F_{p-1}V.
$$
If $V$ and $W$ are two filtered vector spaces, their tensor product $V\otimes
W$ is also a filtered vector space with filtration defined by
$$
F_p(V\otimes W)= \sum_{i+j=p}F_i V\otimes F_j W \subset V\otimes W.
$$
Furthermore, there is a canonical isomorphism of graded   vector spaces
\beqn\label{v-tens-w}
gr\, (V\otimes W) = (gr\, V) \otimes (gr\, W).
\eeqn
Using this general principle  we see that if $V$ is a filtered
vector space then  the filtration of $V$ induces filtrations
of multiple tensor products $\tp{V}{n}$ for any positive integer $n$.
We also have the following analogue of the isomorphism (\ref{v-tens-w})
$$
gr\, V^{\otimes n} = (gr\, V)^{\otimes n}.
$$

Let $A$ be an algebra. A filtration of $A$ is a family $F_{p}A$, $p\in \bz$,
of $k$-submodules of $A$
which gives a filtration of $A$ as a vector space and which satisfies
$$
F_p A\cdot F_qA \subset F_{p +q}A .
$$
The only algebra filtration that we shall see in this chapter will satisfy the
following two conditions:
 $F_{-1}A=0$ and
$\bigcup_p F_pA = A$. 

The product $A\otimes A \rightarrow A$ is thus compatible with filtrations and
induces a graded algebra structure on $gr\, A$, where the grading
is with respect to   nonnegative integers. Moreover, 
 an algebra filtration induces a
filtration by subcomplexes of the Hochschild  complex $C_\bullet(A,A)=C_\bullet(A)$,
and  the cyclic complex  $C^\lambda_{\bullet}(A)$.
Thus in the case of, say, the cyclic complex $C^\lambda_\bullet(A)$ we have
that the space $F_p C^\lambda_n(A)$ is spanned by the elements
$(\astr{0}{n})$, where $a_i\in F_{p_{i}}A$ and $\sum p_i \leq p$.

The grading on the algebra $gr A$ induces a grading $C_\bullet(gr A) =
\bigoplus_pC_\bullet(A)(p)$ of the Hochschild complex by subcomplexes.
A generalization of the isomorphism (\ref{v-tens-w}) to the algebraic
situation states that there is a canonical isomorphism of complexes
$$
gr_p C_\bullet(A) = C_\bullet(gr\, A)(p) .
$$
A similar result holds for the cyclic complex.

\mbox{}

We  now want to apply this general theory together with  Goodwillie's Theorem
\ref{Goodtheorem}  to the case of an algebra
extension $A = R/I$. We regard $R$  as a filtered
algebra, with  filtration given by
$$
F_0 R = I , \;\;\; F_p R = R, \;\; p\geq 1.
$$
It is quite clear that the  associated graded algebra is then just
$$
gr(R)= I\oplus R/I = I\oplus M.
$$
where we  denote $M=R/I$ when we want to stress the $I$-bimodule structure
of the algebra  $A=R/I$. The  algebra  $gr(R)$ has a very simple product
induced from the product on $R$: $I\cdot I\subset I$, $M\cdot I = I\cdot M =
M^2 = 0$. In particular, this means that $I$ acts trivially on $M$.

Let us form the complex 
$M\otimes_{\tilde{I}} \calc(I) \otimes _{\tilde{I}} $
which in degree $n$ is given by 
 $$
\left(M\otimes _{\tilde{I}} \calc(I)\otimes_{\tilde{I}}\right)_n
= M\otimes \tp{I}{n}.
$$
This complex is equipped with the differential $b$ 
given by the usual  formula
\begin{eqnarray*}
b(m, \astr{1}{n}) & = & (ma_1 , \astr{2}{n}) \\
&&\rule{0mm}{5mm}\mbox{}+ \displaystyle{\sum_{i=1}^{n-1}} (-1)^i(m, a_1, \dots,
 a_ia_{i+1}, \dots , a_n)\\
&&\mbox{} + (-1)^n(a_n m , \astr{1}{n-1}).
\end{eqnarray*}
Taking into account  the product on the associated graded algebra $gr(R)$
we see that this differential is simply
$$
b(m, \str{r}{1}{n}) = m\otimes (-b')(\str{r}{1}{n})
$$
since $r_nm= mr_1 = 0$ in this case. Thus we have the following
isomorphism of complexes
\beqn\label{canonical}
M\otimes_{\tilde{I}} \calc(I) \otimes_{\tilde{I}}
 \simeq M\otimes \left (\tp{I}{\bullet},
b'\right).
\eeqn
The complex $(\tp{I}{ \bullet}, b')$ is just the complex $(B(I), b')$ 
of the ideal $I$ considered as a nonunital algebra (compare \ref{b'-complex}). 

To describe the cyclic complex of the algebra $gr(R)$ we note that
from the above and Theorem \ref{Goodtheorem}
$$
C^\lambda_\bullet(I\oplus M) = C^\lambda_\bullet(I)\, \oplus \, \left(M \otimes
B(I)\right)
\, \oplus
(M\otimes_{\tilde{I}} \calc(I) \otimes _{\tilde{I}})_\lambda[1] \oplus \dots.
$$
The general term in filtration degree $p$ is
$$
C^\lambda_\bullet(I\oplus M)(p) = (M\otimes_{\tilde{I}}\calc(I)
\otimes _{\tilde{I}})^{p}_\lambda[p-1]
$$
which we can write using (\ref{canonical})  as follows
$$
\left(M \otimes_{\tilde{I}}\calc(I)\otimes_{\tilde{I}}\right)^{p} =
\left( M\otimes B(I) \otimes \right)^p
$$
Now assume that the ideal $I$, considered as a nonunital algebra, is $H$-unital.
This  means that the bar construction $B(I)$  in the formula above is
a contractible complex.  Thus the following canonical map is a quasi-isomorphism
$$
C^\lambda_\bullet(I\oplus M){(p)} \simeq (M^{\otimes p})_\lambda[p-1]
$$
for $p\geq 1$. When $p=0$ we have that
$$
C^\lambda_\bullet (I\oplus M)(0) = C^\lambda_\bullet (I).
$$

Finally,  let us introduce the following filtration on the cyclic complex
$C^\lambda_\bullet(M)$, where $M=R/I$. We put $F_pC^\lambda_n(M) = 0$ for
 $p\leq n$, and
$F_pC^\lambda_n(M) = C^\lambda_n(M)$, $p\geq n+1$. We note that the  canonical
 surjection
$\pi: C^\lambda_\bullet(R) \rightarrow C^\lambda_\bullet(M)$ carries
 $F_pC^\lambda_\bullet(R)$
onto $F_pC^\lambda_\bullet(M)$.

\blemma\label{basiclemma}
Let  $A=R/I$ be an algebra extension  where we assume that the ideal
$I$ is $H$-unital. Then
 for $p\geq 1$ the map
$$
gr_p C^\lambda_\bullet(R)\rightarrow gr_p C^\lambda_\bullet (A)
$$
induced by the canonical  projection  $\pi: R\ra A$ is a quasi-isomorphism.
\elemma
Proof. We denote again by $M$ the $I$-bimodule $A= R/I$ in the semidirect
 product
algebra $gr(R) = I\oplus M$.  Using that
$$
C^\lambda_\bullet(I\oplus M)(p) = C^\lambda_\bullet(gr(R))(p)=
gr_p C^\lambda_\bullet(R)
$$
we can make the following identifications
$$
 (M \otimes {\cal C}(I) \otimes)
^{p}_{\lambda}[p-1]
\stackrel{\sim}{\rightarrow} C^\lambda_\bullet(gr\, R)(p) =
 gr_pC^\lambda_\bullet(R)
$$
for $p\geq 1$.
For $p=0$, $C^\lambda_\bullet (gr\, R)(0)= C^\lambda_\bullet (I)$.

From the above analysis follows that, when $I$ is an $H$-unital algebra,
the nonzero homology of $gr_pC^\lambda_\bullet(R)$ is simply $M^{\otimes p}_{\lambda}$
located 
in degree $p-1$. This is the same as the homology of $gr_pC^\lambda_\bullet(M)$.
 Furthermore,  for $p=0$, the nonzero
 homology of ${\cal B}(I)$ is $I$.
The lemma then follows from the fact that the isomorphism of Theorem
\ref{Goodtheorem}, the augmentation ${\cal C}(I) \rightarrow I$ and
$\pi$ are compatible with the filtration on $R$.
\qed

\bproposition\label{lastprop}
The projection map $\pi$ induces a quasi-isomorphism
$$
C^\lambda_\bullet(R) / C^\lambda_\bullet(I) \ra C^\lambda_\bullet(R/I)
$$
when the ideal $I$ is an $H$-unital algebra.
\eproposition
Proof.
The lemma above implies by induction that $\pi$ induces a quasi-isomorphism
$$
F_pC^\lambda_\bullet(R)/F_0C^\lambda_\bullet(R) \rightarrow
 F_pC^\lambda_\bullet(R/I)
$$
for all $p$. Since  $F_0 C^\lambda_\bullet(R)=C^\lambda_\bullet(I)$,
 and the  filtration of $R$ is exhaustive,
we have the required result.  \qed

Theorem \ref{excisionthm} now follows from Proposition \ref{lastprop} if we
apply the homology functor to the exact sequence of complexes
$$
0\lrar C^\lambda_\bullet (I) \lrar C^\lambda_\bullet (R) \lrar
C^\lambda_\bullet(R)/C^\lambda_\bullet(I)\lrar 0.
$$
%% excision in algebraic $K$-theory
\subsection{Excision in algebraic $K$-theory}
From the point of view of calculations in the algebraic $K$-theory it is 
important to know how to compare the $K$-groups $K_n(A)$ and  $K_n(A/I)$, 
where $I$ is an ideal in $A$. We have treated the corresponding problem 
for cyclic homology in the previous section. It is possible to define 
abstractly the relative $K$-groups  $K_n(A,I)$ 
so that they fit in the exact sequence
$$
\ra K_{n+1} (A/I) \ra K_n(A,I) \ra K_n(A) \ra K_n(A/I) \ra \cdots \ra K_0(A/I)
$$
We have defined relative $K$-groups for $n=0,1$ in Section \ref{relative-K-section}, 
and we shall indicate how to do that for positive $n$ in a moment. 
The best we can hope for is that if the relative $K$-groups depend only on the 
ideal $I$, then we could substitute the $K$-theory $K(I)$ in place of 
the relative groups $K(A,I)$. Recall from  \ref{nonunital-I}, Section \ref{relative-K-section}
that to define $K_n(I)$ we treat $I$ as a nonunital 
algebra and consider its augmentation $\tilde{I} = \bz \semiright I$
where the product is defined  by 
$(n, u ) ( m, v ) = (nm, nv + mu + uv) $
Then we define 
$$
K_n(I) = \ker[ K_n(\tilde{I}) \ra K_n(\bz)]
$$
This definition should be compatible with the definition of relative $K$-theory in 
the sense that $K_n(I) = K_n(\tilde{I}, I )$. A natural question at this stage
is the following: Is the natural map $K_n(I) \ra K_n(A,I)$ an isomorphism? 
We shall say that an ideal $I$ for which this isomorphism takes place
has the {\em excision property}. 
It has been shown by H. Bass that this is indeed so in degree $0$, and 
we have checked this result in \ref{I-tilde}. It is known that 
for arbitrary ideals this statement is false for $n>0$. It is clear 
then that one needs to find a class of ideals that would satisfy the 
excision property. 
Suslin and Wodzicki have formulated a simple
necessary and sufficient condition for an ideal $I$ to have the excision 
property in the algebraic $K$-theory \cite{SW}. Their result may be stated 
as follows. For any abelian group $M$ we 
shall denote by $M_{\bf Q} = M\otimes _\bz \bq$. 
\btheorem {\em [Suslin-Wodzicki]}
For any nonunital ring $I$ the following conditions are equivalent. 
\begin{enumerate}
\item $I$ has the excision property with respect to $K_{\bf Q}$. 
\item $I_{\bf Q}$ is $H$-unital. 
\end{enumerate}
\etheorem
\bcorollary
If $I$ is a $\bq$-algebra (not necessarily unital) then the following 
conditions are equivalent. 
\begin{enumerate}
\item $I$ has the excision property with respect to $K$. 
\item $I$ is $H$-unital. 
\end{enumerate}
\ecorollary
Among the many important applications of this result there is a proof of the 
following Karoubi conjecture. 
\btheorem 
Let us denote by $\bk$ the $C^*$-algebra of compact operators in a separable 
Hilbert space. Then there is a canonical isomorphism for every $n\geq 0$
$$
K_n^{\rm alg}(A\hat{\otimes}_\pi\bk ) \simeq K^{{\rm top}}_n(A\hat{\otimes }_\pi
\bk)
$$
\etheorem
In the above formula, $\hat{\otimes}_\pi$ denotes the projective completed
tensor product in category of $C^*$-algebras. Moreover, $ K^{{\rm top}}_n$
denotes the topological $K$-theory, i.e. the $K$-theory defined in 
the category of $C^*$-algebras as in the first chapter. 
This  theorem then implies the following  corollary proved by Wodzicki. 
\btheorem
For any $C^*$-algebra $A$ the algebraic $K$-theory 
$K_n^{alg}(A\otimes \bk)$ 
is periodic with period $2$. 
\etheorem

To explain the statement of the theorem, recall that
$$
K_n(A) = \pi_n(BGL(A)^+)
$$
for $n\geq 1$. 
The homology groups of the general linear group $GL(A)$ are related in an 
important way to the algebraic $K$-theory. 
 Consider the homology 
$H_*(GL(A), \bq) = \bigoplus _{n\geq 0} H_n (GL(A), \bq)$ with rational 
coefficients which is  a commutative and cocommutative Hopf algebra and 
so it is completely determined by its primitive part. It turns out that 
the primitive part of this algebra is precisely the algebraic 
$K$-theory of $A$:  
$$
{\rm Prim}\, H_*(GL(A), \bq)\simeq \bigoplus_{n\geq 1} K_n(A)_{\bf Q}
$$
This theorem is  the prototype of the Loday-Quillen 
theorem in cyclic homology. 

Let $I$ be a two-sided ideal in the algebra $A$. The canonical surjection 
$A \ra A/I $ induces a group homomorphism $GL(A) \ra GL(A/I)$. Let us 
denote by $\overline{GL}(A/I)$ the image of  $GL(A)$ inside 
$GL(A/I)$ under this map. There is an induced map 
$BGL(A)^+ \ra B\overline{GL}(A/I)^+$ with a connected homotopy fibre 
$F(A,I)$. 
\bdefinition
The relative $K$-theory groups are defined by 
$$
K_n(A,I) = \pi_n(F(A,I)
$$
for $n\geq 1$. 
\edefinition
It follows that the long exact sequence involving relative $K$-theory stated 
at the beginning of this section is a long exact sequence associated with 
a fibration.

The main theorem proved by Wodzicki and Suslin is the following \cite{SW}. 
\btheorem
Let $I$ be a (not necessarily unital) algebra over $\bq$. Then the following
conditions are equivalent. 
\begin{enumerate}
\item $I$ is $H$-unital. 
\item $I$ has the excision property in cyclic homology. 
\item $I$ has the excision property with respect to $K_{\bf Q}$. 
\end{enumerate}
\etheorem
The equivalence of conditions $1$ and $2$ was 
proved by Wodzicki in \cite{Wodzicki}, who also showed that if $I$ has the excision 
property with respect to $K_{\bf Q}$ then it is $H$-unital. Thus the main result of
Suslin and Wodzicki was to prove that if $I$ is $H$-unital then it has the excision 
property in $K_{\bf Q}$. We refer to the paper \cite{SW} and to Loday's 
exposition \cite{L2}.

\vfill\eject

\section{Further examples}
\subsection{Derivations and cyclic homology}
We shall now devote some time to  calculations and examples. We begin
with another 
application of  Theorem  \ref{Goodtheorem} which will be used to show that the
derivations of an algebra $A$ act trivially on the periodic cyclic
homology of $A$. 
First let us recall that  a  derivation $D$ of an  algebra $A$
is a linear map $D:A\ra A$ satisfying the Leibniz rule
$$
D(ab)= D(a) b+ a D(b)
$$
for all $a,b$ in $A$. This formula uses the fact that the algebra $A$ is a
 bimodule
over itself and  one  defines a derivation with values in a bimodule $M$ over
 $A$
to be a linear map $D: A\ra M$ satisfying the same Leibniz rule.
As we have mentioned in \ref{multiple-tensors}, 
there is a universal model for derivations of $A$ with values in 
an $A$-bimodule, which is provided by the bimodule of $1$-forms
$\Omega^1 A$ and the derivation $d: A\ra \Omega^1A$. 
Recall also, that $\Omega^1$ enters the exact
sequence
\beqn\label{exact}
0\ra \Omega^1 \arr{i} A\otimes A \arr{m} A \ra 0
\eeqn
where the inclusion $i$ is given by 
$$
i(a_0da_1) = a_0(a_1\otimes 1 - 1 \otimes a_1)
$$

The short exact sequence (\ref{exact}) gives rise to a corresponding short
exact sequence of Hochschild chain  complexes
\beqn  \label{compl-seq}
0\lrar C_\bullet(A, \Omega^1 A)\lrar C_\bullet(A,A\otimes A)\arr{m\otimes 1}
 C_\bullet(A,A)\lrar 0
\eeqn
which leads to a long exact sequence in Hochschild homology
\beqn\label{longsequence}
\ra H_{n}(A,A\otimes A)\ra
H_n(A,A)\arr{\delta} H_{n-1}(A,\Omega ^1
A)\ra H_{n-1}(A,A\otimes A)\ra
\eeqn
where $\delta$ denotes the connecting homomorphism.
\bproposition\label{boundary}
The boundary homomorphism
$$
\delta: H_{\bullet +1}(A,A)\lrar H_{\bullet}(A,\Omega^1 A)
$$
is an isomorphism in positive degrees. In degree 0
this map is injective.
\eproposition
Proof. We shall give an explicit formula for the connecting homomorphism
$\delta$. For this purpose let us  define a lifting $s: C_p(A,
A)\ra C_p(A, A\otimes A)$
by
$$
(a_0, a_1, \dots , a_p)\ra (a_0\otimes 1, a_1, \dots , a_p).
$$
Now for any $p$-chain $(\astr{0}{p})$ we have
$$
\begin{array}{rcl}
bs(\astr{0}{p})& = & (a_0\otimes a_1)\otimes
 (\astr{2}{p})+\sum
(-1)^i(a_0\otimes 1)\otimes (\;\dots, a_ia_{i+1}, \dots \;)\\
&& \mbox{}+ (-1)^p(a_pa_0\otimes 1)\otimes (a_2, \dots,
 a_{p-1})\\
\rule{0mm}{8mm}
sb(\astr{0}{p})& = & (a_0a_1\otimes 1)\otimes
 (\astr{2}{p}) +\sum
(-1)^i(a_0\otimes 1)\otimes(\;\dots, a_i a_{i+1},\dots\;)\\
&&\mbox{}+ (-1)^p(a_pa_0\otimes 1)\otimes
 (\astr{2}{p-1})
\end{array}
$$

The difference of the expressions calculated above is
\beqn\label{bs-sb}
\begin{array}{rcl}
(bs-sb)(\astr{0}{p})& = & (a_0a_1\otimes 1)\otimes
 (\astr{2}{p})\\
&&\mbox{} -(a_0 \otimes a_1)\otimes (\astr{2}{p}).
\end{array}
\eeqn
It is clear that the map $m\otimes \tp{1}{p-1}:C_{p-1}(A,
 A\otimes A)\ra
C_{p-1}(A, A)$ in the exact sequence of chain complexes
 (\ref{compl-seq})
sends this expression to zero; from the exactness of the sequence
(\ref{compl-seq}) then follows that there exists a chain in
 $C_{p-1}(A,
\Omega^1 A)$ whose image under this map $(bs-sb)(\astr{0}{p})$. To identify
 this chain we
use the definition of the bimodule $\Omega^1 A$.
We have
$$
\imath(a_0da_1)= b'(a_0\otimes a_1\otimes 1)= a_0a_1\otimes 1-
 a_0\otimes a_1 .
$$
This together with (\ref{bs-sb}) gives the boundary map $\delta$:
$$\begin{array}{rcl}
\delta: C_p(A, A) &\lrar & C_{p-1} (A, \Omega^1 A)\\
a_0\otimes (\astr{1}{p}) &\mapsto & a_0
 da_1\otimes(\astr{2}{p}).
\end{array}
$$

Now we need to show that this map induces an isomorphism in homology in positive
degrees. For this we define a map $f:C_\bullet(A, A\otimes A)\ra
 C_{\bullet +1}(A,A)$
by
$$
f:(x\otimes y)\otimes(\astr{1}{p})\ra (y,\astr{1}{p}, x).
$$
We then calculate
\begin{eqnarray*}
fb((x\otimes y)\otimes (\astr{1}{p})) &=& (ya_1, \dots , a_p,
 x)\\
\rule{0mm}{5mm}
&& \mbox{}+ \sum_{i=1}^{p-1}(-1)^{i}(\;\dots ,
 a_ia_{i+1}, \dots\;)\\
&& \mbox{}+(-1)^p(y,\astr{1}{p-1},a_px)\\
&=& b'f(y, \astr{1}{p}, x).
\end{eqnarray*}
This calculation shows that the complex $C_\bullet (A, A\otimes A)$ with
differential $b$ is isomorphic to the complex $C_\bullet (A,A)$ with the
differential $b'$ with degrees  shifted by one. But the latter complex is
contractible since we assume that the algebra $A$ is unital, so we see that
the homology of the complex \mbox{$C_\bullet (A, A\otimes A)$} is zero in
positive degrees. This means that in the homology long exact sequence
(\ref{longsequence}) every third term is zero in positive degrees, hence the
boundary map $\delta$  indeed  gives an isomorphism
\mbox{$H_{\bullet +1}(A, A)\ra H_\bullet (A,\Omega^1 A)$}. In  lowest
 degrees the long exact
homology sequence has the following form
$$
0\lrar H^1(A,A)\lrar \Omega^1A\lrar A\lrar A/[A,A] \lrar 0
$$
which shows that in this case the map is injective, and the proposition
follows.  \qed

\mbox{}

We are now ready to relate these facts to the Connes long exact
sequence
$$
\ra H_n(A,A) \arr{I} HC_n(A) \arr{S} HC_{n-2}(A) \arr{B} H_{n-1}(A,A)\ra
$$
Let us then consider the semi-direct product algebra $A\oplus \Omega^1A$, where
we assume that the ideal $\Omega^1A$ has zero product. From the
formula of Theorem \ref{Goodtheorem},  the cyclic complex of $A\oplus \Omega^1
 A$
can be written to the first order in $\Omega^1A$
as follows
$$
C^\lambda_\bullet (A\oplus \Omega ^1 A) = C^\lambda_\bullet (A)\oplus
 C_\bullet(A, \Omega^1 A)
\oplus \cdots
$$
where  the complex $C_\bullet (A,\Om^1A)$ of degree one  in $\Om^1 A$ is
the Hochschild complex of $A$ with coefficients in $\Om^1A$.
Let $\pi$ denote the projection
$$
\pi:C^\lambda_\bullet (A\oplus \Omega^1 A)\lrar C_\bullet(A, \Omega^1 A).
$$
There is  a corresponding map from cyclic to Hochschild  homology
$$
HC_{\bullet}(A\oplus \Omega^1A) \ra H_{\bullet}(A, \Omega^1 A).
$$

Let $D:A\ra A$ be a derivation of an algebra $A$. It acts on homology
via the corresponding Lie
derivative  $L_D:C(A,A)\ra C(A,A)$ defined by
$$
L_D(a_0\otimes \cdots \otimes a_p)= \sum_{i=0}^p (a_0,
 \dots, Da_i, \dots , a_p).
$$
\bproposition
The  following diagram
$$
\begin{array}{ccc}
HC_\bullet (A) & \Ar{B} & H_{\bullet +1}(A,A)\\
\Bdal{\tilde{d}} && \Bdar{\delta} \\
HC_\bullet(A,A\oplus\Omega^1 A)& \Ar{p} & H_\bullet (A, \Omega^1 A)
\end{array}
$$
commutes. Here $\delta$ is the boundary map of Proposition \ref{boundary},
and $\tilde{d}$ is
the inclusion
map induced by the canonical map $1\oplus d:A\ra A\oplus \Omega ^1A$.
\eproposition
Proof. Let us assume that $(\astr{0}{p})\in C^\lambda_p(A)$. We have
\begin{eqnarray*}
\pi\tilde{d}(\astr{0}{p}) &=& \pi\point{a_0 +da_0}{ a_p+
 da_p}\\
\rule{0mm}{5mm}
&=& \displaystyle{\sum_{i=0}^p} (a_o, \dots , da_i , \dots , a_p) \\
\rule{0mm}{5mm}
&=& \displaystyle{\sum_{i=0}^p}  (-1)^{pi} (da_i , a_{i+1}, \dots , a_p,
 a_0, \dots , a_{i-1}).
\end{eqnarray*}
From the calculation of the boundary map $\delta $ we see that this element
comes from the chain
$$
\sum (-1)^{pi} (1, a_i , \dots, a_p , a_0, \dots ,
 a_{i-1}) = B(\astr{0}{p})
$$
which proves the proposition.\qed
\bremark
 We note  that this proposition gives a different way to find a formula for
the operator $B$, which appears in the Connes long exact sequence, and 
which we defined in Section \ref{FORMS}. 
\eremark

\btheorem\label{ls=0}
The composite map $L_D\circ S: HC_{\bullet}\ra
HC_{\bullet - 2}$ acts as the zero operator on cyclic homology. 
\etheorem
Proof. Let us consider the following diagram 
$$
\begin{array}{ccccc}
HC_{\bullet}(A) & \Ar{B} & H_{\bullet +1}(A,A) & \Ar{\delta}& H_{\bullet}(A,
 \Omega ^1 A)\\
\Bdal{L_D}&&&& \Bdar{\tilde D}\\
HC_{\bullet}(A) &&\Lla{I}&& H_{\bullet}(A,A)
\end{array}
$$
where $I$  and $B$ are maps from the long exact homology sequence 
of Connes and
$\tilde{D}: H_{\bullet}(A,\Om^1 A)\ra H_{\bullet} (A,A)$ is induced by the map
$$
\tilde D : \Omega^1A \arr{} A, \qquad a_1d a_2 \mapsto a_1Da_2.
$$
 We show that this diagram is commutative. To this end  we calculate
\begin{eqnarray*}
(\astr{0}{p}) &\stackrel{\delta B}{\mapsto} &
 \sum(-1)^{pi}(da_i, a_{i+1},
\dots, a_p, a_0, \dots, a_{i-1})\\
&\stackrel{\tilde D}{\mapsto} & \sum (-1)^{pi} (Da_i,
 a_{i+1}, \dots ,a_p, a_0, \dots,
a_{i-1})\\
&\stackrel{I}{\mapsto} & \sum (a_0, \dots Da_i, \dots a_p)\\
& =& L_D(\astr{0}{p}).
\end{eqnarray*}
Hence we can write
$$
L_D\comp S= I\tilde D\delta BS= 0
$$
since $BS=0$ from the Connes long exact sequence. This completes the proof.
\qed
\bcorollary\label{cor}
Any derivation of the algebra $A$ acts trivially on the periodic cyclic
 homology $HP_\nu(A)$.
\ecorollary
Proof. We recall that the complex which computes the periodic homology
is the inverse limit of the cyclic complex with respect to the action
of the periodicity operator $S$. \qed
\bremark
Theorem \ref{ls=0} was proved by Goodwillie in \cite{Good2} in the setting 
of his theory of cyclic objects. The proof given here is taken 
from \cite{ja1}.
\eremark

\subsection{Cyclic and reduced cyclic homology}

As another application of Goodwillie's theorem \ref{Goodtheorem} we shall 
discuss the reduced cyclic homology of a unital algebra. 
 In this case the long exact sequence derived in the
previous section becomes a long exact sequence relating cyclic and reduced
 cyclic
homology theory of an algebra.

We note first that when $A$ is unital, one has the following short exact
 sequence
\beqn
0\lrar k\lrar A \lrar \bar A \lrar 0
\eeqn
which gives the  sequence of complexes
\beqn
C^\lambda_{\bullet}(k)\arr{\imath} C^\lambda_{\bullet}(A) \arr{\pi}
 \bar{C}^\lambda
_\bullet(A)
\eeqn
where $\bar{C}^\lambda_\bullet(A)$ is the reduced cyclic complex, i.e.
the cyclic complex
of the nonunital algebra $\bar{A}$. As before, this sequence is not exact but
it has a similar property as the sequence discussed in the previous section---%
it gives an exact sequence in homology.

\mbox{}

A unital algebra can be regarded as a filtered algebra with the following simple
filtration. We assume that $F_0A = k$ and that $F_pA = A$ for all $p\geq 1$.
The associated graded algebra $gr\, A$ is in this
case just $gr\, A = k\oplus \bar{A}$. It follows from the definition of our
filtration that the ideal $\bar{A}$ has zero product.
This filtration of the algebra $A$ induces filtrations of the Hochschild and
 cyclic
complexes. We shall be particularly interested in the filtration induced on the
cyclic complex $C^\lambda_\bullet(A)$ of the algebra $A$.

Let us introduce the following filtration of the reduced cyclic complex
$\bar{C}^\lambda(A)_\bullet$. We put
$F_p\bar{C}^\lambda_n(A) = 0$ for  $p\leq n$ and
$F_p\bar{C}^\lambda_n(A) =\bar{C}^\lambda_n(A)$, $p\geq n+1$.
We note that the  canonical surjection $\pi: C^\lambda(A) _\bullet\rightarrow
\bar{C}^\lambda_\bullet(A)$
carries $F_pC^\lambda_\bullet(A)$ onto $F_p\bar{C}^\lambda_\bullet(A)$.
Then we have the following  analogue  of Lemma \ref{basiclemma}.
\blemma
 For $p\geq 1$ the map
$$
gr_p
C^\lambda_\bullet(A)
\rightarrow
gr_p\bar{C}^\lambda_\bullet(A)
$$
induced by $\pi$ is a quasi-isomorphism.
\elemma
Proof.  The associated graded algebra $gr\, A$ of the algebra $A$
given by $gr\, A
= k
\oplus \bar{A}$ is the semi direct product algebra.
From  Theorem \ref{Goodtheorem}  we have the following isomorphism
$$
\Sigma^{p-1} (\bar{A} \otimes {\cal B}(k) \otimes)
^{p}_{\lambda}
\stackrel{\sim}{\rightarrow} C^\lambda(gr\, A)(p) = gr_pC^\lambda(A)
$$
for $p\geq 1$.
For $p=0$, $C^\lambda(gr\, A )(0)= C^\lambda(k)$.
Now the homology  of $gr_pC^\lambda(A)$  consists of just  $\bar{A}^{\otimes
 p}_{\lambda}$
in degree $p-1$ which is the same as the homology of
 $gr_p\bar{C}^\lambda_\bullet
(A)$. On
the other hand,  in degree zero, the nonzero
 homology of ${\cal B}(k)$ is $k$.
The result of the  lemma follows from the fact that the three maps used here: the
 isomorphism
of Theorem \ref{Goodtheorem}, the canonical projection $\pi$ and the
 augmentation
$\calb(k) \ra k$ are compatible with the filtration on $A$.
\qed

Proceeding by induction we see that the map
$$
F_pC^\lambda_\bullet(A)/F_0C^\lambda_\bullet(A) \rightarrow
 F_p\bar{C}^\lambda(A)
_\bullet
$$
is a quasi-isomorphism for all $p$. As $F_0 C^\lambda(A)=C^\lambda(k)$ and
the increasing filtration is exhaustive, we have the following statement.
\bproposition\label{kwiz}
The projection map $\pi$ induces a quasi-isomorphism
$$
C^\lambda_\bullet(A)/C^\lambda_\bullet(k) \rightarrow \bar{C}^\lambda_\bullet(A)
$$
\eproposition

There is a short exact sequence of complexes
$$
0\lrar C^\lambda_\bullet(k) \lrar C^\lambda_\bullet(A) \lrar
 C^\lambda_\bullet(A)/
 C^\lambda_\bullet(k) \lrar 0
$$
so if we use the corresponding long exact sequence in homology and Proposition
\ref{kwiz} above we have the following.
\btheorem\label{reduced-sequence}
There is a long exact homology sequence
$$
\lrar HC_n(k) \lrar HC_n(A)\lrar \bar{H}C_n(A) \lrar HC_{n-1}(k)\lrar
$$
\etheorem

\bremark This theorem is proved in a different way in  \cite{LQ}
with the use of a short
exact sequence of reduced and normalized double complexes.
\eremark

\subsection{Canonical classes in  reduced cyclic cohomology}

It is not difficult to find a cohomology version of the long 
exact sequence derived in Theorem \ref{reduced-sequence}. 
Since the operation $V\mapsto V^*$ of taking a dual
vector space is exact, we obtain a long
 exact sequence relating cyclic and reduced
cyclic cohomology:
$$
\lrar HC^n(k)\lrar \bar{H}C^{n+1}(A)\lrar HC^{n+1}(A)\lrar
HC^{n+1}(k)\lrar
$$
Moreover this sequence is obtained from the sequence of complexes
$$
\bar C_\lambda^{{ \bullet}}(A) \arr{\pi^*}C_\lambda^{{ \bullet}}(A)
\arr{\imath^*} C_\lambda^{{ \bullet}}(k)$$
where * denotes the  transpose map on dual vector spaces and the maps $\imath$
 and
$\pi$ are the same as in the previous section. These maps have the following
properties.
It is clear that  $\imath^* \pi^* = 0$,  $\pi^*$ is  injective, and
$\imath^*$ is
surjective.
Moreover,
 the map induced by $\pi^*$ from $\bar C^\lambda(A)_\bullet$ to the kernel of
$\imath^*$ is a quasi-isomorphism.
Thus we see that, as in the dual case of homology, this sequence, although not
exact, can be regarded as an exact sequence if we are only interested in
its cohomology.

The cyclic cohomology of the field $k$ is particularly simple.
\blemma\label{k}
One has
$$
HC^n(k)=\left\{\begin{array}{cl}
k,& \quad n=2i\\
0, &\quad {\rm otherwise.}
\end{array}
\right.
$$
\elemma
Proof. We note that there are no cyclic cochains in odd degrees and that
the space of cyclic cochains in even degrees  is one dimensional. The result
then  follows from the fact that the differential $b$ is identically  zero in
 the
even degrees.\qed

Hence the cyclic cohomology of a field is generated in each even degree by 
a canonical generator. In degree $2n$ this generator  is the 
image of $ 1\in HC^0(k)$ under $S^n$.
We want to use this fact together with the cohomology long exact sequence to
describe some canonical reduced cohomology classes in
$\bar{H}C^{\bullet}(A)$. For this we shall discuss the connecting map
$HC^{2n}(k) \ra \bar HC^{2n+1}(A)$. We shall use a method similar to the
usual transgression procedure of  Chern-Weil  theory. This means
that  we want
to  construct a cochain $f\in C_\lambda^{2n}(A)$ such that
 $\imath^*f \in C_\lambda^{2n}(k)$ is a cocycle
whose class generates $HC^{2n}(k)$, and  $bf = \pi^*g$ with
$g\in \bar C_\lambda^{2n+1}(A)$. Under these assumptions  $g$ is a
reduced cyclic cocycle, whose class
$[g] \in \bar HC^{2n+1}(A)$ is the image under the above connecting map of the
generator $[\imath^*f] \in HC^{2n}(k)$.

Let us then  consider $B(A)^*$, the dual of the acyclic Hochschild complex defined in 
\ref{b'-complex}. We recall that  $B(A)^*$ is a differential graded algebra, 
with pointwise multiplication of chains and 
with differential $b'$. The cyclic norm operator $N$ gives a trace on
 $B(A)^*$ with values in $C_\lambda^{{ \bullet}}(A)$.  In
degree $n$ the operator $N$ is given by the formula
$$
N= \sum _{i=0} ^{n}\lambda^i.
$$

Let us choose
a $1$-cochain  $\rho :A\ra k$ such that  $\rho(1)=1$. As this is a linear map, 
we can consider its  curvature, which  is  a two-cochain  $\om=b'\rho - \rho^2$.
Consider the Chern-Simons form associated with $\rho$
$$
cs_{2n+1} = N \int _0^1 \rho (tb'\rho - t^2\rho^2)^{n}dt/n!.
$$
which is a cyclic cocycle of degree $2n$.
We have already encountered this definition in formula \ref{Chern-Simons} of 
Section \ref{Higher-Traces}. 

\blemma
One has the following transgression identity
$$
b\, cs_{2n+1} = ch_{2n+2}
$$
where $ ch_{2n+2}= N \omega^{n+1}/(n\!+\!1)!$ is the Chern character form
determined by the cochain $\rho$ and trace $N$.
\elemma
Proof. Let $\rho_t$ be a one-parameter family of cochains as above. Let
$\dot{\om_t}$ denote the derivative of $\om_t$ with respect to the parameter
$t$, and let $\ad\rho_t = [\rho_r, \cdot]$.  We have
$$\begin{array}{rcl}
\pr _t(\omega_t^n)& =& \displaystyle{ \sum_{i=1}^n}
\omega_t^{i-1}\dot{\om_t}\om
_t^{n-i}\\
\rule{0mm}{5mm}
& = &\displaystyle{ \sum_{i=1}^n} \omega_t^{i-1}(b'-\ad
\rho_t)\dot{\rho_t}\om _t^{n-i}\\
\rule{0mm}{5mm}
& = & (b'-\ad  \rho_t)\displaystyle{\sum_{i=1}^n
}\omega_t^{i-1}\dot{\rho_t}\om _t^{n-i}.
\end{array}
$$
Let $\mu_{n,t}= \sum_i \om_t^{i-1}\dot{\rho_t}\omega_t^{n-i}$.
We calculate
$$\begin{array}{rcl}
\pr_t (N(\omega^n_t) )&  = &  N(b'-\ad\rho_t)\mu_{n,t}\\
&=& b\,N (\mu_{n,t})\\
&=& b(N(n\dot{\rho_t}\omega_t^{n-1}))
\end{array}
$$
where we use that $N$  is a trace,  it vanishes on commutators and that
$b\,N= N\, b'$ on cochains. If we apply this formula to the family
$\rho_t = t\rho$, then $\omega_t= tb'\rho-t^2\rho^2$. Integrating from $0$ to
$1$ gives
$$
\begin{array}{rcl}
ch_{2n+2}&=& b \, N \displaystyle{\int _0^1} \rho (tb'\rho -
t^2\rho^2)^{n}dt/n!\\
\rule{0mm}{5mm}
&=& b\, cs_{2n+1}
\end{array}
$$
which proves the lemma. \qed

We note that since  $\rho(1)=1$, $\om$ vanishes if any of its arguments is 1, hence
the cochain $ch_{2n+2}$ is reduced, i.e. it is in the image of $\pi^*$. We now need
to check that the Chern-Simons form descends  to a non-trivial element in
$C_\lambda^{2n}(k)$. For this we calculate
$$\begin{array}{rcl}
f(\underbrace{1, \dots, 1}_{2n+1})&=& cs_{2n+1}(1,\dots,1)\\
& =&
 \displaystyle{\frac{2n+1}{n!}}
\displaystyle{\int_0^1} (t-t^2)^ndt\\
\rule{0mm}{8mm}
&=& \displaystyle{\frac{n!}{(2n)!}}
\end{array}
$$
so  that the Chern-Simons form $cs_{2n+1}$ restricts to
 $C_\lambda^{2n}(k)$. Since the differential $b$ is  zero on
 this
space, this form is also closed.
We note that the resulting cohomology
class is dual to the canonical homology class of the cyclic chain
$((2n)!/n!)\,N \tp{e}{2n+1}$, where $e$ is the identity map on $k$. We have
thus established the following \cite{ja1}.
\btheorem\label{chern}
Let $A$ be a unital algebra and $\rho : A\ra k$ a linear map such that
$\rho(1)=1 $. Let $ch_{2n+2}= N(\omega^{n+1}/(n +1)!)$ be the Chern character
form determined by the map $\rho$ and the trace $N$. Then the reduced cyclic
cohomology class  $[ch_{2n+2}]\in \bar{H}C^{2n+1}(A)$ is induced from the
canonical
generator of $HC^{2n}(k)$ by the connecting map
$HC^{2n}(k)\ra \bar{H}C^{2n+1}(A)$.
\etheorem
\subsection{Reduced cyclic cohomology of the Weyl algebra}
The theorem proved in the last section states that  if $A$ is a unital
algebra then there are canonical reduced cohomology classes in every odd
degree. In this section we want to present an example an algebra whose 
reduced cyclic cohomology can be completely determined  in this way.

Let $n$ be an integer $\geq 1$. Let $V$ be a vector space with basis $\{p_i,
q_j\}$,  $i,j= 1,\dots, n$. We define the {\em Weyl} algebra $A_n$ to be the
quotient of the tensor algebra of the vector space $V$ by the ideal $I$
 generated
by elements
$$
\begin{array}{l}
p_i\otimes q_j - q_j \otimes p_i - \delta_{ij}\\
p_i\otimes p_j - p_j\otimes p_i\\
q_i\otimes q_j - q_j\otimes q_i.
\end{array}
$$

The Weyl algebra $A_n$ is equipped with a canonical filtration
where $F_pA_n$ is the image of $\bigoplus_{n\leq p} V^{\otimes n}$.
As in the case of universal enveloping algebras $gr\, A_n=
S(V)$
is
a (commutative) polynomial ring. We want to  calculate the reduced cyclic
homology of this algebra.
Feigin and Tsygan \cite{ft-LNM} show that the Weyl algebra $A_n$ has the
following cyclic and Hochschild homology.
\bproposition
The Hochschild homology of the Weyl algebra $A_n$ is
$$
H_i(A_n, A_n)=\left\{
\begin{array}{cl}
k,& \quad i=2n\\
0 & \quad \mbox{otherwise}
\end{array}
\right.
$$
The cyclic homology of the algebra $A_n$ is given by
$$
HC_i(A_n)=\left\{\begin{array}{cl}
k,& \mbox{\rm for}\, i\geq 2n\;\;\mbox{\rm and even}\\
0 & \mbox{\rm otherwise.}
\end{array}
\right.
$$
\eproposition

To calculate the reduced cyclic homology of $A_n$, we use the above result
and the long exact sequence relating cyclic homology and reduced
cyclic homology. In particular, we need to show that
 the maps $HC_{2i}(k) \ra HC_{2i}(A_n)$ are nonzero for $i\geq n$.
We check that this is indeed the case using the following theorem of Block
\cite{Block}.
\bproposition
Let $A$ be a filtered algebra and assume that for an integer $n$ the
Hoch\-schild homology $H_i(gr\, A) =0$ for $i>n$. Then the natural map
$ HC_i(F_0 A)\ra HC_i(A)$ is an isomorphism for all $i\geq n$ and a
monomorphism for $i=n-1$. In particular
$HP_i(F_0 A)\simeq HP_i(A)$ for all $i$. 
\eproposition
We now proceed as follows. From the constraint on the Hochschild homology
and Connes' long exact homology sequence follows that the  operator $S$ is
bijective in high degrees. If a map   $HC_{2i}(k) \ra HC_{2i}(A_n)$ were
to vanish, then since these maps are compatible with $S$,
 it would vanish
in degrees $2j$ for all $j\geq i$. Thus the map on periodic homology
would be zero, contradicting the theorem of Block, and we have the following.
\btheorem
The reduced cyclic homology of the Weyl algebra $A_n$ is given by
$$
\bar HC_{i}(A_n)=\left\{ \begin{array}{cl}
k, &{\rm for}\, i=2m+1< 2n\\
0& {\rm otherwise.}
\end{array}
\right.
$$
\etheorem

\bremark The reduced Hochschild homology of the Weyl algebra is the
same as the Hochschild homology of $A_n$ in degrees not less than 2. In
lower degrees we have  an exact sequence
$$
0\lrar H_1(A)\lrar \bar{H}_1(A)\lrar H_0(A)\lrar \bar{H}_0(A)\lrar 0
$$
so that we have $\bar{H}_1(A)= k$.
\eremark

We can easily translate this result into a corresponding statement about
cohomology.
In this case we see from Theorem \ref{chern} that all reduced cyclic cohomology
classes of the Weyl algebra are induced from the cyclic classes of the
field $k$ by means of the transgression procedure explained in the
previous section.
%% this is a short section on the excision in periodic cyclic cohomology 
%% after Cuntz and Quillen
\subsection{Excision in periodic cyclic cohomology}
In the previous sections we have already discussed the results of 
Wodzicki on excision in cyclic homology and Suslin and Wodzicki 
on excision in the algebraic $K$-theory. We have seen 
that we can associate a long exact sequence in cyclic homology 
or rational $K$-theory with the 
short exact sequence of algebras 
$$
0\arr{} I \arr{} R \arr{} R/I\arr{} 0
$$
when the ideal $I$ is $H$-unital. It turns out that, 
in the setting of periodic cyclic cohomology,  a more general 
statement can be proved, and that was done recently by Cuntz and Quillen
\cite{CQExc}. We shall now review their proof. 
First we need to introduce a weaker condition than $H$-unitality. 
\bdefinition
Let $J$ be an algebra over a field of characteristic zero and let 
$\dots \subset  J^n \subset J^{n-1} \subset \dots \subset J^2 \subset J$ be
the filtration of $J$ by its increasing powers. If we denote by 
$C^n(J)$ the space of $n$-cochains on $J$ then we shall say that 
$J$ is {\em approximately $H$-unital} if and only if the complex 
$(\displaystyle{\lim_{\stackrel{\longrightarrow}{\scriptscriptstyle{k}}}} C^\star(J^k), b')$ is 
acyclic. In the following we shall denote by $C^\star(J^\infty) = \displaystyle{
\lim_{\stackrel{\longrightarrow}{\scriptscriptstyle{k}}} }\, C^\star(J^k)$. 
\edefinition
For comparison, let us recall that Wodzicki's notion of $H$-unitality only 
required that the complex $(C^\star(J), b')$ be acyclic. His condition implies 
in particular that $J^2 = J $. 
The following lemma gives a useful criterion for approximate $H$-unitality. 
\blemma\label{appr-H-unit}
Let $J$ be an algebra and let us denote by $m : J\otimes J \ra J $ the 
multiplication map on $J$. Assume that there is a $k\geq 1$ 
and a linear map $\phi: J^k \ra J\otimes J$ which is a section of $m$ in the sense
that $m(\phi(x)) = x $ for $x\in J^k$. Assume moreover that $\phi$ satisfies the following 
formula for $x\in J $ and $y \in J^k$
$$
\phi(xy) = x\phi(y)
$$
Then $J$ is approximately $H$-unital. 
\elemma
Proof. We note that, for every $q$ there is a $p \geq q$ and a $J$-linear splitting
$\psi: J^p \ra J^q \otimes J^q$ of the multiplication $m : J^q \otimes J^q \ra J^{2q}
\supset J^p$ which is obtained by iterating $\phi $ and then multiplying the last $q$ 
variables together. Now define $\Psi: C_n(J^p) \ra C_{n+1}(J^q)$ to be the map 
defined by $\Psi(a_0, \dots , a_n) = (a_0, \dots , \psi(a_n))$. Let $\omega \in C_n(J^p)$. 
We have $(-1)^{n+1} (\Psi b' - b'\Psi)\omega = \omega \in C_n(J^p)$. 
This shows that $J$ is approximately $H$-unital.  \qed

We have the following important example of an application of  this lemma. 
\bproposition\label{semi-free}
\begin{enumerate}
\item Let $J$ be a left ideal in a unital quasi-free algebra $R$ ({\em see Section \ref{RA}}). 
Then there 
exists a left $J$-linear lift $\phi: J^2 \ra J\otimes J$ for the multiplication 
map $m: J\otimes J \ra J^2$ and therefore $J$ is approximately $H$-unital. 
\item Let $e: R\ra J $ be a linear projection onto $J$ such that $e(1) = 0 $. 
Then the formula 
$$
\phi(a) = 1 \otimes a + b'(1\otimes e) b'\alpha (a)
$$
defines an $R$-linear lift $J\ra R\otimes J$ for the multiplication map 
$R\otimes J \ra J$. The restriction of $\phi$ to $J^2$ is therefore a map as 
in the first part. 
\end{enumerate}
\eproposition 

Proof. 1) From \cite[Propn. 5.1]{CQ1} follows that $J$ is a projective left $R$-module
and hence there exists an $R$-linear lift $\phi: J \ra R\otimes J $ for the 
multiplication map $m: R\otimes J \ra J $. 
To obtain the required map, consider the restriction of $\phi $ to $J^2$. 

2) We have that $\phi(a) \in R\otimes J$ if $a\in J $. Since $mb' = 0 $ we 
have that $m\phi(a) = a $. Furthermore, if $a\in R$ and $b\in J$
we have
$$
\begin{array}{rcl}
\phi(ab) & = & 1 \otimes ab + b' ( 1\otimes e) b'(a\alpha(b) + \alpha(a) b + da\;db) \\
& = & 1\otimes ab + b'(1\otimes e)(ab'\alpha(b) + b'(\alpha(a))b - d(a)b\otimes 1
+ da\otimes b ) \\
& = & 1\otimes ab + ab'(1\otimes e) b'\alpha (b) + b'(da\otimes b) 
\\
& = & a\phi(b)
\end{array}
$$
\qed

 Let us now list some examples of 
approximately $H$-unital algebras \cite{CQExc}. It is interesting that these
algebras are not $H$-unital. 
\bexample
Let $ A = \bc[x_1, \dots x_n]$ be the algebra of polynomials in $n$ variables
and let $J= \langle h \rangle$ be the principal ideal $J= hA$ generated 
by a polynomial $h\in A$. Then $J$ is approximately $H$-unital, as we can see by 
defining a map $\phi: J \ra J\otimes J $ by $\phi(gh) = gh\otimes h $ and checking 
that it satisfies the conditions of the lemma above. 
Similarly, principal ideals in the algebra of Laurent polynomials $\bc[z , z^{-1}]$. 
\eexample
\bexample
If $A$ is a Banach algebra and $C^{(k)}(S^1, A)$ is the algebra of $k$-times continuously 
differentiable functions from $S^1$ to $A$, where $k \geq 1$, then let 
$$
S^{(k)}(A) = \{ f \in C^{(k)}(S^1, A) \mid f (1) = 0 \} 
$$
be the corresponding differentiable suspension. Then this algebra is approximately 
$H$-unital.
\eexample
\bexample
An extension of Lemma \ref{appr-H-unit} allows one to show that the Schatten ideals 
$\call^p(\bh)$ in the algebra of bounded linear operators in a separable Hilbert space are 
approximately $H$-unital 
\eexample
\bexample
For an algebra $A$ we denote its  tensor algebra by $TA$ and by $IA$ the 
kernel of the natural surjection homomorphism $TA  \ra A$. From Section \ref{RA} 
(see also \cite{CQ1}) we know that $TA$ may be identified with the nonunital algebra
 $\Omega ^{ev}A$
of even forms over the algebra $A$ equipped with the Fedosov product. 
If we use this identification, then the ideal $IA$ is spanned by forms of degree 
at least two. There is the following $IA$-linear splitting $\phi: (IA)^2 \ra IA \otimes IA$
of the multiplication map: $\phi(\omega da\; db) = \omega \otimes da\; db$. 
This shows that $IA$ is approximately $H$-unital. 
\eexample
\bexample
If $A$ is an algebra, then let $A * A $ be the free nonunital product of $A$ by itself, 
and let $qA$ be the kernel of the folding map $A*A \ra A$. 
The free product $A*A$ is isomorphic to the nonunital algebra $\Omega A$ of 
differential forms over $A$ equipped with the Fedosov product, and the ideal $qA$ 
corresponds to the forms of degree at least one. Then $qA$ is approximately 
$H$-unital as may be seen by defining the splitting $\phi: (qA)^2 \ra qA \otimes qA$
 using the formula
$\phi (\omega da) = \omega \otimes da$. 
\eexample

\btheorem \label{excision-one}
Let $0 \ra J \ra A \ra A/J \ra 0 $ be an exact sequence of algebras and 
suppose that $J$ is approximately $H$-unital. Then there exists a six-term 
exact sequence in periodic cyclic cohomology 
$$
\begin{array}{ccccc}
HP^0(J ) & \larr{} & HP^0(A) & \larr{} & HP^0(A/J) \\
\dal{} &&&&\ual{}\\
HP^1(A/J) & \arr{} & HP^1(A)& \arr{} & HP^1(J)
\end{array}
$$
\etheorem
Proof. (Cf. \cite{CQExc}). The proof relies on the excision properties of 
cyclic cohomology and on the fact that $HP^i = \displaystyle{\lim_{\arr{}}} HC^{2n + i } (A)$
where the direct limit is taken with respect to Connes's periodicity operator $S$. 
For any $k$, let us define the complex of {\em relative cyclic cochains } 
$C^n_\lambda(A,J^k)$ as the cokernel of the canonical map
$C^n_\lambda(A/J^k ) \ra C^n_\lambda (A)$ induced by the 
surjection $ A \ra A/J^k$. Thus for all $k$ we have the following 
exact sequence
\beqn\label{short-k}
0 \arr{} C^n_\lambda(A/J^k) \arr{} C^n_\lambda(A) \arr{} C^n_\lambda(A,J^k) \arr{} 0 
\eeqn
Furthermore, since for each $k$ we have the following short exact sequence
\beqn\label{nilpotent}
0 \ra J^k / J^{k+1} \arr{} A/ J^{k+1} \arr{} A/ J^k \ra 0 
\eeqn
it follows that there is a sequence of surjections
$$
A/J  \larr{} A/J^2 \larr{} A/J^3 \larr{} \cdots \larr{} A/J^k \larr{} \cdots 
$$
Thus there is a direct system of spaces of cyclic $n$-cochains
$$
C^n_\lambda(A/J) \ra C^n (A/J^2) \ra \cdots \; \ra C^n(A/J^k) \ra \cdots 
$$
which, together with the sequence (\ref{short-k}), gives rise to  the following short exact 
sequence
$$
0 \arr{} \lim_{\arr{}} C^n_\lambda(A/J^k) \arr{} C^n_\lambda(A) \arr{} \lim_{\arr{}}
C^n_\lambda(A, J^k)
\arr{} 0 
$$
where the limits are  taken as $k\to \infty$. It now  follows that there is the
following long exact sequence of cyclic cohomology groups
\beqn\label{long-cyclic}
 \ra HC^{n-1}(A,J^\infty) \ra\lim_{\arr{}} HC^n(A,J^k) \ra HC^n(A) \ra
HC^n(A,J^\infty)\ra
\eeqn
Here by $HC^n(A,J^\infty)$ we denote the cyclic cohomology of the complex of 
asymptotic cochains $C^\star_\lambda(A,J^\infty)$. An important step in 
the proof is provided by the following lemma. 
\blemma
The induced map $HC^n(A, J^\infty)\ra HC^n(J^\infty)$ is an isomorphism. 
\elemma

The sequence (\ref{long-cyclic}) is compatible with the operator 
$S: HC^n \ra HC^{n+2}$, 
since $S$ commutes with the restriction map for cochains. Hence we may 
pass to periodic cyclic cohomology by taking the direct limit with respect to the 
action of $S$. Moreover, since for any $k$, the sequence 
$$
0 \arr{} J/J^k \arr{} A/J^k \arr{} A/J \arr{} 0 
$$
represents a nilpotent extension, we have the following fact  
(another theorem of Goodwillie \cite{Good1}): 
The map $HP^\star(A/J) \ra HP^\star (A/J^k)$ is an isomorphism. 
Putting these results together we obtain the following exact sequence of length six. 
$$
\begin{array}{ccccc}
HP^0(J ^\infty) & \larr{} & HP^0(A) & \larr{} & HP^0(A/J) \\
\dal{} &&&&\ual{}\\
HP^1(A/J) & \arr{} & HP^1(A)& \arr{} & HP^1(J^\infty)
\end{array}
$$
If we apply this result to the short exact sequence $0 \ra J \ra J \ra 0 \ra 0 $, we
see that $HP^\star(J^\infty)\simeq HP^\star(J)$.  This finishes the proof of Theorem. \qed

We can now state and prove the general case of the above theorem. 
\btheorem\label{excision-two}
To any short exact sequence of algebras $0 \ra J \ra A \ra A/J\ra 0 $ 
there corresponds an exact sequence in periodic cyclic cohomology as
in Theorem \ref{excision-one}.
\etheorem
Proof. There are the following canonical short exact sequences. First there 
is the universal extension of $A$
\beqn \label{help-one}
0 \arr{} IA \arr{} TA \arr{} A \arr{} 0 
\eeqn
where $TA$ is the nonunital tensor algebra over $A$. Secondly, let $K$ 
be the kernel of the canonical surjection in the following short exact sequence
\beqn\label{help-two} 
0 \arr{} K \arr{} TA \arr{} A/J \arr{} 0 
\eeqn 
The algebras $IA$ and $K$ are ideals in the unital tensor algebra $\widetilde{TA}$ which 
is free and therefore quasi-free. Hence they are both approximately 
$H$-unital by Proposition \ref{semi-free}. Then Theorem \ref{excision-one} applied
to sequences (\ref{help-one}) and (\ref{help-two}) gives that 
$HP^\star(IA)\simeq HP^{\star + 1} (A)$ and $HP^\star(K) \simeq HP^{\star + 1} (A/J)$. 
If we now consider the following short exact sequence of algebras
$$
0 \arr{} IA \arr{} K \arr{} J \arr{} 0 
$$
then the two facts mentioned above together with the long exact sequence 
given by Theorem \ref{excision-one} give the proof of the present Theorem. \qed
\bremark
Cuntz and Quillen provide also a considerable extension of this result 
to the case of bivariant periodic cyclic homology \cite{CQbivariant}. Moreover, 
Cuntz has proved that bivariant periodic cyclic homology of Fr\'{e}chet 
algebras also has the excision property \cite{Cutopology}. 
\eremark

\vfill\eject

\section{Entire cyclic cohomology of Banach algebras}
\subsection{Topology on $\Omega A$ and $QA$}
So far, in our discussion of various cyclic theories associated
with an algebra $A$ we have ignored the question of topology. 
Even if $A$ was a topological algebra, we have disregarded this 
extra piece of information. We now want to discuss the entire cyclic 
cohomology, which is 
a generalization of the periodic cyclic cohomology to the
class of Banach algebras. This 
version of  cyclic cohomology, which  was introduced by Connes in \cite{Connes:Entire}, 
seems better suited for the study of infinite dimensional spaces, 
like the configuration space of the 
constructive field theory or the dual space of higher rank discrete groups. 

We define the entire cochain complex as follows.  A cochain in this complex
is a linear map $f: T(A)\ra \bc$ which is represented by a  possibly infinite 
sequence of multilinear functionals on $A$
$$
(f_0, f_1, \dots , f_{n}, \dots )
$$
where $f_n:A^{\otimes n+1}\ra \bc$. A cochain is called even (odd) if 
its components $f_n$ are nonzero when $n$ is even (odd). This way 
we obtain a $\bz/2$-graded complex equipped with the differential $b+B$
acting between spaces of cochains of opposite parity. This complex 
is not very interesting as its cohomology is trivial. When $A$ is a Banach 
algebra, inside this cochain complex there is a subcomplex of entire 
cochains satisfying the following growth condition
\beqn \label{growth}
\sum_{n\geq 0} \norm{f_{2n}}\, n! \, r^n < \infty
\eeqn
for all positive real $r$. There is an identical condition to be satisfied 
by odd cochains. We shall denote this complex
of unnormalized entire cochains by $C_\epsilon^*(A)$. 
The cohomology of this complex is the {\em entire cyclic 
cohomology}.  We can define a similar theory when we consider 
the same growth condition but we require that the cochains be simplicially 
normalized. In other words, if $f$ is a cochain, we require that
$$ 
f_n(a_0, a_1, \dots, a_n )=0 
$$
whenever $a_i= 1$ for some $1\leq i\leq n$. In other words, 
a simplicially normalized cochain is a linear map $f:\Omega A \ra \bc$. 
The space of simplicially normalized entire cochains will 
be denoted by $\Omega^*_\epsilon A$.

We begin our excursion into the entire cyclic cohomology with a discussion 
of the topological versions of universal algebras defined in Section \ref{FORMS}.  
Let then   $A$ be a unital Banach algebra equipped with a norm
$\emnorm$. Let us denote by $\bar{A} = A/\bc$ the quotient 
space which is a Banach space equipped with the quotient space 
norm. The space of one-forms $\Omega^1A$ becomes a Banach 
space with the completed projective tensor product norm
$\Omega^1 A = A\hat{\otimes}_\pi \bar{A}$. 
 The operator $d: A\ra \Omega^1A$ is now a continuous 
operator of unit norm.  For any $n\geq 1$ we equip the 
space of $n$-forms 
\beqn\label{product-on-forms}
\Omega^n A= \Omega^1 A\otimes_A\stackrel{n}{\cdots} \otimes _A\Omega^1 A
\eeqn
with the topology of the completed projective tensor product 
\cite[p. 260 ff]{Arveson}.
The differential $d$ becomes 
a continuous derivation of unit norm and the product in $\Omega A$
 is now a continuous map.

Given this topological structure, we may now equip
the differential graded algebra $\Omega A = \bigoplus _{n\geq 0}
\Omega ^n A$ with the following family of norms $\emnorm_r$
parametrized by a positive real $r$. If 
$\omega = \sum_{n\geq 0} \omega_n$, we put 
$$
\norm{\omega}_r = \sum_{n\geq 0} \norm{\omega_n}
r^n
$$
The derivation $d$ becomes a continuous map $d: (\Omega A, \emnorm_r)
\ra (\Omega A, \emnorm_r)$. In particular,  $\norm{da}_r \leq 
r \norm{a}$ for $a\in A$. If $\omega_n$ is a homogeneous form 
of degree $n$ then 
$$
\norm{\omega_n}_r = r^n \norm{\omega_n}, \qquad 
\norm{d\omega_n}_r = r^{n+1} \norm{d\omega_n}\leq r^{n+1}
\norm{\omega_n} = r\norm{\omega}_r
$$

For any two $r' < r$, let 
$$
f_{rr'} : (\Omega A, \emnorm_r) \longrightarrow  (\Omega A, \emnorm_{r'}) 
$$ 
be  the identity map, which in this case is norm decreasing. 
This way we obtain a direct system $\{(\Omega A, \emnorm_r), f_{rr'}\}$.
Let $\Omega _r A$ denote the completion of the space 
$(\Omega A \emnorm_r) $
with respect to the norm $\emnorm _r$. 
This is now a Banach algebra, which means that for any two forms $\omega$, 
$\eta $ in $\Omega_r A$ we have \cite[p. 262]{Arveson}
$$
\norm{\omega\eta} _r\leq \norm{\omega}_r\norm{\eta}_r.
$$
where the product is defined using the formula (\ref{product-on-forms}). 

There is an induced 
direct system $(\Omega_r A, f_{rr'})$ of Banach algebras. 
We denote by $\Omega_\epsilon A$ the direct limit of this system
as $r\to 0$: 
$$
\Omega_\epsilon A= \lim_{\lrar} \Omega_r A
$$
$\Omega_\epsilon A$ is a locally convex algebra. 
Finally, we  denote by $f_r$ the canonical maps  
$f_r : \Omega _r A\ra \Omega_\epsilon A$, which are continuous algebra 
homomorphisms.  The derivation $d$  extends  to a continuous derivation 
$d_r: \Omega_r A\lrar \Omega_r A$. The composite 
$f_r\circ d_r: \Omega_r A\lrar \Omega_\epsilon A$ is  a continuous
derivation of $\Omega_r A$ with values in $\Omega_\epsilon A$ and all 
derivations $f_r\circ d_r $ are compatible with the maps $f_{r'r}$. 
Thus this family of derivations is equivalent to a continuous derivation 
$d: \Omega_\epsilon A\lrar \Omega_\epsilon A$ which makes 
$\Omega_\epsilon  A$ into a topological differential graded algebra  
equipped with the locally convex topology induced by the topologies 
on $\Omega_r A$.

In the case when $A$ is a unital Banach algebra, for any positive 
real $r$, the completion $ \Omega_r\tilde{A}$ is 
an augmented Banach algebra
$$
\Omega_r\tilde{A} = \bc \oplus \bar{Q}_r\tilde{A}
$$
from which it follows that also the direct limit of the corresponding
direct system is an augmented algebra 
$$
\Omega _\epsilon \tilde{A} = \bc \oplus \bar{\Omega}_\epsilon 
\tilde{A}
$$

We now define a topology on the Cuntz algebra $QA$ using the families 
of norms $\rnorm{\;}$ on $\Omega A$ together with  the isomorphism 
$\mu$ (\ref{mu-map}) between the two superalgebras. If $\omega$ is an element 
of the Cuntz algebra $QA$ then its norm $\norm{\omega}_{Q_r}$, for any 
positive $r$ is defined by 
$$
\norm{\omega}_{Q_r} = \norm{\mu(\omega)}_r
$$
In particular, for $a\in A$,  $\norm{a}_{Q_r} = \norm{a+ da}_r = 
\norm{a} + r\norm{da}$. This means that the isomorphism $\mu$
induces a new norm on the  the algebra $A$ which is equivalent to the 
original Banach norm on $A$. 

We shall  need the following estimates for the 
Fedosov product. 
\blemma
For any $r> 0$,  and for any $\omega, \eta \in \Omega_rA$   we have
$$
\norm{\omega \circ \eta}_r \leq C_r \norm{\omega}_r\norm{\eta}_r
$$
where $C_r = 1 + r^2$ and  $\circ $ denotes the Fedosov product of forms. If either of the 
forms $\omega$, $\eta$ is closed then this estimate reduces to 
$$
\norm{\omega\circ\eta}_r \leq \norm{\omega}_r \norm{\eta}_r
$$
\elemma
Proof.  If $\omega = \sum \omega_n $ and $\eta = \sum \eta_n$ then 
$$
\begin{array}{rcl}
\omega \circ \eta & = & \displaystyle{\sum_{n\geq 0}\sum_{i+j = n}}
\omega_i \circ \eta_j \\
\rule{0mm}{8mm}
& = & \displaystyle{\sum_{n\geq 0}\sum_{i+j = n}} 
\omega_i \eta_j - (-1)^{|\omega_i|}d\omega_id\eta_j
\end{array}
$$
The norm $\emnorm_{Q_r}$ of this form may then be estimated as 
follows. 
$$
\begin{array}{rcl}
\norm{\omega\circ \eta}_{r} &\leq  & \displaystyle{\sum_{n\geq 0}
\left(\sum_{i+j= n} \norm{\omega_i \eta_j }\right) }r^n 
+ \displaystyle{\sum_{n\geq 0}
\left(\sum_{i+j= n} \norm{d\omega_i d\eta_j} \right)} r^{n+ 2}\\
\rule{0mm}{10mm} 
& \leq  & \displaystyle{\sum_{n\geq 0}
\left(\sum_{i+j= n} \norm{\omega_i }\norm{\eta_j }_\pi\right) }r^n 
+ r^2\displaystyle{\sum_{n\geq 0}
\left(\sum_{i+j= n} \norm{\omega_i}\norm{ \eta_j }
\right)} r^n 
\\
\rule{0mm}{10mm}
& = & (1+ r^2) \norm{\omega}_r\norm{\eta}_r
\end{array}
$$
The last statement of the lemma is now clear. 
\qed

\bcorollary\label{Cuntz-estimate}
If $x,y\in Q_rA$ then 
$$
\norm{xy}_{Q_r} = \norm{\mu(x)\circ \mu(y)}_r 
\leq (1+r^2) \norm{\mu(x)}_r\norm{\mu(y)}_r = 
C_r\norm{x}_{Q_r}\norm{y}_{Q_r}
$$
 \ecorollary
\qed

Hence,  in particular, the multiplication in the Cuntz algebra 
is continuous in this topology.  

We now introduce a direct system of normed spaces $(QA, \qrnorm{\;})$
as  in the case of the algebra of differential forms. For any 
pair $r' < r$ there is a map 
$$
f_{rr'} : (QA, \emnorm_{Q_r}) \lrar (QA, \emnorm_{Q_{r'}})
$$ 
which is the identity map. Using the estimates 
$$
\norm{\omega}_{Q_{r'}} = \norm{\mu(\omega)}_{r'}
< \norm{\mu(\omega)}_r = \norm{\omega}_{Q_r}
$$
we see that for any two $r' < r$ the map $f_{rr'}$ is continuous 
and norm decreasing.  Let us  denote by $Q_rA$ the completion
of the space $(QA, \qrnorm{\;})$,  then there is a direct 
system $( Q_rA, f_{rr'})$ of complete normed spaces. We denote 
by $Q_\epsilon A$ its direct limit taken as $r\to 0$. 

As before, in the case of a unital  Banach algebra $A$,
the Cuntz algebra  $Q_\epsilon \tilde{A}$ is also an augmented 
algebra. This follows from the fact that for any $r$ we
have
$$
Q_r \tilde{A} = \bc \oplus \bar{Q}_r\tilde{A}
$$
so that we have 
$$
Q_\epsilon \tilde{A} = \bc \oplus \bar{Q}_\epsilon\tilde{A}
$$
where $\bar{Q}_\epsilon \tilde{A}$ is the direct limit of 
the direct subsystem $\bar{Q}_r \tilde{A}$. 

\subsection{Simplicial normalization}

We now want to prove the simplicial normalization theorem, 
which states that the entire complexes defined with the 
use of either unnormalized entire cochains $C^*_\epsilon(A)$, 
or normalized cochains $\Omega^*_\epsilon A$, lead to the same
entire cohomology theory. 
\btheorem\label{simplicial}
The inclusion 
$$
\Omega_\epsilon A^* \ra C_\epsilon (A)^*
$$
induces an isomorphism 
$$
HE^*_{sn} (A) \simeq HE^*(A)
$$
\etheorem

\subsubsection{A homotopy of superalgebras}
A first step in the proof of the simplicial normalization theorem for entire cyclic 
cohomology will be to show that the algebras $Q_\epsilon \tilde{A}$ and 
$ \bc \oplus Q_\epsilon A$ are homotopy equivalent. 
At the same time we shall 
establish a homotopy equivalence between the superalgebras 
$\bar{Q}_\epsilon
\tilde{A}$ and $Q_\epsilon A$. The required homotopy is defined using 
a certain family $\{\phi_t\}$ of continuous superalgebra homomorphisms
$\phi_t : Q_\epsilon\tilde{A} \ra Q_\epsilon \tilde{A}$. The definition and properties 
of the family $\phi_t$ (see also
\cite{Cu2}) will depend in an important way on the properties of the 
operators $L$ and $W_t$ established in  section \ref{Cuntz-algebra}. The 
exposition below follows \cite{ja2}.

For any $t\in [0,1]$ we put 
\beqn\label{family}
\begin{array}{rcl}
a & \mapsto & a_t = e^{-tL/2} a e^{tL/2}  = W_{-t} a W_t\\
a^\gamma & \mapsto & a^\gamma _t = e^{tL/2} a^\gamma e^{-tL/2}
= W_t a^\gamma W_{-t}
\end{array}
\eeqn
We note that the first of the two maps is a homomorphism 
from the algebra $A$ to the underlying algebra of $Q_\epsilon \tilde{A}$
and that the other map is induced from the first by applying the 
canonical automorphism $\gamma$ which sends $a_t$ to $a_t^\gamma$.
This pair of algebra homomorphisms extends to the following family of  homomorphisms
$\phi_t: Q \tilde{A} \lrar Q_\epsilon \tilde{A}$. 
For $\omega \in Q_\epsilon\tilde{A}$, 
$\omega = p(a_0) q(a_1)\cdots q(a_n)$, we have
$$
\phi_t(\omega) = \omega_t = p_t(a_0) q_t(a_1) \cdots q_t(a_n)
$$
where 
$$
\begin{array}{rcccl}
p_t(a) & = & p(a_t) & = & 1/2(a_t + a^\gamma_t)\\
q_t(a) & = & q(a_t) & = & 1/2(a_t - a^\gamma_t)
\end{array}
$$
We remark that $\phi_0$ is the identity map whereas $\phi_1$ is the 
map which sends $f$ to $f^\gamma$ by (\ref{conjugation}).

\btheorem\label{main}
The family $\phi_t$ defined above is a family of continuous superalgebra 
homomorphisms from $Q_\epsilon\tilde{A}$ to $Q_\epsilon\tilde{A}$. 
\etheorem

It is clear that the folding map in $Q\tilde{\bc}$ 
\beqn\label{fold}
0 \ra q\tilde{\bc} \ra Q\tilde{\bc} \ra \tilde{\bc}\ra 0
\eeqn
induces a canonical projection
$$
Q\tilde{A} \lrar \bc \oplus QA
$$
The kernel of this projection is an ideal  $I$ generated in $Q\tilde{A}$  
by $q\tilde{\bc}$, the kernel of the folding map (\ref{fold}). 
In the analytic situation, in any of the completions
$Q_r \tilde{A}$ we consider the closed ideal $I_r$, which is the 
completion of the ideal $I$ 
in the topology defined by the norm $\qrnorm{\mbox{}\;\;}$. Thus we 
deduce that there 
are the following  continuous isomorphisms
\beqn\label{rism}
\bc \oplus Q_r A  = Q_r\tilde{A}/ I_r, \qquad  Q_r A = \bar{Q}_r \tilde{A}/
I_r
\eeqn
This observation, together with the properties of the 
family $\{\phi_t\}$ of homomorphisms
of $Q_\epsilon \tilde{A}$ allows us to define a lift of 
$Q_\epsilon A$ back into 
$Q_\epsilon \tilde{A}$. In fact we prove
\bproposition\label{phi1/2}
We have the following continuous isomorphism of locally convex 
superalgebras 
$$
\phi_{1/2} (\bar{Q}_\epsilon\tilde{A}) \simeq 
Q_\epsilon A
$$
\eproposition
Proof. First we define a map $\phi_{1/2}(\bar{Q}_\epsilon \tilde{A}) 
\ra Q_\epsilon A$. 
Since $\phi_t: \bar{Q}_\epsilon \tilde{A}\ra \bar{Q}_\epsilon\tilde{A}$
is a continuous homomorphism for any $t\in [0,1]$, there exists 
a positive real $r$ such that  $\phi_{1/2}(\bar{Q}\tilde{A}) \subset
\bar{Q}_r\tilde{A} $. We then put for any $\omega\in \bar{Q}_\epsilon 
\tilde{A}$ 
$$
\phi_{1/2}(\omega) \mapsto [\phi_{1/2} (\omega)] \in \bar{Q}_r \tilde{A}/I_r 
$$
Since this map is just the restriction of the canonical projection 
$\bar{Q}_r \tilde{A} \ra \bar{Q}_r\tilde{A}/I_r$ to the image of $\phi_{1/2}$, 
it is continuous. 

To define the inverse map, let $[\omega] \in \bar{Q}_r 
\tilde{A}/I_r$ and let us 
choose a representative $\omega_0 \in Q_r A$ of this class. We define
$$
\xi : \bar{Q}_r\tilde{A}/I_r \ni [\omega] \mapsto 
\phi_{1/2}(\omega_0) \in Q_{r'} \tilde{A}
$$
This map is well defined, since for any $r>0$, $\phi_{1/2}$ vanishes 
on the ideal $I_r$, so that the above assignment
does not depend on the choice of $\omega_0$. To show that it is 
continuous we 
note that, as the space $\bar{Q}_r\tilde{A}/I_r$ is equipped with the 
quotient space  norm 
$$
\qrnorm{[\omega]} = \inf_{\eta\in I_r} \qrnorm{\omega + \eta},
$$
for any $r$ there exists $\omega_0 \in [\omega]$ such that 
$\qrnorm{\omega_0}\leq 2\qrnorm{[\omega]}$. 
Since $\phi_t$ is continuous for any $t\in [0,1]$, from the estimate
 follows that there is an $r' <r$ such that 
$$
\norm{\xi([\omega])}_{Q_{r'}} =  \norm{\phi_{1/2}(\omega_0)}_{Q_{r'}}
\leq M\qrnorm{\omega_0} \leq 2M\qrnorm{[\omega]}
$$
This shows that the map $\xi$ is continuous for all $r$ and
so it is continuous  in the locally convex topology on the union of 
$\bar{Q}_r\tilde{A}/I_r \simeq Q_r A $.  This finishes the proof. \qed

With this proposition we have constructed the following homotopy
\beqn\label{homotopy}
Q_\epsilon\tilde{A} \lrar  \bc \oplus \phi_{1/2}(\bar{Q}_\epsilon\tilde{A}) 
\simeq \bc\oplus Q_\epsilon A \stackrel{l}{\hookrightarrow}
 Q_\epsilon \tilde{A}
\eeqn
where the  map $l$  is the lift provided by the previous theorem. 
There is a corresponding homotopy of the reduced part $\bar{Q}_\epsilon
\tilde{A}$ onto $Q_\epsilon A$
$$ 
\bar{Q}_\epsilon\tilde{A} \lrar \phi_{1/2}(\bar{Q}_\epsilon\tilde{A})
\simeq Q_\epsilon A \stackrel{l}{\hookrightarrow}
 \bar{Q}_\epsilon \tilde{A} 
$$

\subsubsection{Homotopy of supertraces}
We now proceed with  the discussion of 
the action of the family $\phi_t$ on the entire cyclic cohomology. We shall do this 
extending  a relation between cocycles in the $\bz_2$-graded complex 
and supertraces which was discussed
in  Section \ref{Supertraces} to the entire case. We shall also need to show that 
supertraces that are homotopic in a sense explained below correspond to 
cohomologous cocycles.  

To describe 
homotopy of supertraces on $QA$ we use the $X$-complex approach of Cuntz and 
Quillen \cite{CQ3} which we already encountered in the section \ref{RA}. We therefore need to discuss supertraces
on the superbimodule $\Omega^1QA$ over the Cuntz algebra $QA$. 
As before, supertraces on $\Omega ^1QA$ are linear functionals $T$ which 
vanish on the supercommutator quotient space $\Omega^1 QA_{\natural s}
= \Omega^1QA/ [QA, \Omega^1QA]_s$. 
First of all, it turns out that we may identify the 
supercommutator quotient space 
$\Omega^1Q_{\natural s}A$ with the space 
$\Omega ^1 A \oplus \Omega A$ \cite{CQ4}. 
We shall use the following notation. Let $\theta $ be the 
inclusion $\theta : A \ra QA$
and $\theta^\gamma$ the composition of $\theta $ with the canonical 
automorphism $\gamma$ of $Q A$. We let $p$ denote the cochain on $A$ given 
by $p  = (\theta +  \theta^\gamma) /2$ and similarly 
$q = (\theta - \theta^\gamma)/2$. 
We then have the following identification of the even 
part $(\Omega^1 Q_{\natural s} A)^{ev}$
of the supercommutator space
$$
\bigoplus _{n\geq 0 } A \otimes \bar{A}^{2n+1} \arr{\simeq} 
(\Omega^1 Q_{\natural s} A)^{ev}
$$
This 
isomorphism is given by the cochains $(\theta^\gamma q ^{2n} d\theta)^{ev}$. 
The odd part of the supercommutator space is 
$$
\Omega^1A \oplus \bigoplus _{n\geq 1} A \otimes 
\bar{A}^{\otimes 2n} \arr{\simeq}
(\Omega ^1QA_{\natural s})^{od}
$$
where this time the isomorphism is given by cochains 
$(\theta d\theta)^{od}$, $(\theta 
q^{2n-1}d\theta)^{od}$. 
Using this coordinate system, Cuntz and Quillen 
show \cite{CQ4} that  the components of an even 
supertrace $T$  on $\Omega^1 Q A$ are given by the following cochains
$$
T_{2n+1} = T(\theta^\gamma q^{2n} d\theta): A\otimes \bar{A}^{2n} \ra \bc
$$
On the other hand, an odd supertrace $T$ has the following cochains as 
its components
$$
T_{2n} = T(\theta q^{2n-1}d\theta ), \qquad T_1 = 
T\theta d\theta, \qquad T_1 b = 0 
$$
To summarize, using the above isomorphisms we can 
represent any supertrace 
$T$ on $\Omega^1 Q$ in terms of cochains on $\Omega^1 A \oplus \Omega A$. 
 From now on components of  functionals on $\Omega^1 QA$ 
will always be written in terms of these isomorphisms. We can 
now state the following \cite{CQ4}.
\bproposition\label{supertraces}
 Let $T$ be a supertrace on $\Omega ^1QA$ and 
let $\delta $ be the differential $\delta : QA \ra \Omega^1 QA$. 
Then we have 
\beqn\label{homo-ident}
\begin{array}{lcc}
(T\circ \delta) _{2n}& = &   - nPbT_{2n-1} + BT_{2n+1}\\
\rule{0mm}{5mm}
(T\circ\delta)_{2n+1} & = &  (n +1/2) PbT_{2n} - BT_{2n+2}
\end{array}
\eeqn
\eproposition

We say that two supertraces
on $QA$ are homotopic if and only if there exists a supertrace $T$ on 
$\Omega^1 QA$  such that 
$$
\tau_0 - \tau _1 = T\circ \delta .
$$
\bproposition\label{homot-straces}
If $\tau_0$ and $\tau_1$ are homotopic supertraces on $QA$ then the 
corresponding
cocycles $\tau_0^z$ and $\tau^z_1$ are cohomologous in the cochain
complex $(\Omega A^*, b+B)$. 
\eproposition
Proof.  Let us first assume that  $\tau_0$ 
and $\tau_1$  are even homotopic supertraces on $QA$. Since the 
difference of two 
even supertraces is an even supertrace, we can rescale it and 
use the formula (\ref{homo-ident}) to find the corresponding cocycle. 
We get
$$
\begin{array}{rcl}
(\tau_1 - \tau_0)_{2n}^z &  = & \displaystyle{ \frac{(-1)^n}{n!}
(\tau_1 - \tau_0)}\\
\rule{0mm}{8mm} & = &\displaystyle{ \frac{(-1)^n}{n!}(T\circ \delta )_{2n}}\\
\rule{0mm}{8mm} & = & \displaystyle{ \frac{(-1)^n}{n!}
( - nPbT_{2n-1} + BT_{2n+1})}
\end{array}
$$
The rescaled cochain $(\tau_1 - \tau_0)^z $ is a $\kappa$-invariant 
cocycle in the complex $(\Omega A, b+B)$. Moreover, 
it is invariant under the operator $P$, and so the
 above identity indicates that therefore the right hand side 
must also be $P$-invariant. We thus have
$$
(\tau_1 - \tau_0)^z = (\tau_1 - \tau_0)^z P 
= - n bP T_{2n- 1} + BPT_{2n+1}
$$
where we have used that $P$ commutes with $b$ and $B$. If we denote
 by $\tilde{T} _{2n+1}
= \displaystyle{(-1)^n P T_{2n+1}/n!}$ we may write
$$
(\tau_1 - \tau_0)^z_{2n} = b\tilde{T}_{2n- 1} + B \tilde{T}_{2n+1}
$$
which demonstrates that the cocycle $(\tau_1 - \tau_0)^z$ is a coboundary. 

Let us now assume that   $\tau_0 $ and $\tau_1$ are odd supertraces. 
In this case, the rescaled cocycle $(\tau_1 - \tau_0)^z$ satisfies 
the following homotopy identity 
$$
\begin{array}{rcl}
(\tau_1 - \tau_0)^z_{2n+1} &  = &  \displaystyle{\frac{(-1)^n 2^n }{(2n + 1)!!}
 ( (n + \frac{1}{2} ) PbT_{2n} - B T_{2n+2}) }\\
\rule{0mm}{5mm} & = & b\tilde{T}_{2n} + B \tilde{T}_{2n+2}
\end{array}
$$
where $\tilde{T}_{2n} = \displaystyle{\frac{(-1)^n 2^n}{(2n+1) !!}} PT$. 
We have used again the fact that the cocycle 
$(\tau_1 - \tau_0)^z $ is $P$-invariant. \qed

Finally, in the case when $A$ is a Banach algebra, one shows as in 
\cite[Prop. 2.4]{Connes:Entire} that cocycles coming from a continuous supertrace 
on $Q_\epsilon A$ satisfy the entire growth condition \ref{growth}. 
Moreover, any entire cohomology class may be represented by an entire 
cyclic cocycle that comes from a continuous supertrace on $Q_\epsilon A$. 

Our next step is to prove that the family $\phi_t$ of 
homomorphisms establishes an isomorphism between the 
homotopy classes of supertraces on $\bar{Q}_\epsilon\tilde{A}$ and 
$Q_\epsilon A$. 
\bproposition\label{super-homo}
Let $\tau$ be a continuous supertrace on $\bar{Q}_\epsilon\tilde{A}$. 
Then $\tau_t= \tau\circ\phi_t$ is a continuous 
supertrace on $\bar{Q}_\epsilon\tilde{A}$
whose homotopy class does not depend on $t$. 
\eproposition
Proof (Sketch).  Since $\phi_t$ is a superalgebra homomorphism, it is 
clear that $\tau_t = \tau \circ \phi_t$ is a continuous supertrace
for all $t\in [0,1]$. Let us define the following map $T$ on 
$\Omega^1Q_\epsilon \tilde{A}$. For any $\omega$, $\eta $ 
in $\bar{Q}_\epsilon\tilde{A}$ we put 
$$
T(\omega\delta \eta) = \int^1_0 \tau(\phi_t(\omega)\dot{\phi}_t
(\eta))dt
$$
One checks that $T$ is a well-defined continuous supertrace on 
$\Omega^1Q_\epsilon\tilde{A}$. \qed

\bcorollary\label{koniec}
The family $\phi_t$ of homomorphisms of $Q_\epsilon \tilde{A}$
induces an isomorphism of homotopy classes of supertraces 
on $\bar{Q}_\epsilon\tilde{A}$ and $Q_\epsilon A$. 
\ecorollary

Hence, using Connes's correspondence between entire cyclic cohomology 
classes and continuous supertraces, we may use Corollary \ref{koniec}
to deduce that the entire cochain complexes $\Omega_\epsilon A^*$ and 
$\bar{\Omega}_\epsilon A^*$ are quasi-isomorphic. Moreover, Cuntz and 
Quillen have proved that the entire cochain complexes 
$\bar{\Omega}_\epsilon A^*$ and $C_\epsilon^* (A)$ are also quasi-isomorphic. 
We thus have the isomorphisms:
$$
HE^*_{sn}(A)\simeq H^*(\bar{\Omega}_\epsilon A^*, b+B) \simeq HE^*(A)
$$
which finishes the proof of simplicial normalization theorem \ref{simplicial}. 
See \cite{ja2}. 

\subsection{Homotopy invariance of entire cyclic cohomology}
As another application of the theory of supertraces to the entire cyclic cohomology 
we present a different proof of the homotopy invariance of entire cyclic 
cohomology \cite{Masoud}.

Let us then consider two unital  Banach  algebras
$A$ and $B$ and a continuous family of
homomorphisms $f_t: A \ra B$. We assume that each of the homomorphisms $f_t$
is continuous for $t\in [0,1]$ and that  the family is uniformly bounded,
i.e. that there exists a constant $M$ such that $\norm{f_t}\leq M$ for
$t\in [0,1]$. We suppose, moreover, that the corresponding family of
derivatives $\dot{f}_t$ with respect to the parameter $t$ is continuous
and uniformly bounded by a constant $N$ on the interval $[0,1]$. We then
have the following.
\btheorem
Let $f_t : A \ra B$ be a family of homomorphisms of Banach algebras
satisfying the properties given above. Then $f_0$ and $f_1$ induce the same
map on the entire cyclic cohomology $HE^*(B) \ra HE^*(A)$.
\etheorem
Proof. The family  of homomorphisms $f_t: A\ra B$  extends to a family of continuous
homomorphisms  of the corresponding  Cuntz algebras
$\Phi_t:Q_r A\ra Q_r B$ for any positive real $r$ if we define
$$
\Phi_t(p(a_0)q(a_1) \cdots q(a_n)) = p(f_t(a_0))q(f_t(a_1))\cdots q(f_t(a_n))
$$
To prove that this map is continuous for any positive real $r$ we first need to find
estimates for $\qrnorm{p(f_t(a))}$ and $\qrnorm{q(f_t(a))}$. From the definition
of the norm on the Cuntz algebra we have
$$
\norm{p(f_t(a))} = \norm{f_t(a)} \leq M \norm{a} = M \norm{p(a)}
$$
and
$$
\begin{array}{rcl}
\norm{q(f_t(a))} & = & \norm{d(f_t(a))} = \displaystyle{\inf_{\lambda\in \bc}}\;
\norm{f_t(a) + \lambda f_t(1)} \\
\mbox{} \\
&  \leq  & \norm{f_t} \displaystyle{\inf_{\lambda \in \bc}}\; \norm{a+\lambda}
= \norm{f_t} \norm{da} \\
\mbox{} \\
&  \leq  & M \norm{q(a)}
\end{array}
$$
Thus if $\omega= \sum \omega_n \in \Omega_r A$ then using the above estimates
we have that
$$
\begin{array}{rcl}
\qrnorm{\omega} & = & \displaystyle{\sum \norm{\Phi_t(\omega_n)}}
 r^n \leq M^{n+1} \norm{\omega_n}r^n \\
\mbox{} \\
& = &  M\qrnorm{\omega}
\end{array}
$$
This shows that $\Phi_t$ is a continuous homomorphism $\Phi_t: Q_r A \ra Q_rB$
for all positive $r$  and so it extends to  a continuous homomorphism
$\Phi_t : Q_\epsilon A \ra Q_\epsilon B$. This means that $\Phi_t$ induces a
map $\Phi_t^*$ from the space of continuous supertraces on $Q_\epsilon B$
to the space of continuous supertraces on $Q_\epsilon A$.

We now want to prove that $f_0$ and $f_1$ induce the same map on the entire
cyclic cohomology.
Let us then assume that $T$ is a continuous supertrace on
$Q_\epsilon B$, we want to show that the supertraces
on $Q_\epsilon A$ $T\circ \Phi_0$ and $T\circ \Phi_1$ are
homotopic. But we can define a supertrace on $\Omega^1 Q_\epsilon B$
by the formula
$$
\tau(\omega\delta \eta) = \int^1_0 T(\Phi_t(\omega)\dot{\Phi}_t(\eta))dt
$$
It is clear that we have $T\delta \eta = T\Phi_1(\eta) -
T\Phi_0(\eta)$ which means that the supertraces
$T\Phi_0$ and $T\Phi_1$ are homotopic.
We now need to check that $\Phi_t$ and $\dot{\Phi}_t $ are continuous.
But since we have  for any $a\in A$
$$
\begin{array}{rcl}
d/dt(p(f_t(a))) &  =  & p(\dot{f_t}(a)) \\
d/dt(q(f_t(a))) & = &  q(\dot{f_t}(a))
\end{array}
$$
we can estimate the norms of these  elements as before, using the fact that the
derivative $\dot{f_t}$ is uniformly bounded by $N$. We find that
$$
\begin{array}{rcl}
\norm{p(\dot{f}_t(a))} & \leq  &  N\norm{p(a)} \\
\rule{0mm}{8mm}
\norm{q(\dot{f}_t(a))} & \leq & N \norm{q(a)}
\end{array}
$$
From this point on the proof is the same as in the case of simplicial normalization
theorem.
\qed

\subsection{Derivations and entire cyclic cohomology}

As a final application of the technique of supertraces on the Cuntz algebra to entire 
cyclic cohomology we shall derive Khalkhali's result \cite{Masoud} 
on the triviality of the 
action of derivations on the entire complex \cite{ja3}. 

Let us assume that $D: A\ra A$ is a derivation. There is an induced Lie derivative $L_D$
acting on the differential forms $\Omega A$ which is defined using the 
following formula
\beqn\label{Lie-derivative}
L_D(a_0da_1 \dots da_n) = \sum_{i=0}^n a_0da_1\cdots d(Da_i)\cdots da_n
\eeqn

We check that $L_D$ is a well-defined map of $b+B$-complexes. 
\blemma 
The Lie derivative commutes with $b$, $\kappa$ and $d$. Consequently, it commutes 
with $B$. 
\elemma
Proof. The first statement follows from the fact that $L_D$ is a derivation on $\Omega A$. 
Let us consider 
$$
\begin{array}{rcl}
L_Db(\omega da) & = & L_D((-1)^{\deg(\omega)}(\omega da - da \omega))\\
& = & (-1)^{\deg(\omega)}(L_D(\omega)da- L_D(da)\omega + \omega L_Dda - daL_D\omega)\\
& = & b(L_D(\omega)da) + b(\omega(L_Dda))\\
& = & bL_D(\omega da)
\end{array}
$$
Hence $L_Db = bL_D$. We  check the remaining two cases in the same way. 
\qed

Another simple calculation, using the definition of the Fedosov product and 
the fact that $L_D$ commutes with $d$ establishes the following lemma. 
\blemma
$L_D$ is a derivation on  the Cuntz algebra. 
\elemma
Proof. 
$$
\begin{array}{rcl}
L_D(\omega\circ \eta) & = & L_D(\omega\eta) - (-1)^{\deg(\omega)}L_D(d\omega d\eta)\\
& = & L_D(\omega)\eta + \omega L_D(\eta) -(-1)^{\deg(\omega)}(d(L_D\omega)d\eta
- d\omega dL_D\eta)\\
& = & (L_D\omega) \circ \eta +  \omega\circ L_D(\eta)
\end{array}
$$\qed

From the universal property of $\Omega^1QA$ follows that there exists a 
unique $QA$-bimodule map $\Phi: \Omega^1QA \ra QA$ such that $L_D = \Phi\circ\delta$
where $\delta : QA\ra \Omega^1 QA$ is the canonical  differential. 
If $\tau$ is a supertrace on $QA$ we can define a map $T: \Omega^1QA \ra \bc$ by 
$$
T(\omega\delta \eta) = \tau(\omega\Phi(\delta\eta)) 
$$
\bproposition
$T$ is a supertrace on $\Omega^1QA$
\eproposition
Proof. We need to show that $T$ vanishes on supercommutators. Let $\alpha$, $\omega$ 
and $\eta $ be homogeneous elements of $QA$. Then the supercommutator $[\alpha, 
\omega\delta\eta]_s$ 
with respect to the Fedosov product is given by the following formula
$$
[\alpha, \omega\delta\eta]_s = \alpha\omega\delta\eta - (-1)^{\deg(\alpha)\deg(\omega)\deg(\eta)}
(\omega\delta(\eta\alpha) - \omega\eta\delta\alpha)
$$
so that 
$$
\begin{array}{rcl}
T[\alpha, \omega\delta\eta]_s & = & \tau\left\{\alpha\omega L_D\eta \right. \\
& & \mbox{} \left. - (-1)^{\deg(\alpha)\deg(\omega)\deg(\eta)}
(\omega L_D(\eta\alpha) - \omega\eta L_D\alpha)\right\}\\
&= & \tau(\alpha\omega L_D\eta) - \tau(\alpha\omega L_D\eta)
\\
& = & 0
\end{array}
$$ \qed

\btheorem 
If $D$ is a derivative $D: A\ra A$ on a Banach algebra $A$ then the 
corresponding Lie derivative $L_D$ acts as zero on the entire cyclic cohomology 
of $A$. 
\etheorem
Proof. As we have seen in the proof of the simplicial normalization theorem, 
we may represent entire cohomology classes by normalized entire cocycles
which correspond to continuous supertraces on $Q_\epsilon A$. 
Let $\tau$ be such a trace. Then for any $\eta \in Q_\epsilon A$
$$L_D( \tau ) (\eta) = \tau L_D\eta = T\delta\eta
$$
which shows that 
$$
L_D(\tau) = T\circ \delta
$$
so that $L_D(\tau)$ represents a null-homotopic supertrace on 
$\Omega^1Q_\epsilon A$. \qed

\subsection{The JLO cocycle} 

The entire cyclic cohomology allows one to extend the 
theory of finitely summable Fredholm modules to the cases 
where such assumption is too restrictive. Explicit examples
of such situations are given in Connes' book \cite{Connes:Book}. 
We shall content ourselves here with a brief exposition of 
the so called JLO-cocycle and Chern characters of 
$\theta$-summable modules. 

Let $T$ be a compact operator in a Hilbert space $\bh$, and 
let us denote by $\mu_n(T)$ its $n$-th eigenvalue. 
The following  notion has been introduced by Connes in 
\cite{Connes:Entire}. 

\bdefinition
A $\theta$-summable Fredholm module over an algebra $A$ is 
given by a Fredholm module $(\bh, F)$ over $A$ such that 
for all $a\in A$ $[F,a]\in J^{1/2}$, 
where $J^{1/2}$ is an ideal of compact operators defined by 
$$
J^{1/2} = \{ T\in {\bk}(\bh)\mid \mu_n(T) = O((\log n)^{-1/2})\}
$$
\edefinition

$\theta$-summable Fredholm modules may be described in 
terms of unbounded operators as follows. Connes proves 
that the conditions of the $\theta$-summable Fredholm 
module imply the existence of an unbounded self-adjoint operator $D$ 
satisfying the following properties: 
\begin{enumerate}
\item $\sign D = D/\mid D\mid = F$. 
\item For any $a\in A$, the operator $[D,a]$ is bounded. 
\item $\tr(\exp (-D^2))$ is finite. 
\end{enumerate}

The resulting pair $(\bh , D)$ is called a $\theta$-summable 
$K$-cycle. We want to define a character associated with 
this object. As we recall, the notion of a character, 
which was defined in terms of cyclic cocycles, depended on the
existence of a certain trace. In the present case, the 
character will be defined as an entire cyclic cocycle. 

It follows from  our discussion of the entire cyclic cohomology
that we may represent entire cohomology classes by cocycles 
coming from supertraces on the Cuntz algebra $QA$. Unfortunately, 
in applications, it is not always easy to construct such 
supertraces. To circumvent this problem, Connes has introduced
a certain involution algebra $\call$
 together with a 
continuous algebra homomorphism
$$
\rho: Q_\epsilon A \arr{} \call
$$
It is possible to construct a trace $\tau$ on $\call$, which 
then pulls back to a trace on $Q_\epsilon A$, allowing us 
to construct entire cyclic cocycles. We shall give 
the required formula in the even case; both cases are 
treated in Connes' paper and also in his book. 

\btheorem
Let $(\bh, D, \gamma)$ be an even $\theta$-summable $K$-cycle 
on a locally convex algebra $A$, and let $\tau$ be a trace
on the involution algebra $\call$. Then the following 
equality defines a normalized even entire cocycle 
$\{\phi_{2n}\}_{n\in \bn}$:
$$
\phi_{2n}(a_0, \dots , a_{2n})= 
\Gamma(n+\frac{1}{2})\tau\rho(Fa_0[F,a_1]\cdots [F, a_{2n}])
$$
for all $a_j\in A$. 
\etheorem
 
This is the $\theta$-summable analogue of the Chern character
of a Fredholm module. Unfortunately, the formula given above 
is not very easy to use in applications. Jaffe, Le\'{s}niewski
and Osterwalder have defined a simpler entire cyclic cocycle
\cite{JLO} which was then shown to be cohomologous  to 
the Chern character by Connes in \cite{Connes:Cohomol}. 
The construction of the JLO-cocycle, as it is now known, 
was motivated by $\theta$-summable $K$-cycles arising 
in the supersymmetric quantum field theory. 

Let $\triangle_n$ be the set 
$$
\triangle_n = \{s_i \geq 0, i= 0,  \dots, n\mid \sum s_i = 1\}
$$
Then the even JLO-cocycle is defined by the following formula
$$
\psi_{2n}(a_0, \dots , a_{2n}) = 
\int_{\triangle_n} ds_0\cdots ds_n
\tr(\gamma a_0 e^{-s_0D^2}[D,a_1]\cdots 
[D, a_{2n}]e^{-s_{2n}D^2})
$$
There is an interesting discussion of the algebraic properties 
of this cocycle using superconnections in Quillen's paper
\cite{Q:Cochains}.

\vfill\eject

%% bibiliografia do wykladow w Poznaniu

\end{document}